\documentclass[11pt]{article} 

\usepackage[utf8]{inputenc} 

\usepackage[margin=1in]{geometry} 
\usepackage{graphicx} 

\usepackage{booktabs} 
\usepackage{array} 
\usepackage{paralist} 
\usepackage{verbatim} 
\usepackage{mathrsfs}
\usepackage{amssymb}
\usepackage{amsthm}
\usepackage{amsmath,amsfonts,amssymb,esint}
\usepackage{graphics}
\usepackage{enumerate}
\usepackage{mathtools}
\usepackage{xfrac}
\usepackage{bbm}
\usepackage{subcaption}
\usepackage{times}

\usepackage[usenames,dvipsnames]{xcolor}
\usepackage[colorlinks=true, pdfstartview=FitV, linkcolor=blue, citecolor=blue, urlcolor=blue]{hyperref}

\usepackage{pgfplots}

\numberwithin{equation}{section}

\newtheorem{theorem}{Theorem}[section]
\newtheorem{conjecture}[theorem]{Conjecture}

\newtheorem{proposition}[theorem]{Proposition}
\newtheorem{lemma}[theorem]{Lemma}
\newtheorem{definition}[theorem]{Definition}
\theoremstyle{definition}
\newtheorem{remark}[theorem]{Remark}
\newtheorem{open_problem}{Problem}

\newcommand{\norm}[1]{\left\|#1\right\|}
\newcommand{\abs}[1]{\left|#1\right|}

\newcommand{\average}[1]{\left\langle#1\right\rangle}

\newcommand*{\supp}{\ensuremath{\mathrm{supp\,}}}
\newcommand*{\Id}{\ensuremath{\mathrm{Id}}}
\renewcommand*{\div}{\ensuremath{\mathrm{div\,}}}

\newcommand*{\N}{\ensuremath{\mathbb{N}}}

\newcommand*{\T}{\ensuremath{\mathbb{T}}}
\newcommand*{\Z}{\ensuremath{\mathbb{Z}}}
\newcommand*{\R}{\ensuremath{\mathbb{R}}}
\newcommand*{\CC}{\ensuremath{\mathbb{C}}}
\newcommand*{\tr}{\ensuremath{\mathrm{tr\,}}}
\newcommand{\eps}{\varepsilon}

\newcommand{\RR}{\mathring R}
\newcommand{\Rbar}{\mathring {\overline R}}
\newcommand{\RSZ}{\mathcal R}
\renewcommand*{\Re}{\ensuremath{\mathrm{Re\,}}}
\renewcommand*{\tilde}{\widetilde}
\renewcommand*{\hat}{\widehat}
\newcommand*{\curl}{\ensuremath{\mathrm{curl\,}}}

\newcommand{\Proj}{\ensuremath{\mathbb{P}}}
\newcommand{\les}{\lesssim}

\usepgfplotslibrary{fillbetween}
\pgfmathdeclarefunction{wave}{1}{\pgfmathparse{#1+(1/6)*sin(2*pi*deg(x))}}
\pgfmathdeclarefunction{line}{1}{\pgfmathparse{#1}}

\title{Convex integration and phenomenologies in turbulence}

\author{Tristan Buckmaster
\thanks{Department of Mathematics, Princeton University.
{\footnotesize \href{mailto:buckmaster@math.princeton.edu}{buckmaster@math.princeton.edu}.}
}
\and 
Vlad Vicol
\thanks{Courant Institute for Mathematical Sciences,  New York University.
{\footnotesize \href{mailto:vicol@cims.nyu.edu}{vicol@cims.nyu.edu}.}
}
}
\date{} 

\usepackage[nottoc,notlot,notlof]{tocbibind}

\begin{document}

\maketitle

\begin{abstract}
In this review article we discuss a number of recent results concerning {\em wild} weak solutions of the incompressible Euler and Navier-Stokes equations. These results build on the groundbreaking works of De Lellis and Sz\'ekelyhidi Jr., who extended Nash's fundamental ideas on $C^1$ flexible isometric embeddings, into the realm of fluid dynamics. These techniques, which go under the umbrella name {\em convex integration}, have fundamental analogies with the phenomenological theories of hydrodynamic turbulence~\cite{DeLellisSzekelihidi12,Szekelyhidi12,DLSZ17,DLSZ19}. Mathematical problems arising in turbulence (such as the Onsager conjecture) have not only sparked new interest in convex integration, but certain experimentally observed features of turbulent flows (such as intermittency) have also informed new convex integration constructions. 

First, we give an elementary construction of nonconservative $C^{0+}_{x,t}$ weak solutions of the Euler equations, first proven by De~Lellis-Sz\'ekelyhidi Jr.~\cite{DeLellisSzekelyhidi12a,DLSZ13}. Second, we  present Isett's~\cite{Isett16} recent resolution of the flexible side of the Onsager conjecture. Here, we in fact follow the joint work~\cite{BDLSV17} of De~Lellis-Sz\'ekelyhidi Jr.~and the authors of this paper, in which weak solutions of the Euler equations in the regularity class $C^{\sfrac 13-}_{x,t}$ are constructed, attaining any energy profile. Third, we give a concise proof of the authors' recent result~\cite{BV}, which proves the existence of infinitely many weak solutions of the Navier-Stokes in the regularity class $C^0_t L^{2+}_x \cap C^0_t W^{1,1+}_x$. We conclude the article by mentioning a number of open problems at the intersection of convex integration and hydrodynamic turbulence.
\end{abstract}

\setcounter{tocdepth}{1}
\tableofcontents

\section{Introduction}
The \emph{Navier-Stokes} equations, written down almost 200 years ago \cite{Navier1823,Stokes1845}, are thought to be the fundamental set of equations governing the motion of  viscous fluid flow. In their homogenous incompressible form these equations predict the evolution of the velocity field $v$ and scalar pressure $p$ of the fluid by
\begin{subequations}
\label{eq:NSE}
\begin{align}
\partial_t v + \div (v \otimes v) + \nabla p - \nu \Delta v &=0 \, 
\label{eq:NSE:1}\\
\div v &= 0 \,.
\label{eq:NSE:2}
\end{align}
\end{subequations}
Here $\nu>0$ is the kinematic viscosity. One may rewrite the nonlinear term in non-divergence form as $\div (v \otimes v) = (v \cdot \nabla) v$. The Navier-Stokes equations may be derived rigorously from the Boltzmann equation~\cite{BardosGolseLevermore93,LionsMasmoudi01b,GolseStRaymond04}, or from lattice gas models~\cite{QuastelYau98}. In three dimensions, the global in time well-posedness for \eqref{eq:NSE} remains famously unresolved and is the subject of one of the \emph{Millennium Prize} problems~\cite{F2006}. 

Formally passing to the {\em inviscid limit $\nu \to 0$} we arrive at the {\em Euler} equations,  
which are the classical model for the motion of an incompressible homogenous inviscid fluid, and were in fact derived a century earlier
by Euler~\cite{Euler1757}. The equations for the unknowns $v$ and $p$ are
\begin{subequations}
\label{eq:Euler}
\begin{align}
\partial_t v + \div (v \otimes v) + \nabla p  &=0 \, ,
\label{eq:Euler:1}\\
\div v &= 0 \, .
\label{eq:Euler:2}
\end{align}
\end{subequations}
As for their viscous counterpart, the global in time well-posedness for the three-dimensional Euler equations remains an outstanding open problem, arguably of greater physical significance~\cite{Constantin07}. Indeed, an Euler singularity requires infinite velocity gradients and is thus intimately related to the anomalous dissipation of energy for turbulent flows~\cite{Frisch95}. 

When considering the Cauchy problem, the Navier-Stokes and Euler equations are to be supplemented with an incompressible initial datum $v_0$. For simplicity, throughout this paper the systems \eqref{eq:NSE} and \eqref{eq:Euler} are posed on the periodic box $\T^3 = [-\pi,\pi]^3$, and the initial condition $v_0$ is assumed to have zero mean on $\T^3$. Since solutions $v(\cdot,t)$ preserve their mean, we have $\int_{\T^3} v(x,t) dx = \int_{\T^3} v_0(x) dx = 0$ for all $t>0$. The pressure is uniquely defined under the normalization $\int_{\T^3} p(x,t) dx = 0$ by solving $-\Delta p = \div \div (v\otimes v)$. In order to ensure a nontrivial long-time behavior, it is customary to add a zero mean forcing term $f(x,t)$ to the the right side of the Navier-Stokes equations \eqref{eq:NSE:1}. 

In the bulk of this paper (cf.~Sections~\ref{sec:Euler:C0}--\ref{sec:NSE:L2}) we consider {\em weak, or distributional solutions} of the Navier-Stokes and Euler equations (defined precisely in Section~\ref{sec:math}). The motivation for considering weak solutions of~\eqref{eq:NSE} and~\eqref{eq:Euler} stems from the Kolmogorov and Onsager theories of {\em hydrodynamic turbulence}~\cite{Frisch95}. The fundamental ansatz of these phenomenological theories is that in the vanishing viscosity limit solutions of the Navier-Stokes equations {\em do not remain smooth uniformly with respect to $\nu$} (in a sense to be made precise in Section~\ref{sec:physics}), and thus may only converge to distributional solutions of the Euler equations.  Therefore, in an attempt to translate predictions made by turbulence theory into mathematically rigorous questions, it is natural to work within the framework of weak solutions. 

\subsection*{Organization of the paper}
\begin{itemize}

\item[Section~\ref{sec:physics}:] 
In this section we discuss some of the fundamental features of the Kolmogorov~\cite{Kolmogorov41a,Kolmogorov41b,Kolmogorov41c} and Onsager~\cite{Onsager49} phenomenological theories of fluid turbulence. This topic is too vast to review here in detail. For a detailed study of turbulence theories we refer the reader to the books~\cite{Batchelor53,MoninYaglom71,MoninYaglom13,Frisch95,FoiasManleyRosaTemam01,Tsinober01},  the surveys papers with a mathematical perspective~\cite{Constantin94,Eyink95,Robert03,Constantin06,EyinkSreeniviasan06,Shvydkoy10,BardosTiti13,Eyink18}, and  to the references therein.  In Section~\ref{sec:physics} we focus on the 
features that relate most to the convex integration constructions considered in later sections: the anomalous dissipation of energy, energy fluxes, scalings of structure functions, and  intermittency. 

\item[Section~\ref{sec:math}:] 
The phenomena modeled by the Euler and Navier-Stokes equations are not just important, but also fascinating; see e.g.~the images in van Dyke's ~{\em Album of fluid motion}~\cite{VanDyke82}. Consequently, the literature concerning the analysis of these equations is vast. For an overview of the field we refer the reader to the classical texts~\cite{ConstantinFoias88,Temam95} and also to the books~\cite{DoeringGibbon95,Lions96,Chemin98,Temam01,FoiasManleyRosaTemam01,MB2002,LemarieRieusset02,MP12,Lemarie16,RSR16,Tsai18}. 
In order to place the convex integration constructions in context, in this section we recall only a few of the rigorous mathematical results known about \eqref{eq:NSE} and \eqref{eq:Euler}. We focus on the definition(s) and regularity of weak solutions, we discuss the results which have led to the resolution of the Onsager conjecture, and present the recent results concerning the non-uniqueneess of weak solutions for the Navier-Stokes equations and related hydrodynamic models.

\item[Section~\ref{sec:convex:integration}:]
In this section we summarize some of the key aspects of Nash-style convex integration schemes in fluid dynamics.  A number of excellent surveys articles on this topic are already available in the literature, by De Lellis and Sz\'ekelyhidi Jr.~\cite{DeLellisSzekelihidi12,Szekelyhidi12,DLSZ17,DLSZ19}. These surveys discuss in detail the Nash-Kuiper theorem~\cite{Nash54,Kuiper55}, Gromov's $h$-principle~\cite{Gromov73}, convex integration constructions for flexible differential inclusions inspired by the work of M\"uller and \v{S}ver\'ak~\cite{MullerSverak03}, the Scheffer~\cite{Scheffer93}-Shnirelman~\cite{Shnirelman00} constructions, leading to the constructions of non-conservative H\"older continuous solutions of the Euler equations. Our goal here is to discuss some of the intuition behind Nash-style convex integration schemes for the Euler equations, and to provide the mathematical intuition behind the intermittent building blocks the authors have introduced~\cite{BV} for the Navier-Stokes equations. 

These heuristic arguments are made precise in Sections~\ref{sec:Euler:C0},~\ref{sec:Euler:C1/3}, and~\ref{sec:NSE:L2}, below where we give the proofs for some of the more recent developments in the field. Our aim in these subsequent sections is to present concise proofs, rather than the most general results. 

\item[Section~\ref{sec:Euler:C0}:]
We present the main result of De Lellis and Sz\'ekelyhidi Jr.'s paper~\cite{DLSZ13}, cf.~Theorem~\ref{thm:Euler:C0:DLSZ} below. This work gave the first proof for the existence of a $C^{0+}_{x,t}$ weak solution of the 3D Euler equations which is non-conservative, following a Nash scheme with Beltrami building blocks. To simplify the presentation we only show the existence of a non-conservative weak solution in this regularity class, cf.~Theorem~\ref{thm:Euler:C0} below.

\item[Section~\ref{sec:Euler:C1/3}:]
The $C^{0+}_{x,t}$ construction discussed in the previous section may be viewed as the start in the race towards proving the flexible side of the Onsager conjecture. In this section we present the resolution of this conjecture on the H\"older scale, recently established by Isett~\cite{Isett16}, cf.~Theorem~\ref{thm:Onsager:critical:Isett} below. We discuss the papers on which Isett's construction relies (e.g.~the Mikado flows by Daneri-Sz\'ekelyhidi Jr.~\cite{DSZ17}) and the main ideas in Isett's work~\cite{Isett16}. The proof we present in this section is that of Theorem~\ref{thm:Onsager:critical:dissipative}, established by De Lellis, Sz\'ekelyhidi Jr.~and the authors of this paper in~\cite{BDLSV17}. This work extends and simplifies~\cite{Isett16}, allowing one to construct {\em dissipative weak solutions in the regularity class $C^{\sfrac 13-}_{x,t}$}. The exposition follows~\cite{BDLSV17} closely, but the presentation of Mikado flows is slightly different, as to be consistent with the intermittent jets which we introduce in the next section. 

\item[Section~\ref{sec:NSE:L2}:] 
We discuss the main ideas of the authors' recent result~\cite{BV}, cf.~Theorem~\ref{thm:BV:main} below. To simplify the presentation we only give the proof of Theorem~\ref{thm:NSE:main} which establishes the existence of a weak solution to the Navier-Stokes equations in the regularity class $C^0_t( L^{2+}_x \cap W^{1,1+}_x)$, with kinetic energy that is {\em not monotone decreasing}. This result directly implies the non-uniqueness of weak solutions: compare the solution of  Theorem~\ref{thm:NSE:main} to the Leray solution with the same initial condition; the later's energy is non-increasing, hence they cannot be the same. The main idea in the proof is to use intermittent building blocks in the convex integration construction, such as {\em intermittent Beltrami flows}~\cite{BV} or {\em intermittent jets}~\cite{BCV18}.
The development of intermittent building blocks was for the first time announced in the context of the authors' unpublished work with Masmoudi.

\item[Section~\ref{sec:open}:]
We conclude the paper with a number of open problems. Most of these problems are well-known in the field and concern the regularity of weak solutions for the Euler and Navier-Stokes equations. We additionally discuss open problems regarding convex integration constructions in fluid dynamics.
\end{itemize}

\subsection*{Acknowledgments} 
The work of T.B. has been partially supported by the NSF grant DMS-1600868.
V.V. was partially supported by the NSF grant DMS-1652134. 
The authors are grateful to Raj Beekie, Theodore Drivas, Matthew Novack, and Lenya Ryzhik for suggestions and stimulating discussions concerning aspects of this review.

\section{Physical motivation}
\label{sec:physics}

Hydrodynamic turbulence remains the greatest challenge at the intersection of  mathematics and physics. During the past century our understanding of this phenomenon was greatly enriched by the predictions of Prandtl, von Karman, Richardson, Taylor, Heisenberg, Kolmogorov, Onsager, Kraichnan, etc. The success of their theories in modeling the statistics of turbulent flows has been astounding~\cite{Frisch95}. Nevertheless, to date no single mathematically rigorous (unconditional) bridge between the incompressible Navier-Stokes equations at high Reynolds number and these phenomenological theories  has been established (cf.~Section~\ref{sec:0th:law}). 

To fix the notation in this section, let us denote by $v^\nu = NSE^\nu(v_0,t)$ a solution of the Cauchy problem for the {\em forced} version of the Navier-Stokes equations \eqref{eq:NSE} with viscosity\footnote{Throughout this paper we abuse notation and denote also by $\nu$ the inverse of the Reynolds number $\Re^{-1} = \nu /(UL)$, where $L = 2\pi $ is the characteristic length scale of the domain $\T^3$, and $U$ is an average r.m.s.\ velocity, e.g.~$U = (\fint_{\T^3} |v_0(x)|^2 dx )^{\sfrac 12}$. The infinite Reynolds number limit $\Re \to \infty$ is used interchangeably with the vanishing viscosity limit $\nu \to 0$, keeping $U$ and $L$ fixed.}  $\nu$, at time $t$, and with initial datum $v_0 \in L^2$ which is taken to be incompressible, zero mean, and sufficiently smooth. The 
forcing\footnote{While in this paper we restrict ourselves to the deterministic framework, in turbulence theory one typically considers a stochastic forcing term~\cite{Novikov65,BensoussanTemam73,VishikKomechFusikov79, Eyink96}: a wave-number localized, gaussian and white in time forcing as a source which drives turbulent cascades. In this setting, one may sometimes rigorously establish the existence, uniqueness, and mixing properties of invariant measures for the underlying Markov semigroup, e.g.~\cite{HairerMattingly06, HairerMattingly2011,KuksinShirikian12} for two dimensional flows and~\cite{BarbatoFlandoliMorandin2012,FGHV16} for dyadic shell models (see also~\cite{Debussche13,GlattHoltz14}). These invariant states are expected to encode the statistics of turbulent flows at high Reynolds number~\cite{Frisch95}.} term $f^\nu$ is taken to have zero mean, is (statistically) stationary, and injects energy into the system at {\em low frequencies}.\footnote{To make this precise, one may for instance assume that there exists an inverse length scale $\kappa_I$, independent of $\nu$, such that $\Proj_{\leq \kappa} f^\nu = \Proj_{\leq \kappa_I} f^\nu$ for all $\kappa \geq \kappa_I$. Here and throughout the paper $\Proj_{\leq \kappa}$ denotes a Fourier multiplier operator, which projects a function onto its Fourier frequencies  $\leq \kappa$ in absolute value. Equivalently, $\Proj_{\leq \kappa}$ is a mollification operator at length scales $\leq \kappa^{-1}$.} 
The equations are posed on $\T^3$ with periodic boundary conditions.\footnote{Here we leave aside the physically extremely important, but mathematically very challenging issue of fluid motion in bounded domains~\cite{Schlichting60,MaekawaMazzucato18}. In laboratory experiments the generation of a turbulent flow involves the presence of a solid boundary, such as flat plate or a grid mesh. Classically, the Navier-Stokes system \eqref{eq:NSE} is supplemented with no-slip Dirichlet boundary conditions for the velocity field at the solid wall, whereas for the Euler system \eqref{eq:Euler} the non-penetrating boundary conditions are imposed. The vanishing viscosity limit $\nu\to 0$ leads to the consideration of boundary layers which typically separate from the wall; one of the fundamental driving mechanisms for the transition from a laminar to a turbulent regime~\cite{Kato84b,DrazinReid04}. See also the discussion in~\cite{BardosTiti13,ConstantinVicol18}.}

Given the complex nature of turbulent flows, it is unreasonable to expect to make predictions about individual solutions $v^\nu$ to the Cauchy problem for the forced system~\eqref{eq:NSE}. Indeed, theories of fully developed turbulence typically  attempt to make {\em statistical predictions} about the behavior of fluid flow at high Reynolds numbers, away from solid boundaries, for length scales in the inertial range, and under certain assumptions -- for instance, ergodicity, statistical homogeneity, isotropy, and self-similarity. Note that typically it is not possible to rigorously prove these assumptions  directly from first principles (e.g.\ from the Navier-Stokes equations), and so certain ambiguities arise.  One of these ambiguities lies in the definition of a {\em statistical average}, denoted in this section by $\average{\cdot}$. 

Given a suitable {\em observable $F$} of the solution $v^\nu$, theoretical physics considerations typically use $\average{F(v^\nu)}$ to denote an {\em ensemble average} with respect to a putative probability measure $\mu^\nu$ on $L^2$ which is time independent.\footnote{Following the pioneering work of Foias~\cite{Foias72,Foias73}, in a purely deterministic setting one may consider the concept of a {\em stationary statistical solution} to the Navier-Stokes equation. Stationary statistical solutions are probability measures on $L^2$ which satisfy a stationary Liouville-type equation, integrated against cylindrical test functions. Their existence may be rigorously established using the concept of a generalized Banach limit from long time averages, but their uniqueness remains famously open. This notion of solution has been explored quite a bit in the past decades, see e.g. the books~\cite{VishikKomechFusikov79,FoiasManleyRosaTemam01}.} That is, ones assumes to be at statistical equilibrium, and that the probability measure $\mu^\nu$ encodes the macroscopic statistics of the flow. On the other hand, in laboratory experiments a measurement of the turbulent flow is usually a {\em long time average} at fixed viscosity, in order to reach a stationary regime. That is, one observes solutions which are close to, or on, the attractor of the system.\footnote{It is typical in certain laboratory measurements to recast measurements from the time domain into the space domain by appealing to the Taylor hypothesis~\cite{Frisch95}.}  In analogy with classical statistical mechanics, turbulence theories deal with the possible discrepancy between ensemble averages and statistical averages  by making an impromptu {\em ergodic hypothesis}. The implication of the ergodic hypothesis is that averages against an ergodic invariant measure (possibly also mixing), are the same as long time averages, giving a meaning to $\average{\cdot}$. Lastly, we note that in this section  $\average{\cdot}$ sometimes includes a spatial average over $\T^3$, which may be justified under the assumption of statistical homogeneity.

\subsection{Anomalous dissipation of energy}
\label{sec:0th:law}
The fundamental ansatz of Kolmogorov's 1941 theory of fully developed turbulence~\cite{Kolmogorov41a,Kolmogorov41b,Kolmogorov41c}, sometimes called the {\em zeroth law of turbulence}, postulates the {\em anomalous dissipation of energy} -- the non-vanishing of the rate of dissipation of kinetic energy of turbulent fluctuations per unit mass, in the limit of zero viscosity (cf.~\eqref{eq:anomaly} below). The zeroth law of turbulence is verified experimentally to a tremendous degree~\cite{Sreenivasan98,PearsonEtAl02,KanedaEtAl03}, but to date we do not have a single example where it is rigorously proven to hold, directly from \eqref{eq:NSE}.

To formulate this ansatz, we use the aforementioned notation and denote by  $v^\nu$ be a solution of \eqref{eq:NSE} with stationary force $f^\nu$. We start with the balance of kinetic energy in the Navier-Stokes equation, in order to derive a correct formula for the energy dissipation rate per unit mass. By taking an inner product of $v^\nu$ with the forced \eqref{eq:NSE} system, and {\em assuming the functions $v^\nu$ are sufficiently smooth}, one arrives at the pointwise energy balance
\begin{align}
\partial_t \frac{|v^\nu|^2}{2} + \nabla \cdot \left( v^\nu \left(\frac{|v^\nu|^2}{2} + p^\nu \right) - \nu \nabla \frac{|v^\nu|^2}{2}\right) = f^\nu \cdot v^\nu -  \nu |\nabla v^\nu|^2.
\label{eq:NSE:formal:energy:balance}
\end{align}
Integrating over the periodic domain we obtain the kinetic energy balance
\begin{align}
\frac{d}{dt} \fint_{\T^3} \frac{|v^\nu|^2}{2} dx = \fint_{\T^3} f^\nu \cdot v^\nu dx - \nu \fint_{\T^3} |\nabla v^\nu|^2 dx \, ,
\label{eq:NSE:global:energy:balance}
\end{align}
where the first term on the right side  denotes the total work of the force and the second term denotes the energy dissipation rate per unit mass. Note that all the terms in \eqref{eq:NSE:global:energy:balance} have dimensional units of $U^3 L^{-1}$.
Estimate \eqref{eq:NSE:global:energy:balance} is the only known coercive a-priori estimate for the 3D Navier-Stokes equations, and it gives an a-priori bound for the solution $v^\nu$ in the so-called  {\em energy space} $L^\infty_t L^2_x \cap L^2_t H^1_x$. Leray~\cite{Leray34} used the energy balance for a suitable approximating sequence,  combined with a compactness argument,  to prove the existence of a  global in time weak solution  to \eqref{eq:NSE} which lies in $L^\infty_t L^2_x \cap L^2_t H^1_x$, and obeys \eqref{eq:NSE:global:energy:balance} weakly in time with an inequality  instead of the equality. See Definition~\ref{d:leray-hopf} below for the definition of a Leray solution for \eqref{eq:NSE}.
A-posteriori one may ask the question of whether the local energy balance \eqref{eq:NSE:formal:energy:balance} may be actually justified when $v^\nu \in L^\infty_t L^2_x \cap L^2_t H^1_x$ is a weak solution of equation \eqref{eq:NSE}. To date this question remains open (see however the works~\cite{Lions60,Shinbrot74,Kukavica06,CCFS08,LeslieShvydkoy18} for sufficient conditions). Instead, following the work of Duchon-Robert~\cite{DuchonRobert00} -- equivalently, the commutator formula of Constatin-E-Titi~\cite{ConstantinETiti94} -- one may prove that for a weak solution $v^\nu$ in the energy class (by interpolation   $v^\nu$ also lies in $L^{\sfrac{10}{3}}_{x,t} \supset L^3_{x,t}$), the following equality holds in the sense of distributions
\begin{align}
\partial_t \frac{|v^\nu|^2}{2} + \nabla \cdot \left( v^\nu \left(\frac{|v^\nu|^2}{2} + p^\nu \right) - \nu \nabla \frac{|v^\nu|^2}{2}\right) = f^\nu \cdot v^\nu -  \nu |\nabla v^\nu|^2 - D(v^\nu)
\label{eq:NSE:local:energy:balance}
\end{align}
where the $(x,t)$-distribution $D(v^\nu)$ is defined by a weak form of the K\'arm\'an-Howarth-Monin relation~\cite{KarmanHowarth38,Monin59} (see also~\cite{MoninYaglom71,Frisch95})
\begin{align}
D(v^\nu)(x,t) = \lim_{\ell \to 0}  \frac 14 \int_{\T^3} \nabla \varphi_\ell(z) \cdot \delta v^\nu(x,t;z) |\delta v^\nu(x,t;z)|^2 dz \, .
\label{eq:KHM}
\end{align}
In \eqref{eq:KHM} we have denoted the {\em velocity increment in the direction $z$} by 
\begin{align}
\delta v^\nu(x,t;z) = v^\nu(x+z,t) - v^\nu(x,t)
\label{eq:increment}
\end{align}
and the approximation of the identity $\varphi_\ell$ is given by $\varphi_\ell(z) = \frac{1}{\ell^3} \varphi\left(\frac{z}{\ell}\right)$, where $\varphi \geq 0$ is an even bump function with mass equal to $1$. The limit in \eqref{eq:KHM} is a limit of $L^1_{x,t}$ objects in the sense of distributions, and it is shown in~\cite{DuchonRobert00} that $D(v^\nu)$ is independent of the choice of $\varphi$. When compared to \eqref{eq:NSE:formal:energy:balance}, identity~\eqref{eq:NSE:local:energy:balance} additionally takes into account the possible dissipation  of kinetic energy, due to possible singularities of the flow $v^\nu$, encoded in the {\em defect measure} $D^\nu$. Note that if $v^\nu$ is sufficiently smooth to ensure that $\lim_{|z| \to 0} \frac{1}{|z|} \int_0^T \int_{\T^3} | \delta v^\nu(x,t;z)|^3 dx dt = 0$, then one may directly show that $D(v^\nu) \equiv 0$ (cf.~\cite{CCFS08}). See Sections~\ref{sec:turbulence:Onsager} and~\ref{ss:euler_results} below for more details. Similarly to \eqref{eq:NSE:global:energy:balance}, once we average the local energy balance \eqref{eq:NSE:local:energy:balance} over $\T^3$, the divergence term on the left side vanishes, and we are left with 
\begin{align}
\frac{d}{dt} \fint_{\T^3} \frac{|v^\nu|^2}{2} dx = \fint_{\T^3} f^\nu \cdot v^\nu dx - \nu \fint_{\T^3} |\nabla v^\nu|^2 dx  - \fint_{\T^3} D(v^\nu) dx \, ,
\label{eq:NSE:global:energy:balance:2}
\end{align}
which yields a balance relation between energy input and energy dissipation. 

With \eqref{eq:NSE:global:energy:balance:2}, we define the {\em mean energy dissipation rate per unit mass} by
\begin{align}
\eps^\nu = \nu \average{|\nabla v^\nu|^2} + \average{D(v^\nu)} \label{eq:epsilon:nu} \, ,
\end{align}
where as discussed before, $\average{\cdot}$ denotes a suitable ensemble/long-time and a space average. The quantity $\eps^\nu$ has physical units of $U^3L^{-1}$. The zeroth law of turbulence, or the anomalous dissipation of energy, postulates that in the inviscid  limit $\nu\to 0$ (keeping $U$ and $L$ fixed) the mean energy dissipation rate per unit mass does not vanish, and moreover that there exists an $\eps  \in (0,\infty)$ such that
\begin{align}
\eps = \liminf_{\nu \to 0} \eps^\nu  > 0.
\label{eq:anomaly}
\end{align}
Figures~\ref{fig:Sreeni:1} and~\ref{fig:Pearson} below present classical experimental evidence which is consistent with the positivity of $\eps$.  Further compelling experimental support for \eqref{eq:anomaly} is provided by the more recent numerical simulation~\cite{KanedaEtAl03}, see also the recent review of experimental and numerical evidence~\cite{Vassilicos15}.\footnote{It is worth emphasizing that $\eps>0$ implies that the sequence of Navier-Stokes solutions $\{v^\nu\}_{\nu > 0}$, say of Leray-Hopf kind, cannot remain uniformly bounded (with respect to $\nu$) in the space $L^3_t B^{s}_{3,\infty,x}$ for any $s>\sfrac 13$.  In fact, in~\cite{DrivasEyink17} it is shown that even if $\eps = 0$, but the rate of vanishing of $\eps^\nu$ is slow, say $\liminf_{\nu\to0} \frac{\log (\eps^\nu) }{\log(\nu)} = \alpha \in (0,1]$, then the sequence of Leray solutions $v^\nu$ cannot remain uniformly in the space $L^3_t B^{s}_{3,\infty,x}$ with $s > \frac{1+\alpha}{3-\alpha}$. Thus, experimental evidence robustly points towards Euler singularities.} 
\begin{figure}[h!]
\centering
 \begin{subfigure}[b]{0.55\textwidth}
    \includegraphics[width=0.9\textwidth]{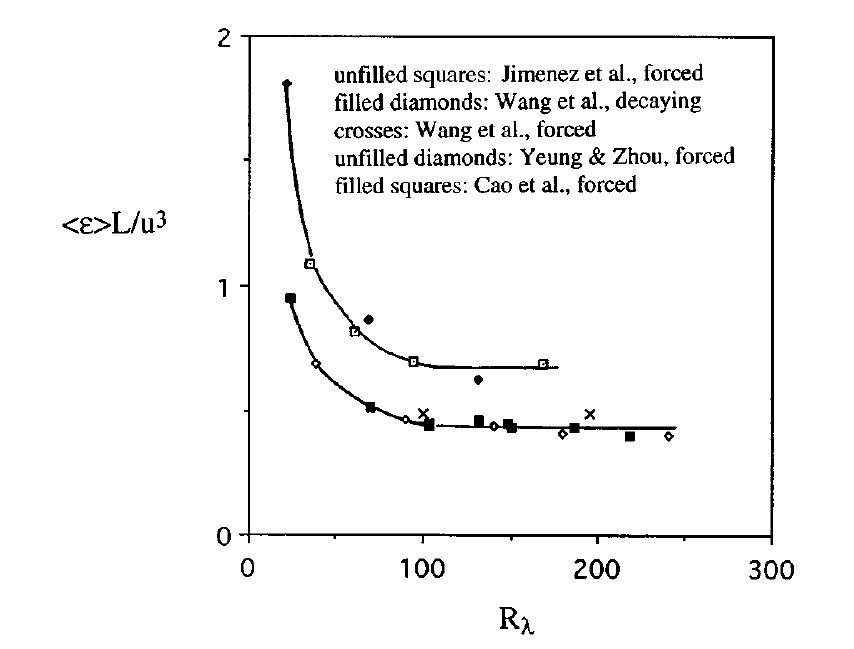}
    \caption{{\small [Sreenivasan~\cite{Sreenivasan98}]:  The variation of the quantity $  \eps^\nu L U^{-3}$ with the Taylor microscale Reynolds number $\Re_\lambda$ in simulations of homogeneous and isotropic turbulence in periodic box.}}
        \label{fig:Sreeni:1}
  \end{subfigure}
  \begin{subfigure}[b]{0.1\textwidth}
  \, 
  \end{subfigure}
  \begin{subfigure}[b]{0.35\textwidth}
     \includegraphics[width=0.9\textwidth]{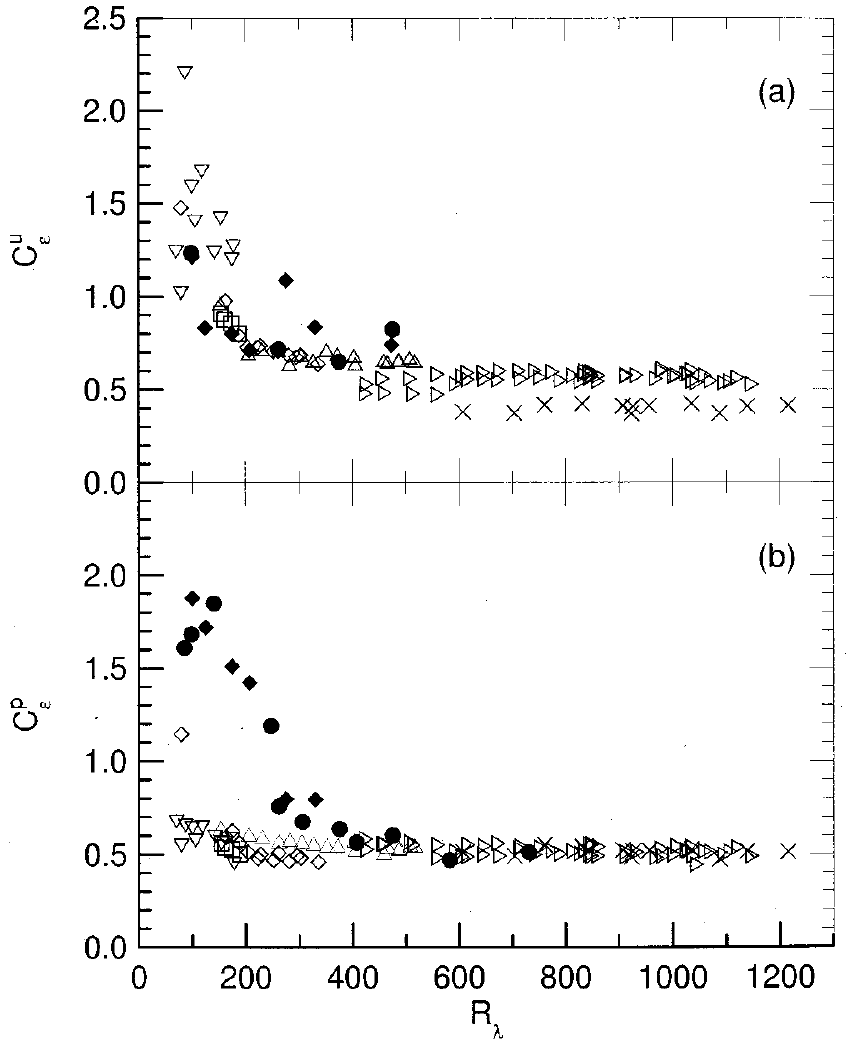}
    \caption{{\small [Pearson, Krogstad, de Water~\cite{PearsonEtAl02}]:  $C^u_\eps =\eps^\nu  L_u U^{-3}$, with $L_u =$ streamwise integral length scale; and $C^p_\eps =  \eps^\nu L_p U^{-3}$, with $L_p =$ inverse wavenumber at which peak in energy spectrum occurs.}}
     \label{fig:Pearson}
  \end{subfigure}
  \caption{{\small Experimental and numerical evidence for the anomalous dissipation of energy.}}
\end{figure}

\subsection{Basics of the Kolmogorov ('41) theory}
\label{sec:turbulence}

Based on  the anomalous dissipation of energy and certain scaling arguments, in 1941 Kolmogorov~\cite{Kolmogorov41a,Kolmogorov41b,Kolmogorov41c} proposed a theory for homogenous isotropic turbulence, whose key predictions we summarize below. We follow the presentation  in~\cite{Frisch95,MoninYaglom13,EyinkSreeniviasan06,Shvydkoy10}, to which we refer the reader for further details.

Besides the zeroth law of turbulence \eqref{eq:anomaly}, the assumptions of Kolmogorov's theory are {\em homogeneity}, {\em isotropy}, and  {\em self-similarity}. Let $\hat z \in {\mathbb S}^2$ be a unit direction vector and let $\ell > 0$ be a length scale in the {\em inertial range}, meaning that $\ell_D \ll \ell \ll \ell_I$, where $\ell_I$ is the integral scale of the system (the inverse of the maximal Fourier frequency $\kappa_I$ active in the force), and $\ell_D = \nu^{\sfrac 34} \eps^{-\sfrac 14}$ is the  {\em Kolmogorov dissipative length scale} (the only object  which has the physical unit  of $L$ and may be written as $\nu^a \eps^b$; recall that $\nu$ has units of $L U$). Recall the notation~\eqref{eq:increment} for velocity increments. Homogeneity is the assumption that the statistics of turbulent flows is shift invariant: at large Reynolds numbers the velocity increment $\delta v^\nu(x,t;\ell \hat z)$ has the same probability distribution for every $x \in \T^3$. Isotropy is the assumption that the statistics of turbulent flows is locally rotationally invariant: the probability distribution for $\delta v^\nu(x,t;\ell \hat z)$ is the same for all $\hat z \in {\mathbb S}^2$. Lastly, self-similarity postulates the existence of an exponent $h>0$, such that  $\delta v^\nu(x,t;\lambda \ell \hat z)$ and  $\lambda^h \delta v^\nu(x,t;\ell \hat z)$ have the same law, for  $\lambda>0$ such that both $\ell$ and $\lambda \ell$ lie in the inertial range.
Based on these assumptions,\footnote{The assumptions listed here are not {\em minimal}, in the sense that one can deduce a number of the predictions of the Kolmogorov theory by assuming less. We refer the reader for instance to~\cite{ConstantinFefferman94,Frisch95,NieTanveer99,ConstantinNieTanveer99,DuchonRobert00,Eyink02,CS2014,Drivas18} in the deterministic setting, and to~\cite{FGHR08,FGHV16,BCZPSW18} in the stochastic one.} the theory makes predictions about structure functions  and the energy spectrum. 

\begin{figure}[h!]
\begin{centering}
    \includegraphics[width=0.5\textwidth]{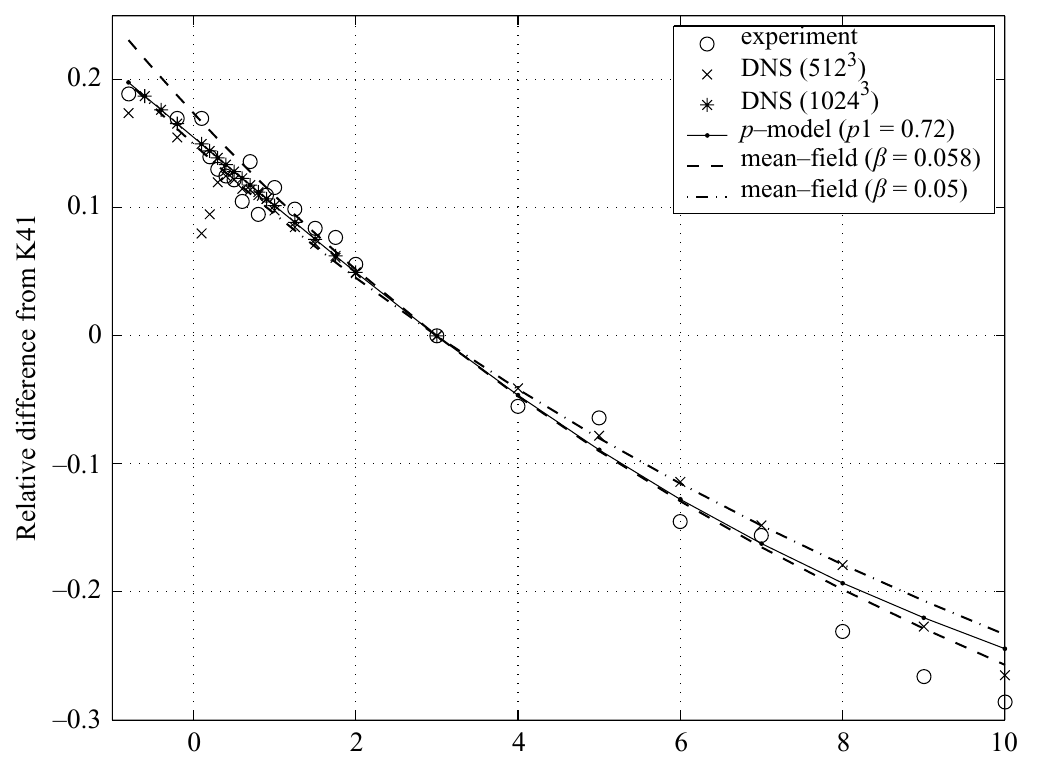}
    \caption{ {\small Mild deviations of $\zeta_p$ from $p/3$. [Chen, Dhruva, Kurien, Sreenivasan, Taylor~\cite{ChenEtAl05}]. The graphic plots $p$ versus $(\zeta_{|p|}-\frac{p}{3}) \frac{3}{p}$ for various nonzero moments $p$, as gathered from experiments (circles) and direct numerical simulations (crosses and stars). The solid line is the outcome of the $p$-model of [Meneveau-Sreenivasan~\cite{MeneveauSreenivasan87}]. The dashed and dotted lines are predictions of the $\beta$ model of [Yakhot~\cite{Yakhot01}].}}
    \label{fig:Chen}
\end{centering} 
\end{figure}

For $p\geq 1$ one may define the $p^{th}$ order {\em longitudinal structure function}
\begin{align*}
S_p^{\|}(\ell) = \average{(\delta v^\nu(x,t;\ell \hat z) \cdot \hat z)^p} 
\end{align*}
where the ensemble/long-time $\average{\cdot}$ takes into account homogeneity and isotropy, so that we do not have to explicitly write averages in $x$ over $\T^3$ and in $\hat z$ over ${\mathbb S}^2$. Note that for $p$ which is odd, $S_p^{\|}(\ell)$ need not a-priori have a sign. Instead, one may define the $p^{th}$ order {\em absolute structure function}  
\begin{align*}
S_p(\ell) = \average{|\delta v^\nu(x,t;\ell \hat z)|^p} 
\end{align*}
which is intimately related to the definition of a Besov space.\footnote{Recall that $v \in B^{s}_{p,\infty}$ means that $v \in L^p$ and that $\sup_{|z |>0} \frac{1}{|z|^s} \norm{\delta v(x;z)}_{L^p_x} < \infty$, for $s\geq 0$ and $p<\infty$.}
$S_p(\ell)$ scales in the same way as $S_p^{\|}(\ell)$, and they both have physical units of $U^p$. Notice that since $\eps \ell$ has units of $U^3$, it follows that $(\eps\ell)^{\sfrac p3}$ has the same physical units as $S_p(\ell)$. Consequently, the only the value of the self-similarity exponent which is consistent with physical units as  $\ell \to 0$ is $h=\sfrac 13$, and thus the Kolmogorov theory predicts the asymptotic behavior
\begin{align}
S_p(\ell) \sim (\eps \ell)^{\sfrac p3}\label{eq:Sp:Kolmogorov}
\end{align}
for $\ell$ in the inertial range, in the infinite Reynolds number limit. Denoting by $\zeta_p$ the limiting {\em structure function exponent}
\begin{align}
\zeta_p = \lim_{\ell \to 0} \lim_{\nu \to 0} \frac{\log \left( S_p(\ell)  \right)}{\log \left(\eps \ell \right)} \, ,
\label{eq:zeta:p:def}
\end{align}
the relation \eqref{eq:Sp:Kolmogorov} indicates that in Kolmogorov's theory of homogenous isotropic turbulence we have \begin{align} 
\zeta_p = \frac{p}{3} \, , \qquad \mbox{for all} \qquad  p\geq 1 \, ,
\label{eq:K:zeta:p}
\end{align} 
in view of the assumption of self-similarity of the statistics at small scales. 
Figure~\ref{fig:Chen} shows that this heuristic argument for the value of $\zeta_p$ yields a surprisingly small deviations, at least for $p$ close to $3$. 
Nonetheless, as seen in Figure~\ref{fig:Chen}, except for $p=3$, when the Kolmogorov prediction {\em $\zeta_3=1$ is indeed supported by all the experimental evidence}, for $p\neq 3$ experiments do indeed deviate from the Kolmogorov prediction. This is related to the phenomenon of intermittency discussed in Section~\ref{sec:intermittent} below.

\begin{figure}[h!]
  \begin{centering}
    \includegraphics[width=0.45\textwidth]{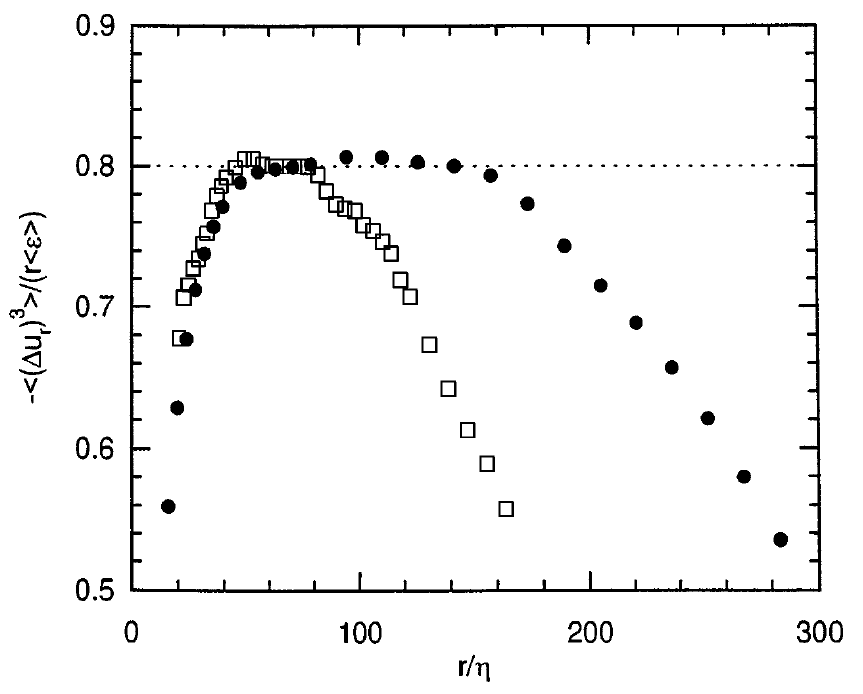}
    \caption{{\small Experimental evidence for $S_{3}^{\|}(\ell) = - \frac{4}{5} (\eps \ell)$. The quantity $- S_3^{\|}(\ell)/(\eps \ell )$ is plotted as a function of $\ell$. Squares denote experimental observations of centerline in pipe flow at $\Re = 230000$. Circles indicate data from a  $512^3$ DNS of homogenous turbulence at $\Re = 220$. Dots indicate the $4/5$ law. [K. R. Sreenivasan et. al.~\cite{SreenivasanEtAl96}]}}
    \label{fig:Sreeni:2}
  \end{centering}
\end{figure}

For the third order longitudinal structure function $S_3^{\|}$, Kolmogorov derived what is considered an {\em exact result in turbulence}, the famous {\em $4/5$-law}, which states that 
\begin{align}
S_3^{\|}(\ell) \sim -\frac{4}{5} \eps \ell
\label{eq:4/5:law}
\end{align}
holds in the infinite Reynolds number limit, for $\ell \ll \ell_I$. Identity \eqref{eq:4/5:law} is remarkable because a-priori, there is no good reason for the cubic power of the longitudinal increments to have a sign,\footnote{We also note here that the $4/5$-law has an analogue in Lagrangian variables, the so-called Ott-Mann-Gaw\c{e}dzki relation~\cite{OttMann00,Falkovich01}. It relates the anomalous dissipation rate $\eps$ to the time-asymmetry in the rate of dispersion of Lagrangian particles in a turbulent flow. This {\em Lagrangian arrow of time} may be proven rigorously under mild assumptions, see the recent work~\cite{Drivas18}.} on average. Moreover, in addition to claiming that $\zeta_3 = 1$, \eqref{eq:4/5:law} predicts the {\em universal pre-factor} of $-\sfrac 45$. Compelling experimental support of the $4/5$-law is provided for instance by the measurement in Figure~\ref{fig:Sreeni:2}. 
From a mathematical perspective the $4/5$-law is particularly intriguing because under quite mild assumptions one may establish it rigorously. We refer the reader to the results and excellent discussions in~\cite{NieTanveer99,Eyink02}, where evidence is provided that in the inviscid limit, i.e.\ for the Euler equations, the \eqref{eq:4/5:law} should hold  with just local space-time averaging in $(x,t)$ and angular averaging over the direction of the separation vector $\hat z$ (without the assumption of isotropy). This viewpoint is intimately related to Onsager's predictions discussed in Section~\ref{sec:turbulence:Onsager} below.\footnote{See also~\cite{BCZPSW18} for a derivation of the $4/5$-law in the context of the stochastic Navier-Stokes equations, with forcing which is white in time and colored in space, under the seemingly very mild assumption of {\em weak anomalous dissipation}: $\lim_{\nu \to 0} \nu {\mathbb E} \norm{v^\nu}_{L^2}^2 = 0$. Here $v^\nu$ is a stationary martingale solution.}

For $p=2$, from \eqref{eq:Sp:Kolmogorov}--\eqref{eq:K:zeta:p} the Kolmogorov prediction yields $\zeta_2 = \sfrac 23$. One may translate this scaling of the second order structure function into the famous $-\sfrac 53$ energy density spectrum, defined in terms of Fourier projection operators as follows. For $\kappa >0$ the mean kinetic energy per unit mass carried by wavenumber $\leq \kappa$ in absolute value is given by $\frac 12 \average{|\Proj_{\leq \kappa} v^\nu|^2}$. The energy spectrum is then defined as 
\begin{align}
{\mathcal E}(\kappa) = \frac 12 \frac{d}{d\kappa} \average{|\Proj_{\leq \kappa} v^\nu|^2} 
\end{align}
so that the total kinetic energy may be written as $\frac 12 \average{|v^\nu|^2} = \int_0^\infty {\mathcal E}(\kappa) d\kappa$. The Kolmogorov prediction $\zeta_2 = \sfrac 23$ then translates into 
\begin{align}
	{\mathcal E}(\kappa) \sim \epsilon^{2/3} \kappa^{-5/3},
	\label{eq:K:41:spectrum}
\end{align}
for $\kappa^{-1}$ in the inertial range, and in the infinite Reynolds number limit. See~\cite{Frisch95} for experimental support for \eqref{eq:K:41:spectrum}. This power law requires however that velocity fluctuations are uniformly
distributed over the three dimensional domain, which as discussed in Section~\ref{sec:intermittent} below, is not always justified (see Figure~\ref{fig:Sreeni:3}).

\subsection{Basics of the Onsager ('49) theory}
\label{sec:turbulence:Onsager}
In his famous paper on statistical hydrodynamics, Onsager~\cite{Onsager49} considered the possibility that ``turbulent energy dissipation [...] could take place just as readily without the final assistance of viscosity [...] because the velocity field does not remain differentiable''. The pointwise energy balance for smooth solutions $v$ of the Euler equations~\eqref{eq:Euler} is
\begin{align}
\partial_t \frac{|v|^2}{2} + \nabla \cdot \left( v \left(\frac{|v|^2}{2} + p \right) \right) = f  \cdot v \ .
\label{eq:Euler:formal:energy:balance}
\end{align}
Integrating over the periodic domain we obtain the kinetic energy balance
\begin{align}
\frac{d}{dt} \fint_{\T^3} \frac{|v|^2}{2} dx = \fint_{\T^3} f \cdot v  dx \, ,
\label{eq:Euler:global:energy:balance}
\end{align}
which becomes a conservation law when $f\equiv 0$. Onsager is referring to the fact that if the solution $v$ of \eqref{eq:Euler} is not sufficiently smooth, i.e. it is a {\em weak solution}, then the energy balance/conservation \eqref{eq:Euler:global:energy:balance} cannot be justified. Onsager's remarkable analysis went further and made a precise statement about the necessary regularity of $v$ which is required in order to justify \eqref{eq:Euler:global:energy:balance}. This has been phrased in mathematical terms as the {\em Onsager Conjecture} (see Conjecture~\ref{c:Onsager} below). We refer to the review articles~\cite{EyinkSreeniviasan06,Shvydkoy10,Eyink18} for a detailed account of the Onsager theory of {\em ideal turbulence}, and present here only some of the ideas (in terms of Fourier projection operators, as in Onsager's work~\cite{Onsager49}).

We regularize a weak solution $v$ of the Euler equations~\eqref{eq:Euler}, by a smooth cutoff in the Fourier variables at frequencies $ \leq \kappa$, and consider the kinetic energy of $\Proj_{\leq \kappa} v$. Then, similarly to \eqref{eq:Euler:formal:energy:balance}--\eqref{eq:Euler:global:energy:balance} we obtain that 
\begin{align}
\frac{d}{dt} \fint_{\T^3} \frac{|\Proj_{\leq \kappa} v|^2}{2} dx = \fint_{\T^3} \Proj_{\leq \kappa}f \cdot \Proj_{\leq \kappa} v  dx - \Pi_{\kappa}
\label{eq:Euler:global:energy:balance:2}
\end{align}
where as in~\cite{Onsager49,Eyink94,ConstantinETiti94,Frisch95,CCFS08} we denote by $\Pi_\kappa$ the mean energy flux through the sphere of radius $\kappa$ in frequency space, i.e.
\begin{align}
\Pi_\kappa = - \fint_{\T^3} \Proj_{\leq \kappa}(v \otimes v) \colon \nabla \Proj_{\leq \kappa} v dx =    \fint_{\T^3} \underbrace{\Big((\Proj_{\leq \kappa} v \otimes \Proj_{\leq \kappa} v) -  \Proj_{\leq \kappa}(v \otimes v) \Big) \colon \nabla \Proj_{\leq \kappa} v}_{=: \pi_\kappa (x,t)} dx  \,.
\label{eq:Onsager:flux:def}
\end{align}
The above defined mean energy flux $\Pi_\kappa$, and corresponding density $\pi_\kappa$ may also be computed as in the right side of \eqref{eq:CET:main} below, with $\ell \approx \kappa^{-1}$.\footnote{Recalling the notation $D(v)$ from \eqref{eq:KHM} for the measure obtained from the K\'arm\'an-Howarth-Monin relation, we note that Duchon-Robert~\cite{DuchonRobert00} proved that if $v \in L^{3}_{x,t}$ is a weak solution of the Euler equations, then $D(v) = \lim_{\kappa \to \infty} \pi_\kappa$ (in the sense of distributions). Thus, setting $\nu = 0$ in \eqref{eq:NSE:local:energy:balance} we obtain a pointwise balance relation which is valid for weak solutions of the Euler equation. In fact, in \cite{DuchonRobert00} it is shown that if $v$ is a strong limit (in $L^3_{x,t}$) of Leray weak solutions $v^\nu$ of \eqref{eq:NSE}, then $\lim_{\nu \to 0} \nu |\nabla v^\nu|^2 + D(v^\nu) = D(v)$, and thus a-posteriori we obtain that $D(v) \geq 0$ in the sense of distributions.} From \eqref{eq:Euler:global:energy:balance:2} we deduce upon passing $\kappa \to \infty$ that the energy balance \eqref{eq:Euler:global:energy:balance} is holds {\em if and only if} the {\em total energy flux}  
\begin{align}
\Pi = \lim_{\kappa \to \infty} \Pi_\kappa \, ,
\end{align}
{\em vanishes}. Onsager's prediction is that in order for $\Pi$ to be nontrivial, and thus for the weak solution $v$ of the Euler equation to be non-conservative, it should not obey $|\delta v(x;z)| \les |z|^{\theta}$ with $\theta > \sfrac 13$ (see Part (a) of Conjecture~\ref{c:Onsager} below).

We emphasize that in 3D turbulent flows the {\em energy transfer} from one scale/frequency to another is observed to be mainly {\em local}, i.e. the principal contributions to $\Pi_\kappa$   come from $\Proj_{\approx \kappa'} v$, with $\kappa' \approx \kappa$. A rigorous estimate on the locality of the energy transfer arises in~\cite{CCFS08}, where it is proven that 
\begin{align}
  | \Pi_{2^j}| \les  \sum_{i=1}^\infty 2^{-2/3 |j - i|} 2^{i} \|\Proj_{\approx 2^i} v\|_{L^3}^3 \, .
  \label{eq:LP:piece:est}
\end{align}
Estimate \eqref{eq:LP:piece:est} gives the best known condition on $v$ which ensures $\Pi = 0$, namely $v \in L^3_t B^{\sfrac 13}_{3,c_0,x}$ (cf.~\cite{CCFS08}), a condition which is for instance sharp in the case of the 1D Burgers equation~\cite{Shvydkoy10}.

\begin{figure}[h!]
 \begin{centering}
\includegraphics[width=0.45\textwidth]{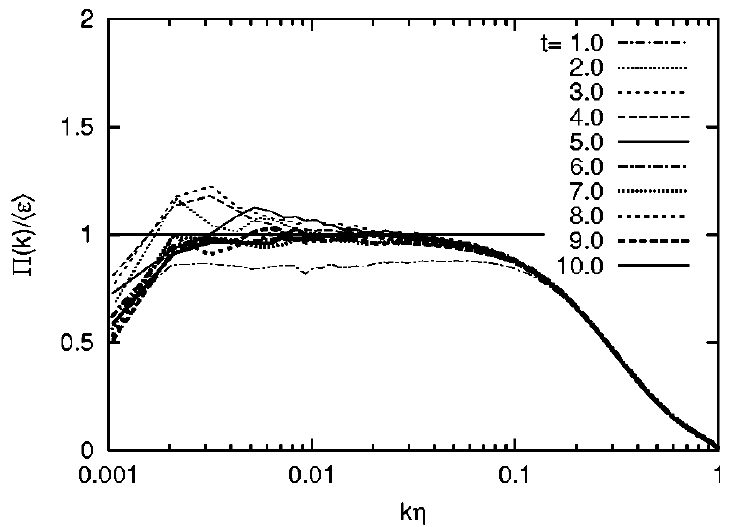}
\caption{{\small Kolmogorov's anomalous energy dissipation rate $\eps$ and Onsager's energy flux $\Pi_\kappa$ appear to agree, at least for $\kappa$ in the inertial range.   [Kaneda et al~\cite{KanedaEtAl03}]: $2048^3$ DNS runs on Earth Simulator computing system.}}
\label{fig:Kaneda}
\end{centering}
\end{figure}

It is not an accident that the $\sfrac 13$-derivative singularities required by Onsager for a dissipative anomaly $\Pi \neq 0$, matches Kolmogorov's assumed $\sfrac 13$ local self-similarity exponent required for $\eps >0$. As already observed by Onsager~\cite{Onsager49}, if $v$ is a weak solution of the Euler equations which is a strong limit  of a sequence $\{v^\nu\}$ of Navier-Stokes solutions for which the anomalous dissipation of energy \eqref{eq:anomaly} holds, then the total energy flux associated to $v$ must match this dissipation anomaly in the vanishing viscosity limit:
\begin{align}
\eps =   \langle \Pi \rangle.
\label{eq:Kolmogorov:Onsager}
\end{align}
On the experimental side, the evidence for \eqref{eq:Kolmogorov:Onsager} is quite convincing, see e.g.~Figure~\ref{fig:Kaneda}.
The energy flux  provides a connection between the Kolmogorov and Onsager theories, and a physics derivation of \eqref{eq:Kolmogorov:Onsager} is as follows~\cite{Frisch95}. Assume for simplicity that $f^\nu = f$ is statistically stationary, and that $\Proj_{\kappa_I} f  = f$ for some integral frequency $\kappa_I$. Denote by $\Pi_\kappa^\nu$ the energy flux through the frequency ball of radius $\kappa$ for a solution $v^\nu$ of the Navier-Stokes equation, i.e. replace $v$  in \eqref{eq:Onsager:flux:def} with $v^\nu$. Then similarly to \eqref{eq:Euler:global:energy:balance:2},   since the ensemble/long-time average $\average{\cdot}$ is stationary, we obtain that 
\begin{align}
  \average{\Pi_\kappa^\nu}  + \nu \average{|\nabla \Proj_{\leq \kappa} v^\nu|^2} = \average{f \cdot \Proj_{\leq \kappa} v^\nu}   \label{eq:weird:energy:val:*}
\end{align}
for $\kappa \geq \kappa_I$.
On the other hand, assuming that the Euler solution is statistically stationary,  \eqref{eq:Euler:global:energy:balance:2} yields 
\begin{align}
\average{\Pi_\kappa}    = \average{f \cdot \Proj_{\leq \kappa} v} 
\label{eq:weird:energy:val:**}
\end{align}
To conclude, we recall that from the definition \eqref{eq:KHM} we have $\average{D(v^\nu)} = \lim_{\kappa \to \infty} \average{\Pi_\kappa^\nu}$ (cf.~\cite{DuchonRobert00}), and with $\eps^\nu$ as given by \eqref{eq:epsilon:nu}, we pass $\kappa\to \infty$ in \eqref{eq:weird:energy:val:*} and \eqref{eq:weird:energy:val:**}, to arrive at 
\begin{align}
\eps - \average{\Pi} = \lim_{\nu \to 0} \left( \eps^\nu - \average{\Pi} \right) = \lim_{\nu \to 0} \average{f \cdot (v^\nu - v)} = 0 
\end{align}
since we assumed $v^\nu \to v$. Note that here we have made a number of assumptions which are not justified.  

\subsection{Intermittency}
\label{sec:intermittent}
In this last part of Section~\ref{sec:physics} we consider the intermittent nature of turbulent flows. This is a topic of significant debate~\cite{Kolmogorov62,mandelbrot1974intermittent,FSN78,FrischParisi85,Sreeni85,MeneveauSreenivasan91,SheLeveque94,Yakhot01}. Large parts of the books~\cite{Frisch95,MoninYaglom13} are dedicated to this mystery, and we refer to these texts, and to the recent papers~\cite{GibbonDoering03,CS2014} for a more detailed discussion. Our interest stems from the fact that this phenomenon is the primary motivation for the {\em intermittent convex integration scheme} discussed in Section~\ref{sec:NSE:L2}.

\begin{figure}[h!]
\begin{centering}
    \includegraphics[width=0.45\textwidth]{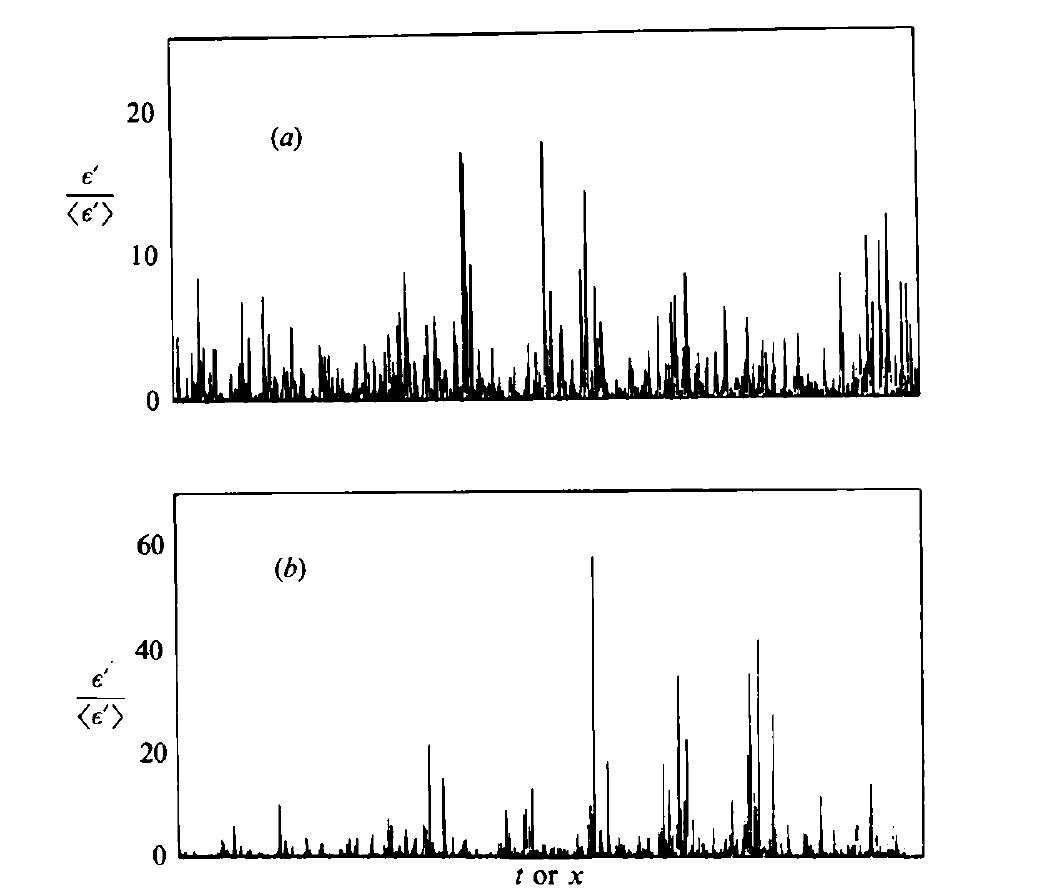}
    \caption{{\small [Meneveau \& Sreenivasan~\cite{MeneveauSreenivasan91}]. Two typical signals, with $\eps' = (du_1/dt)^2$, serving as the surrogate of the energy dissipation
rate, are plotted here upon normalizing by their mean values. Graph (a) was in a laboratory boundary layer at a moderate Reynolds number. Graph (b) was obtained in the atmospheric surface layer at high Reynolds number ($10^4$). The increased intermittency at the higher Reynolds number is clear.}}
    \label{fig:Sreeni:3}
    \end{centering}
\end{figure}

In a broad sense, intermittency is characterized as a deviation from the Kolmogorov 1941 laws. Already in 1942 Landau remarked that the rate of energy dissipation in a  fully developed turbulent flow  is observed to be spatially and temporally inhomogeneous, and thus Kolmogorov's homogeneity and isotropy assumptions need not be valid (cf.~\cite{Landau59,Frisch95}). Figure~\ref{fig:Sreeni:3} shows a typical signal used in experiments to measure $\eps^\nu$. The main feature seems to be the presence of sporadic dramatic events, during which there are large excursions away from the average.

\begin{figure}[h!]
\begin{centering}
    \includegraphics[width=0.5\textwidth]{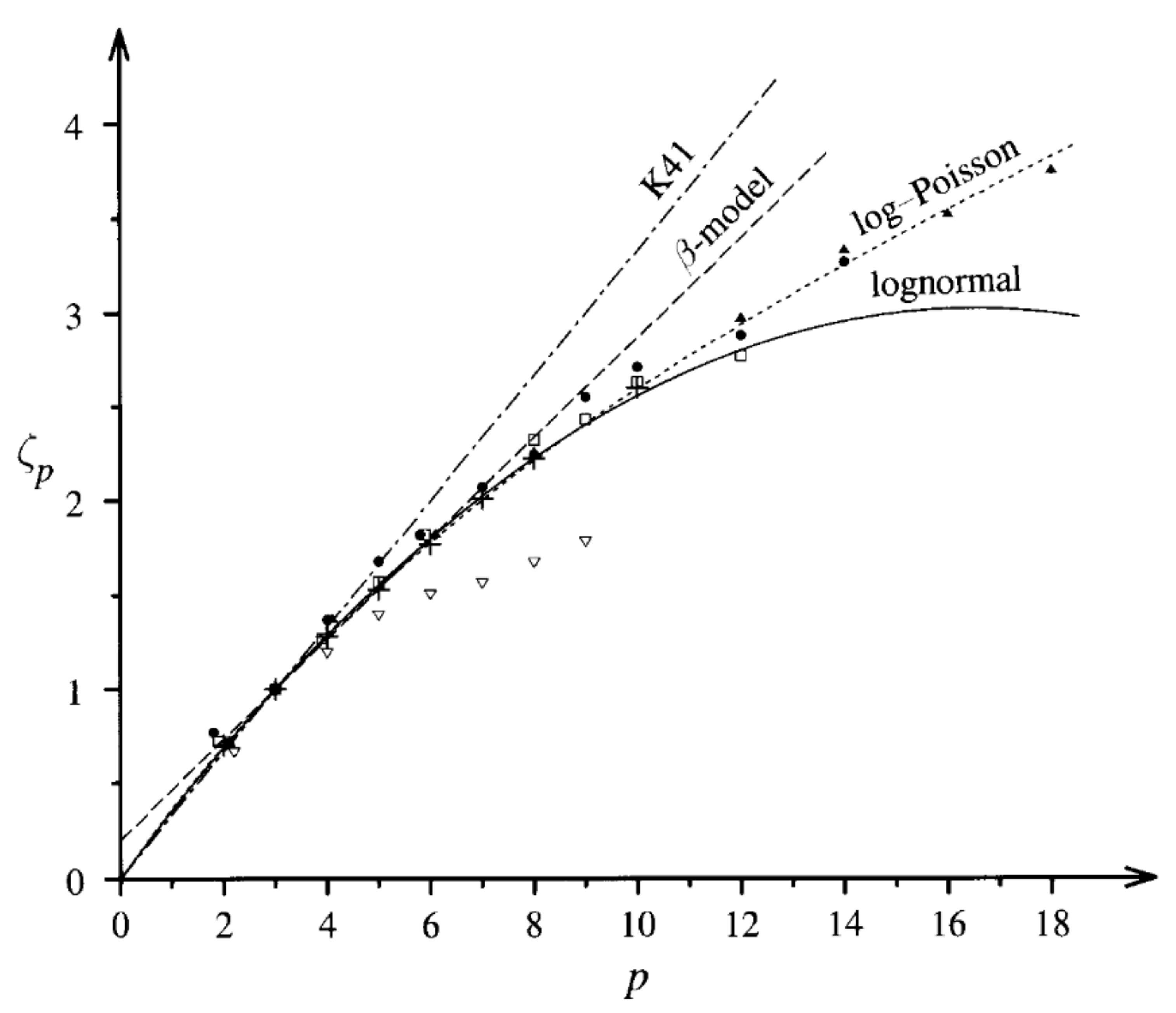}
    \caption{{\small [Frisch~\cite{Frisch95}]: Anomalous scaling of structure functions. Data: inverted white triangles: [Van Atta-Park~\cite{VanAttaPark72}]; black circles, white squares, black triangles: [Anselmet, Gagne, Hopfinger, Antonia~\cite{AnselmetEtAl84}] at $\Re = 515, 536, 852$; $+$ signs: from S1 ONERA wind tunnel.}}
    \label{fig:Zeta:p}
\end{centering}
\end{figure}

A common signature of intermittency is that the structure function exponents $\zeta_p$ deviate from the Kolmogorov predicted value of $\sfrac{p}{3}$, and moreover, for $p\neq 3$ they do not appear to be universal. Figure~\ref{fig:Zeta:p}, compiled by Frisch in~\cite{Frisch95}, highlights this fact. We again see in Figure~\ref{fig:Zeta:p} that the prediction $\zeta_3= 1$ seems to be confirmed by all experimental data, but for $p>3$ we have $\zeta_p < \sfrac p3$, while for $p<3$ we have $\zeta_p > \sfrac p3$.
While there are many phenomenological theories\footnote{For instance: the log-normal model of~\cite{Kolmogorov62}, the $\beta$-model of~\cite{FSN78}, the multifractal model~\cite{FrischParisi85}, the log-Poisson model of~\cite{SheLeveque94}, or the mean-field there of~\cite{Yakhot01}. Interestingly, all these models predict $\zeta_2 > 0.694$. See~\cite[Chapter 8]{Frisch95} for a detailed discussion.} for predicting the structure function exponents $\zeta_p$ in intermittent turbulent flows, none of them seem  to be able to explain all experimental data, and their connection to dynamical evolution of the underlying Navier-Stokes/Euler equations seems to be limited.\footnote{As noted in~\cite{Foias97}, the deterministic bounds on the structure function exponents which one may rigorously establish from the Navier-Stokes equations~\cite{Constantin94,ConstantinFefferman94,ConstantinNieTanveer99} always seem to be bounded from above by the phenomenological predictions.}

A particularly appealing intermittency model is the $\beta$-model of Frisch-Sulem-Nelkin~\cite{FSN78}, which was revisited recently by Cheskidov-Shvydkoy~\cite{CS2014} from a modern analysis perspective. In order to make a connection with the measure-theoretic support of the defect measure $D(v)$, the authors in~\cite{CS2014}  define active regions $A_q$ whose volumes $V_q$ are given in terms of an $L^2 - L^3$ skewness factor which measures the saturation of the Bernstein
inequalities at frequencies $\approx 2^q$. More precisely, using the active volumes defined as
\[ V_q =  L^3 \frac{\average{|\Proj_{\approx 2^q} v|^2}^3}{\average{|\Proj_{\approx 2^q} v|^3}^{2}}, \] 
the intermittency dimension $D$ is defined as
\[ D = 3 - \liminf_{q\to \infty} \frac{\log_2(L^3 V_q^{-1})}{q}, \]
and then the $\beta$-model yields corrections~\cite{FSN78,Frisch95,CS2014} to the Kolmogorov predictions \eqref{eq:K:zeta:p} and \eqref{eq:K:41:spectrum} as
\[ \zeta_p = \frac{p}{3} + (3-D) \left(1 - \frac{p}{3} \right) \qquad \mbox{and} \qquad {\mathcal E}(\kappa) \sim \eps^{- \sfrac 23} \kappa^{-\sfrac 53} \left(\frac{\kappa_I}{\kappa}\right)^{1- \sfrac D3} . \] 
Note that the Kolmogorov theory corresponds to $D=3$, in which the turbulent events fill space. On the other hand, the simulation of~\cite{KanedaEtAl03} estimates $D \approx 2.7$. 

We note in closing that in the intermittent convex integration construction of Section~\ref{sec:NSE:L2}, it is essential that the building blocks concentrate  on a set with dimension strictly less than one. This translates (cf.~estimate \eqref{e:W_Lp_bnd}) into the fact that the skewness ratio of our intermittent building blocks $\frac{\norm{\cdot}_{L^1}}{\norm{\cdot}_{L^2}}$ scales better than the frequency of the building blocks to the power $-1$. This is one of the essential aspects of the construction, and is discussed in detail in Sections~\ref{sec:convex:integration} and~\ref{sec:NSE:L2}.

\section{Mathematical results}
\label{sec:math}
\subsection{The Euler Equations}\label{ss:euler_results}

Local well-posedness for smooth solutions to the Euler equations is classical \cite{Lichtenstein29} (cf.\ \cite{MB2002}).\footnote{In some critical spaces, the Euler equations are known to be ill-posedness in the sense of Hadamard \cite{BourgainDong15,ElgindiMasmoudi14,ElgindiJeong17}.}  By the Beale-Kato-Majda criterion, global well-posedness for the Euler equations is known to hold under the assumption that the $L^{\infty}$ norm of the vorticity is $L^1$ integrable in time \cite{BKM84}. In 2D, vorticity is transported, leading to global well-posedness of smooth solutions \cite{Hoelder33,Wolibner33} as well as weak solutions with $L^{\infty}$ bounded vorticity \cite{Yudovich63,MB2002,Vishik99}. Global well-posedness for smooth solutions to the 3D Euler equation is famously unresolved, and is intimately related to the Clay Millennium problem \cite{F2006}. Indeed, there exists numerical evidence to suggest that the 3D Euler equations may develop a singularity \cite{GuoHou14,GuoHou14b} (cf. \cite{PS92,Gibbon08,GBK08,ChoiHouEtAl17}). Recently, Elgindi and Jeong demonstrated the formation of a singularity  in the presence of a conical hourglass-like boundary~\cite{ElgindiJeong18}.

Within the class of weak solutions, the Euler equations are known to display paradoxical behavior. In the seminal work \cite{Scheffer93}, Scheffer demonstrated the existence of non-trivial weak solutions with compact support in time (cf.~\cite{Shnirelman97,Shnirelman00}). These results represent  a quite drastic demonstration of non-uniqueness for weak solutions to the Euler equations. For the purpose of this article, we define a weak solution for \eqref{eq:Euler} as:
\begin{definition}
\label{eq:def:weak:Euler}
A vector field $v \in C^0_tL^2_x$ is called a weak solution of the Euler equations if for any $t$ the vector field $v(\cdot,t)$ is weakly divergence free, has zero mean, and satisfies the Euler equation distributionally:  
\[
\int_{\R} \! \int_{\T^3} v \cdot (\partial_t \varphi + (v \cdot \nabla) \varphi ) dx dt = 0\,,
\]
for any divergence free test function $\varphi$.\footnote{Note the pressure can be recovered by the formula $-\Delta p=\div\div(v\otimes v)$ with $p$ of zero mean.} For a weak solution to the Cauchy problem this definition is modified in the usual way.
\end{definition}

As mentioned in Section~\ref{sec:physics}, one motivation for studying weak solutions to the Euler equations, is that in the inviscid limit, turbulent solutions exhibiting a dissipation anomaly are necessarily weak solutions. In~\cite{Onsager49} Onsager conjectured the following dichotomy:
\begin{conjecture}[Onsager's conjecture]\label{c:Onsager}\ 
\begin{enumerate}[(a)]
\item\label{Onsager:a} Any weak solution $v$ belonging to the H\"older space $C^\theta_{x,t}$ for $\theta>\sfrac{1}{3}$ conserves kinetic energy.
\item\label{Onsager:b} For any $\theta<\sfrac{1}{3}$ there exist weak solutions $v\in C^\theta_{x,t}$ which dissipate kinetic energy.
\end{enumerate}
\end{conjecture}
Part~\eqref{Onsager:a} of this conjecture was partially established by Eyink in \cite{Eyink94}, and later proven in full by Constantin, E and Titi in \cite{ConstantinETiti94} (see also \cite{DuchonRobert00,CCFS08}, and the more recent work~\cite{Shvydkoy18}, for refinements).  The proof follows by a simple commutator argument: Suppose $v$ is a weak solution to the Euler equations, and let  $v_{\ell}$ be the spatial mollification of $v$ a length scale $\ell$. Then, $v_{\ell}$ satisfies
\[\int_{\T^3}\abs{v_{\ell}(x,t)}^2dx-\int_{\T^3}\abs{v_{\ell}(x,0)}^2dx=2\int_0^t\int_{\T^3} \tr((v\otimes v)_{\ell}(\nabla v_{\ell}))dxds~.\]
Applying the identity 
\[\int_{\T^3} \tr((v_{\ell}\otimes v_{\ell})(\nabla v_{\ell}))dx\equiv 0~,\]
yields
\begin{align}
\int_{\T^3}\abs{v_{\ell}(x,t)}^2dx-\int_{\T^3}\abs{v_{\ell}(x,0)}^2dx=2\int_0^t\int_{\T^3} \tr\left(((v\otimes v)_{\ell}-(v_{\ell}\otimes v_{\ell}))(\nabla v_{\ell})\right)dxd\tau \, .
\label{eq:CET:main}
\end{align}
Applying the commutator estimate Proposition \ref{p:CET}, in Section \ref{sec:Euler:C0}, we deduce
\[\abs{\int_{\T^3}\abs{v_{\ell}(x,t)}^2dx-\int_{\T^3}\abs{v_{\ell}(x,0)}^2dx}\leq C \ell^{3\theta-1}\norm{v}_{C^\theta}^3~.\]
Thus, if $\theta>\sfrac13$, the right hand side converges to zero as $\ell\rightarrow 0$.

Concerning part~\eqref{Onsager:b} of the Onsager conjecture, strictly speaking the weak solutions of Scheffer are not dissipative, as dissipative solutions are required to have non-increasing energy. The existence of dissipative weak solutions to the Euler equations was first proven by Shnirelman in~\cite{Shnirelman00} (cf.\ \cite{DLSZ09,DeLellisSzekelyhidi10}). In the groundbreaking papers \cite{DLSZ13,DeLellisSzekelyhidi12a}, De Lellis and Sz\'ekelyhidi Jr.\ made significant progress towards Part (b) of Onsager's conjecture by proving the first construction of dissipative H\"older continuous weak solutions to the Euler equations (see~Theorem~\ref{thm:Euler:C0:DLSZ}). After a series of  advancements \cite{IsettThesis,BDLISZ15,Buckmaster15,BDLSZ16,DSZ17}, part (b) of the Onsager conjecture was resolved by Isett in \cite{Isett16}:
\begin{theorem}[Theorem~1,~\cite{Isett16}]
\label{thm:Onsager:critical:Isett}
For any $\beta \in (0,\sfrac 13)$ there exists a nonzero weak solution $v \in C^\beta(\T^3 \times \R)$, such that $v$ vanishes identically outside of a finite interval.
\end{theorem}

Like the original paper of Scheffer~\cite{Scheffer93}, the weak solutions constructed by Isett~\cite{Isett16} are not strictly dissipative. This technical issue was resolved in the paper \cite{BDLSV17}, in which the precise statement of part (b) was proven:
\begin{theorem}[Theorem~1.1,~\cite{BDLSV17}]
\label{thm:Onsager:critical:dissipative}
Let $e\colon[0,T]\to \R$ be a strictly positive smooth function. For any $\beta \in (0,\sfrac 13)$ there exists a weak solution $v\in C^\beta(\T^3 \times [0,T])$ of the Euler equations \eqref{eq:Euler}, whose kinetic energy at time $t \in [0,T]$ equals $e(t)$.
\end{theorem} 
 
The exponent $\sfrac13$ in Onsager's conjecture can be viewed in terms of a larger class of \emph{threshold exponents} at which  a dichotomy in the behavior of solutions arises. In a recent expository paper~\cite{Klainerman17} on the work of Nash, Klainerman considered various threshold exponents in the context of non-linear PDE (see also \cite{BSV16} for a discussion of thresholds exponents in the context of hydrodynamic equations). In order to simply the discussion, consider Banach spaces of the form $X^\alpha  = C_t^0 C^\alpha_x$. As in~\cite[Page 11]{Klainerman17}, let us define the following exponents: 
\begin{itemize}
\item The {\em scaling exponent} $\alpha_*$ determines the norm for which the $X^{\alpha_*}$  is invariant under the natural scalings of the equation.
\item The {\em Onsager exponent} $\alpha_O$ determines the norm for which the Hamiltonian of a PDE is conserved.
\item The
 {\em Nash exponent} $\alpha_N$ determines the threshold for which the PDE is flexible or rigid in the sense of the $h$-principle. 
 \item The {\em uniqueness exponent} $\alpha_U$ determines the threshold for which uniqueness of solutions holds.
  \item The {\em well-posedness exponent} $\alpha_{WP}$ determines the threshold for which local well-posedness holds. 
\end{itemize}
Since flexibility implies non-uniqueness, and well-posedness implies uniqueness, we have $\alpha_N \leq \alpha_U \leq \alpha_{WP} $. For the Euler equations,  $\alpha_{WP} = 1$ (cf.~\cite{Hoelder33,BaTi2010,BourgainDong15,ElgindiMasmoudi14}), $\alpha_U$ is  conjectured to be $1$ (the Beale-Kato-Majda criterion implies that $\alpha_U \leq 1$), $\alpha_O = 1/3$ (cf.~\cite{ConstantinETiti94,CCFS08,Isett16,BDLSV17}), and $\alpha_*= 0$. In general, one expects the ordering  $\alpha_* \leq \alpha_O \leq \alpha_N \leq \alpha_U \leq \alpha_{WP}$ (cf.~\cite[Equation (0.7)]{Klainerman17}). 

\subsection{The Navier-Stokes Equations}

The global well-posedness for the 3D Navier-Stokes equations, is one of the most famous open problems in mathematics, subject to one of seven \emph{Clay Millennium Prize} problems \cite{F2006}. Local well-posedness in various scale\footnote{Recall that if $v(x,t)$ is a solution of \eqref{eq:NSE}, then so is $v_\lambda(x,t) = \lambda v(\lambda x,\lambda^2 t)$ for every $\lambda > 0$.} invariant spaces follows by classical contraction mapping arguments \cite{FujitaKato64,Kato84,KochTataru01,Lemarie16} and global well-posedness typically follows when the datum is small in these spaces. If one relaxes one's notion of solutions and considers instead \emph{weak solutions}, then Leray~\cite{Leray34} and later Hopf~\cite{Hopf51} proved that for any finite energy initial datum there exists a global weak solution to the Navier-Stokes equation. More precisely, Leray proved the global existence in the following class of weak solutions:
\begin{definition}\label{d:leray-hopf} A vector field $v \in C^0_{\rm weak}([0,\infty);L^2(\mathbb T^3) )\cap L^2([0,\infty);\dot{H}^1(\mathbb T^3))$ is called a Leray-Hopf weak solution of the Navier-Stokes equations if for any $t \in \mathbb R$ the vector field $v(\cdot,t)$ is weakly divergence free, has zero mean, satisfies the Navier-Stokes equations distributionally: 
\begin{align*}
\int_{\R} \! \int_{\T^3} v \cdot (\partial_t \varphi + (v \cdot \nabla) \varphi + \nu \Delta \varphi ) dx dt +  \int_{\T^3} v(\cdot,0) \cdot \varphi(\cdot,0) dx = 0\,,
\end{align*} for any divergence free test function $\varphi$,
and satisfies the 
\emph{energy inequality}:
\begin{equation}\label{e:energy_ineq}
\frac{1}{2}\int_{\mathbb T^3}\abs{v(x,t)}^2\,dx+\int_{\mathbb T^3\times [0,t]}\abs{\nabla v(x,s)}^2\,dxds
\leq \frac{1}{2}\int_{\mathbb T^3}\abs{v(x,0)}^2\,dx\,.
\end{equation}
\end{definition}
Leray-Hopf solutions are known to be regular and unique under the additional assumption that one of the  Lady{\v{z}}enskaja-Prodi-Serrin conditions are satisfied, i.e.\ the solution is bounded in a scaling invariant space $L^p_t L^q_x$ for  $\sfrac 2p + \sfrac 3q = 1$ \cite{KiselevLadyzhenskaya57,Prodi59,Serrin62,EscauriazaSerginSverak03}. 
One possible strategy to proving that the Navier-Stokes equation is well-posed is then to show that the weak solutions are smooth \cite{F2006}.  Since smooth solutions are necessarily unique, such a result would imply the uniqueness of weak solutions. 

In recent work by the authors \cite{BV}, another class of weak solutions was considered, namely:
 \begin{definition}\label{d:weak_sol} A vector field $v \in C^0_tL^2_x$ is called a weak solution of the Navier-Stokes equations if for any $t $ the vector field $v(\cdot,t)$ is weakly divergence free, has zero mean, and 
\begin{align*}
\int_{\R} \! \int_{\T^3} v \cdot (\partial_t \varphi + (v \cdot \nabla) \varphi + \nu \Delta \varphi ) dx dt = 0\,,
\end{align*}
for any divergence free test function $\varphi$.
\end{definition}
The above class is weaker than Definition \ref{d:leray-hopf} in the sense that solutions need not satisfy the energy inequality \eqref{e:energy_ineq}; however, they are stronger in the sense that the $L^2$ norm in space is required to be strongly continuous in time. Such solutions satisfy the integral equation~\cite{FabesJonesRiviere72} 
\begin{align*}
v(\cdot,t) =  e^{\nu \Delta t}v(\cdot,0)+\int_0^t e^{\nu \Delta(t-s)}\mathbb P \div(v(\cdot,s)\otimes v(\cdot,s)) ds \, ,
\end{align*}
and are sometimes called {\em mild} or \emph{Oseen} solutions (cf.~\cite[Definition 6.5]{Lemarie16}).  As is the case for Leray-Hopf solutions, weak solutions of the form described in Definition \ref{d:weak_sol} are known to be regular under the additional assumption that one of  Lady{\v{z}}enskaja-Prodi-Serrin conditions is satisfied \cite{FabesJonesRiviere72,FurioliLemarieRieussetTerraneo00,LionsMasmoudi01,LemarieRieusset02,Kukavica06b,Germain06}. The principal result of \cite{BV} is:
\begin{theorem}[Theorem 1.2,~\cite{BV}]
\label{thm:BV:main}
There exists $\beta > 0$, such that the following holds. For any nonnegative smooth function $e(t) \colon [0,T] \to [0,\infty)$, and any $\nu \in (0,1]$, there exists a weak solution
of the Navier-Stokes equations $v\in C^0([0,T];H^\beta(\T^3)) \cap C^0([0,T];W^{1,1+\beta}(\T^3))$, such that $\int_{\T^3} |v(x,t)|^2 dx = e(t)$ holds for all $t\in [0,T]$. 
\end{theorem}

Since the energy profile may be chosen to have compact support, and $v\equiv 0$ is a solution, the result implies the non-uniqueness of weak solutions to the Navier-Stokes equations, in the sense of Definition~\ref{d:weak_sol}. Theorem~\ref{thm:BV:main} represents a failure of the strategy of proving global well-poseness via weak solutions, at least for the class of weak solutions defined in Definition~\ref{d:weak_sol}. 

One may naturally ask if such non-uniqueness holds for Leray-Hopf weak solutions. This problem remains open. Non-uniqueness of Leray-Hopf solutions were famously conjectured by Lady{\v{z}}enskaja \cite{Lady67}. More recently, \v{S}ver\'ak and Jia proved the non-uniqueness of Leray-Hopf weak solutions assuming that a certain spectral assumption holds \cite{JiaSverak15}. While Guillod and \v{S}ver\'ak have provided in \cite{GuillodSverak17} compelling numerical evidence that the assumption of~\cite{JiaSverak15} may be satisfied, a rigorous proof remains to date elusive.\footnote{If one considers the analog of a Leray-Hopf solutions for the fractional Navier-Stokes equation, where the Laplacian is replaced by the fractional Laplacian $(-\Delta)^{\alpha}$, then non-uniqueness is known to hold for $\alpha < \sfrac13$ in view of the recent works~\cite{CDLDR17,DR08}.}

An alternate, stronger version of Leray-Hopf solutions is often considered in the literature:
\begin{definition}\label{d:loc_leray-hopf} A Leray-Hopf weak solution to the Navier-Stokes equation satisfying the \emph{local energy equality}, is a vector field satisfying the same conditions as detailed in Definition \ref{d:leray-hopf}; however, with the energy equality \eqref{e:energy_ineq} replaced with the local energy inequality
\begin{equation}\label{e:loc_energy_ineq}
\frac{1}{2}\int_{\mathbb T^3}\abs{v(x,t)}^2\,dx+\int_{\mathbb T^3\times [t_0,t]}\abs{\nabla v(x,s)}^2\,dxds
\leq \frac{1}{2}\int_{\mathbb T^3}\abs{v(x,t_0)}^2\,dx\,.
\end{equation}
for almost every $t_0 \geq 0$ and all $t>t_0$.
\end{definition}
The advantage of Definition \ref{d:loc_leray-hopf} over Definition \ref{d:leray-hopf} is that from the localized energy inequality \eqref{e:loc_energy_ineq}, one can deduce that the solutions possess epochs of regularity, i.e.\ many time intervals on which they are smooth. Indeed, in \cite{Leray34},  Leray proved that such solutions are  almost everywhere in time smooth since the   \emph{singular set of times} $\Sigma_T$ has Hausdorff dimension $\leq \sfrac 12$. Improving on this, Scheffer~\cite{Scheffer76} proved that the $\sfrac 12$-dimensional Hausdorff measure of $\Sigma_T$ is $0$. More detailed results, concerning the Minkowski dimension have been obtained in~\cite{RobinsonSadowski07,Kukavica09}.

A curious consequence of the partial regularity result of Leray \cite{Leray34}, the local-wellposedness theory and the \emph{weak-strong uniqueness} result of Prodi-Serrin~\cite{Prodi59,Serrin62,Lemarie16,Wiedemann17}, is that if 
a Leray-Hopf solution in the sense of Definition \ref{d:loc_leray-hopf} is not smooth for some time $t>0$, then on an open interval in time the solution would be in fact a strong solution that blows up, implying a negative answer to Millennium prize question. Thus assuming for the moment that the Millennium prize question is out of reach, one is left to prove the non-uniqueness result of Leray-Hopf solutions in the sense of Definition \ref{d:loc_leray-hopf} via a bifurcation at $t=0$ -- this is indeed the strategy employed by \v{S}ver\'ak and Jia \cite{JiaSverak15}. Unfortunately, this suggests that convex integration is perhaps ill-suited for the task of proving non-uniqueness of Leray-Hopf solutions in the sense of Definition \ref{d:loc_leray-hopf}. However, the above argument does not apply in the context of the Leray-Hopf solutions defined in Definition \ref{d:leray-hopf}, and thus the argument does not rule out a proof of non-uniqueness of such solutions via the method of convex integration.

The partial regularity theory for Leray-Hopf solutions leads to the natural question of whether there exists weak, singular solutions to the Navier-Stokes equations that are smooth outside a suitably small set in time. In \cite{BCV18}, jointly with M.~Colombo, the following result was established:
\begin{theorem}[Theorem 1.1,~\cite{BCV18}]
\label{thm:RT}
There exists $\beta>0$ such that the following holds.
For $T>0$, let $u^{(1)}, u^{(2)} \in C^0([0,T];\dot{H^3}(\T^3))$ be two strong solutions  of the Navier-Stokes equations. There exists a weak solution $v \in C([0,T];H^\beta(\T^3))$ and is such that 
\[v(t) = u^{(1)}(t)~for~ t\in [0,\sfrac{T}{3}] \quad\mbox{and}\quad v (t)\ = u^{(2)}(t)~for~ ~t\in [\sfrac{2T}{3},T]\]
Moreover, there exists a zero Lebesgue measure set of times $\Sigma_T \subset (0,T]$ with Hausdorff  dimension less than $1-\beta$, such that $v \in C^\infty((0,T) \setminus \Sigma_T \times \T^3)$.
\end{theorem}
Theorem \ref{thm:RT} represents the first example of a mild/weak solution to the Navier-Stokes equation whose singular set of times $\Sigma_T \subset (0,T]$ is both \emph{nonempty}, and has Hausdorff dimension strictly less than $1$.

In additional to localizing the energy inequality in time, as was done in \eqref{e:loc_energy_ineq}, one can also localize it in space, leading to the \emph{generalized energy inequality} of Scheffer \cite{Scheffer77}:
\begin{equation}
2\iint \abs{\nabla v}^2 \varphi\,dxdt\leq \iint  \left( \abs{v}^2(\varphi_t+\Delta \varphi)+(\abs{v}^2 +2p) v\cdot \nabla \varphi \right)\,dxdt
\label{eq:suitable}
\end{equation}
for any non-negative test function $\varphi $. Weak solutions satisfying the generalized energy inequality are known as {\em suitable weak solutions} \cite{Scheffer77,CaffarelliKohnNirenberg82}. Note that \eqref{eq:suitable} is a restatement of the condition that the defect measure $D[v]$ in \eqref{eq:NSE:local:energy:balance} is non-negative. Following the pioneering work of Scheffer \cite{Scheffer77,Scheffer80}, Cafferelli, Kohn and Nirenberg famously proved that the singular set of suitable weak solutions has zero parabolic 1D Hausdorff measure. Analogously to the case of Leray-Hopf solutions in the sense of Definition \ref{d:loc_leray-hopf}, convex integration methods seem ill-suited for proving the non-uniqueness of suitable weak solutions to the Navier-Stokes equation.

In view of the discussion of Section~\ref{sec:physics}, we are led to consider the question of whether the nonconservative weak solutions to the Euler equations obtained in~\cite{Isett16,BDLSV17} arise as vanishing viscosity limits of weak solutions to the Navier-Stokes equations.\footnote{Vanishing viscosity limits of Leray-Hopf solutions to the Navier-Stokes equations are known to be Lions   {\em dissipative measure-valued} solutions of the Euler equations -- these solutions however do not necessarily satisfy the Euler equations in the sense of distributions. Under additional assumptions, it was in fact shown earlier by Di Perna-Majda~\cite{DiPernaMajda87} that vanishing viscosity limits are measure-valued solutions for \eqref{eq:Euler}. See~\cite{BDLSZ11,Wiedemann17} for the weak-strong uniqueness property in this class. On the other hand, if one assumes an estimate on velocity increments in the inertial range,  which amounts to $\zeta_2>0$, it was shown in~\cite{ConstantinVicol18} that weak limits of Leray solutions are weak solutions of the Euler equation.} In this direction, as a direct consequence of the proof of Theorem~\ref{thm:BV:main}, one obtains:
\begin{theorem}[Theorem~1.3, \cite{BV}]\label{thm:NSE:Euler}
For $\bar \beta >0$ let $v \in C^{\bar \beta}_{t,x}(\T^3 \times [-2T,2T])$ be a zero-mean weak solution of the Euler equations. Then there exists $\beta>0$, a sequence $\nu_n \to 0$, and a uniformly bounded sequence $v^{(\nu_n)} \in C^0_t ([0, T]; H^\beta_x(\mathbb T^3))$ of weak solutions to the Navier-Stokes equations in the sense of Definition~\ref{d:weak_sol}, with $v^{(\nu_n)} \to v$ strongly in $C^0_t([0, T]; L^2_x(\mathbb T^3))$.
\end{theorem}
The above result shows that  being a strong $L^2$ limit of weak solutions to the Navier-Stokes equations, in the sense of Definition~\ref{d:weak_sol}, cannot serve as a selection criterion for weak solutions of the Euler equation.  See also Remark~\ref{rem:vanishing:viscosity} below.

Lastly, in relation to the threshold exponents considered in Section \ref{ss:euler_results}, if
one considers the family of Banach spaces $C_t H^\alpha_x$, then $\alpha_*=\alpha_{WP}=\sfrac 12$ \cite{FujitaKato64,Kato84,KochTataru01,Lemarie16,BP08}. If we relabel $\alpha_{O}$ the exponent in which the energy equality holds, then as a consequence of Theorem \ref{thm:BV:main}, and the simple observation that regular solutions obey the energy equality, $0<\alpha_{O},\alpha_N\leq \sfrac12$. As a consequence of the expected ordering  $\alpha_* \leq \alpha_O \leq \alpha_N \leq \alpha_U \leq \alpha_{WP}$, one would naturally conjecture that $ \alpha_*= \alpha_O = \alpha_N= \alpha_U =\alpha_{WP}=\sfrac12$.

\section{Convex integration schemes in incompressible fluids}
\label{sec:convex:integration}
The method of convex integration can be traced back to the work of Nash, who used it to construct exotic counter-examples to the $C^1$ isometric embedding problem \cite{Nash54} -- a result that was cited in awarding Nash the Abel prize in 2015 (cf.\ \cite{Kuiper55}). The method was later refined by Gromov \cite{Gromov73} and it evolved into a general method for solving \emph{soft/flexible} geometric partial differential equations \cite{EliashbergMishachev02}. In the influential paper \cite{MullerSverak03}, M\"uller and \v{S}ver\'ak adapted convex integration to the theory of differential inclusions (cf.\ \cite{KircheimMullerSverak03}), leading to renewed interest in the method as a result of its greatly expanded applicability.

\subsection{Convex integration schemes for the Euler equations}

Inspired by the work \cite{MullerSverak03,KircheimMullerSverak03}, and building on the plane-wave analysis introduced by Tartar \cite{tartar_compensated_1979,tartar_compensated_1983,diperna_compensated_1985}, De Lellis and Sz\'ekelyhidi Jr., in \cite{DLSZ13}, applied convex integration in the context of weak $L^\infty$ solutions of the Euler equations, yielding an alternative proof  of Scheffer's~\cite{Scheffer93} and Schnirelman's~\cite{Shnirelman00} famous non-uniqueness results. The work \cite{DLSZ13}, has since been extended and adapted by various authors to various problems arising in mathematical physics \cite{DeLellisSzekelyhidi10,CordobaFaracoGancedo11,Shvydkoy11,Wiedemann11,ChiodaroliDeLellisKremlOndvrej15}, see the reviews~\cite{DeLellisSzekelihidi12,Szekelyhidi12,DLSZ17,DLSZ19} and references therein.

In a first attempt at attacking Onsager's famous conjecture on energy conservation, De Lellis and Sz\'ekelyhidi Jr.\  in their seminal paper \cite{DLSZ13} developed a new convex integration scheme, motivated and resembling in part the earlier schemes of Nash and Kuiper \cite{Nash54,Kuiper55}. In \cite{DLSZ13}, De Lellis and Sz\'ekelyhidi Jr.\ demonstrated the existence of continuous weak solutions $v$ to the Euler equations satisfying a prescribed kinetic energy profile, i.e.\ given a smooth function $e:[0,T]\rightarrow \R^+$, there exists a weak solution $v$ such that 
\begin{equation}\label{e:energy_profile}
\frac{1}{2}\int_{\mathbb T^3} \abs{v(t,x)}^2\,dx=e(t)\,.
\end{equation}
See~Theorem~\ref{thm:Euler:C0:DLSZ} below.  The proof proceeds via induction. At each step $q\in\mathbb N$,  a pair $(v_q, \mathring{R}_q)$ is constructed solving the \emph{Euler-Reynolds system}  
\begin{subequations}\label{e:euler_reynolds}
\begin{align}
\partial_t v_q + \div (v_q\otimes v_q) + \nabla p_q =\div\mathring{R}_q\\ 
\div v_q = 0 \,.
\end{align}
\end{subequations} 
such that as $q \to \infty$ the sequence $\mathring{R}_q$ converges uniformly to $0$ and the sequence $v_q$ converges uniformly to a weak solution to the Euler equations \eqref{eq:Euler} satisfying \eqref{e:energy_profile}.

The Euler-Reynolds \eqref{e:euler_reynolds} system arises naturally in the context of computational fluid mechanics.
As mentioned in~\cite{Frisch95}, via~\cite{Lamb93}, the concept of eddy viscosity and microscopic to macroscopic stresses may be traced back to the work of Reynolds~\cite{Reynolds1985}. 
Given a solution  $v$ to \eqref{eq:Euler}, let $\bar v$ be the velocity obtained  through the application of a filter (or averaging operator) that commutes with derivatives, ignoring the unresolved small scales. Then $(\overline v, R)$ is a solution to $\eqref{e:euler_reynolds}$ for $ R=\overline {v\otimes v}-\overline v\otimes \overline v = \overline{ (v-\overline v) \otimes (v-\overline v)}$.   In this context the $3\times 3$ symmetric tensor $R$ is referred to as the \emph{Reynolds stress}.

For comparison, the iterates $(v_q, \RR_q)$ constructed via a convex integration scheme are approximately spatial averages  of the final solution $v$ at length scales decreasing with $q$. Owing to the analogy to computational fluid mechanics, we refer to the symmetric tensor $\mathring{R}_q$ as the Reynolds stress. Without loss of generality, we will also assume $\mathring{R}_q$ to be traceless.

At each inductive step, the perturbation $w_{q+1} = v_{q+1}-v_{q}$ is designed such that the new velocity $v_{q+1}$ solves the Euler-Reynolds system
\begin{equation*}\begin{split}
\partial_t v_{q+1} + \div (v_{q+1}\otimes v_{q+1})+\nabla p_{q+1} &= \div\mathring R_{q+1}\\
\div v_{q+1} &= 0\, .
\end{split}
\end{equation*}
with a smaller Reynolds stress $\RR_{q+1}$. Using the equation for $v_q$ we obtain the following decomposition of $\RR_{q+1}$:
\begin{equation*}\begin{split}
\div  \RR_{q+1} 
&=  \div(w_{q+1}\otimes w_{q+1} - \mathring R_{q})  + \nabla(p_{q+1} - p_q)\\
&\quad +\partial_t w_{q+1} + v_{q} \cdot \nabla w_{q+1} \\
&\quad + w_{q+1} \cdot \nabla v_{q}
\end{split}\,,
\end{equation*}
which we denote (line-by-line) as the \emph{oscillation error}, \emph{transport error} and \emph{Nash error} respectively.  The Reynolds stress $\mathring R_{q+1}$ can then be defined by solving the above divergence equation utilizing an $-1$ order linear differential operator $\mathcal R \approx \div^{-1}$ (see \eqref{eq:RSZ}).

The perturbation $w_{q+1}=v_{q+1}-v_{q}$ is constructed as a sum of highly oscillatory \emph{building blocks}. In earlier papers \cite{DLSZ13,DeLellisSzekelyhidi12a,BDLISZ15,Buckmaster15,IsettOh16,BDLSZ16}, Beltrami waves were used as the building blocks of the convex integration scheme (see Section \ref{sec:Beltrami} for a discussion). In later papers \cite{DSZ17,Isett16,BDLSV17,Isett17}, Mikado waves were employed (see Section~\ref{sec:Mikado}). These building blocks are used in an analogous fashion to the \emph{Nash twists} and \emph{Kuiper corrugations} employed in the $C^1$ embedding problem \cite{Nash54,Kuiper55}. The perturbation $w_{q+1}$ is designed in order to obtain a cancellation between the low frequencies of the quadratic term $\div(w_{q+1}\otimes w_{q+1})$ and the old Reynolds stress error $\mathring R_q$, thereby reducing the size of the oscillation error.  Roughly speaking, the principal part of the perturbation, which we label $w_{q+1}^{(p)}$, will be of the form
\begin{equation}\label{e:perturbation_form}
w_{q+1}^{(p)}\sim \sum_{\xi} a_{\xi}(\mathring R_q )W_{\xi}\,,
\end{equation}
where the $W_{\xi}$ represent the \emph{building blocks} oscillating at a prescribed high frequency $\lambda_{q+1}$, and the coefficient functions $a_{\xi}$ are chosen such that
\begin{equation}\label{e:low_freq_R_cancellation}
\sum_{\xi} a_{\xi}^2(\mathring R_q )  \fint_{\T^3} W_{\xi}\mathring\otimes W_{\xi} = -\mathring R_q\,,
\end{equation}
where $\mathring\otimes$ denotes the trace-free part of the tensor product. As we will see in Section~\ref{eq:compose:with:flow:map:and:get:a:coffee}, the principal part will need to be modified from the form presented in \eqref{e:perturbation_form} in order to minimize the transport error. This will be achieved by flowing the building blocks $W_{\xi}$ along the flow generated by $v_q$ (see Section~\ref{sec:flow:maps:1}). Additionally, in order to ensure that $w_{q+1}$ is divergence free, we will need to introduce a divergence free corrector $w_{q+1}^{(c)}$ such that
\[w_{q+1}=w_{q+1}^{(p)}+w_{q+1}^{(c)}\,,\]
is divergence free.

Heuristically, let us assume for the moment that the frequencies scale geometrically,\footnote{In practice,  it is convenient to use a super-exponentially growing sequence $\lambda_q$ which obeys $\lambda_{q+1} \approx \lambda_q^b$, where $b>1$.} i.e.\
\[\lambda_q=\lambda^q\]
for some large $\lambda\in \mathbb N$. In order that to ensure that the inductive scheme converges to a H\"older continuous velocity $v$ with H\"older exponent $\beta>0$, then by a scaling analysis, the perturbation amplitude is required to satisfy the bound
\begin{equation}\label{e:perturbation_scaling}
\norm{w_{q+1}}_{C^0}\leq \lambda_{q+1}^{-\beta}\,.
\end{equation}
In view of \eqref{e:low_freq_R_cancellation}, this necessitates that the Reynolds stress $\mathring R_q$ obeys the bound
\begin{equation}\label{e:Reynolds_scaling}
\norm{\mathring R_{q}}_{C^0}\leq \lambda_{q+1}^{-2\beta}\,.
\end{equation}
As a demonstration of the typical scalings present in convex integration schemes for the Euler equations, let us consider the Nash error.  Heuristically, since $v_q$ is defined as the sum of perturbations of frequency $\lambda^{q'}$ for $q'\leq q$ and  $w_{q+1}$ is of frequency $\lambda_{q+1}\gg \lambda_{q'}$ for every $q'\leq q$ we have
\begin{align*}
\norm{\mathcal R\left(w_{q+1} \cdot \nabla v_{q}\right)}_{C^{0}}&\les \frac{\norm{w_{q+1}}_{C^{0}}\norm{v_{q}}_{C^{1}}}{\lambda_{q+1}}
\end{align*}
where we recall that $\mathcal R$ is a $-1$ order linear differential operator solving the divergence equation. Applying \eqref{e:perturbation_scaling} and assuming that $\beta<1$ then we obtain
\begin{align*}
\norm{\mathcal R\left(w_{q+1} \cdot \nabla v_{q}\right)}_{C^{0}}&\les \lambda_{q+1}^{-\beta-1}\sum_{q'\leq q} \lambda_{q'}^{1-\beta}\\
&\les \lambda_{q+1}^{-\beta-1}\lambda_{q}^{1-\beta}\\
&\les \lambda_{q+2}^{-2\beta}\lambda^{3\beta-1}
\end{align*}
Thus in order to ensure that $\mathring R_{q+1}$ satisfying the bound \eqref{e:Reynolds_scaling} we $q$ replaced by $q+1$, we require that $\beta \leq \sfrac 13$. Thus, from this simple heuristic, we recover the Onsager-critical H\"older regularity exponent $\sfrac 13$. 

\subsection{Convex integration schemes for the Navier-Stokes equations}

Analogously to the case of the Euler equation, in order to construct the weak solutions of the Navier-Stokes equations, one proceeds via induction: for each $q \geq 0$ we assume we are given a solution $(v_q, \mathring R_q)$ to the Navier-Stokes-Reynolds system:
\begin{subequations}\label{e:Navier-Reynolds}
\begin{align}
\partial_t v_q + \div (v_q\otimes v_q)+\nabla p_q - \nu \Delta v_q &= \div\mathring R_q\\
\div v_q &= 0\, .
\end{align}
\end{subequations}
where the stress $\mathring R_q$ is assumed to be a trace-free symmetric matrix.

The main difficultly in implementing a convex integration scheme for the Navier-Stokes equations, compared to the Euler equations, is ensuring that the dissipative term $\nu \Delta w_{q+1}$ can be treated as an error in comparison to the quadratic term $\div(w_{q+1}\otimes w_{q+1})$. 

As in the case for Euler, the principal part of perturbation, $w_{q+1}^{(p)}$ is of the form \eqref{e:perturbation_form}, satisfying the low mode cancellation \eqref{e:low_freq_R_cancellation}. The principal difference to the Euler schemes is that the building blocks are chosen to be intermittent. In \cite{BV}, \emph{intermittent Beltrami waves} were introduced for this purpose, and in \cite{BCV18} the \emph{intermittent jets} were introduced (see~Section~\ref{s:intermittent}), which have a number of advantageous properties compared to intermittent Beltrami waves. 

In physical space, intermittency causes concentrations that results in the formation of intermittent peaks. In frequency space, intermittency smears frequencies. Analytically, intermittency has the effect of saturating Bernstein inequalities between different $L^p$ spaces \cite{CS2014}. In the context of convex integration, intermittency  reduces the strength of the linear dissipative term $\nu\Delta w_{q+1}$ in order to ensure that the nonlinear term $\div(w_{q+1}\otimes w_{q+1})$ dominates. 

For the case of intermittent jets, in order to parameterize the concentration, we introduce two parameters $r_{\|}$ and $r_{\perp}$ such that
\begin{equation}\label{e:ell_lambda_heur}
\frac{\lambda_q}{\lambda_{q+1}}\ll r_{\perp} \ll r_{\|}\ll 1\,.
\end{equation}
Each jet $W_{\xi}$ is defined to be supported on $\sim (r_{\perp}\lambda_{q+1})^3$ many cylinders of diameter $\sim \frac{1}{\lambda_{q+1}}$ and length $\sim \frac{r_{\|}}{r_{\perp}\lambda_{q+1}}$. In particular, the measure of the support of $W_{\xi}$  is $\sim r_{\|} r_{\perp}^2$. We note that such scalings are consistent with the jet $W_{\xi}$ being of frequency $\sim \lambda_{q+1}$. Finally, we normalize $W_{\xi}$ such that its $L^2$ norm is $\sim 1$. Hence by scaling arguments, one expects an estimate of the form
\begin{equation}\label{e:W_heuristic_est}
\norm{W_{\xi}}_{W^{N,p}}\les r_{\perp}^{\sfrac{2}{p}-1} r_{\|}^{\sfrac{1}{p}-\sfrac12}\lambda_{q+1}^N\,.
\end{equation}
In contrast to the Euler equations schemes, the inductive schemes for the Navier Stokes equations measure the perturbations $w_{q+1}$ and Reynolds stresses $R_q$, in $L^2$ and $L^1$ based spaces respectively. Assuming the bounds 
\begin{equation}\label{e:perturbation_scaling_NS}
\norm{w_{q+1}}_{L^{2}}\leq \lambda_{q+1}^{-\beta}\,,
\end{equation}
in order to achieve \eqref{e:low_freq_R_cancellation}, heuristically this requires that the the Reynolds stress $\mathring R_q$ obeys the bound
\begin{equation}\label{e:Reynolds_scaling_NS}
\norm{\mathring R_{q}}_{L^{1}}\leq \lambda_{q+1}^{-2\beta}\,.
\end{equation}
We note that \eqref{e:perturbation_scaling_NS} is suggestive that the final solution $v=\sum_{q} w_q$ converges in $H^{\beta}$; however in this review paper (as well as in the papers \cite{BV,BCV18}) we are not interested in obtaining the optimal regularity, we actually obtain a worse regularity exponent.

Using the $-1$ order linear operator $\mathcal R$, \eqref{e:W_heuristic_est} and \eqref{e:Reynolds_scaling} we are able to heuristically estimate the contribution of the dissipative term resulting from the principal perturbation $\nu\Delta w_{q+1}^{(p)}$ to the Reynolds stress error:
\begin{align*}
\norm{\mathcal R\left(\nu \Delta w_{q+1}^{(p)}\right)}_{L^1}&\lesssim \nu\lambda_{q+1}\sum_{\xi} \norm{a_{\xi}(\mathring R_q )}_{L^{\infty}}\norm{W_{\xi}}_{L^1} \notag \\
&\lesssim \nu r_{\perp} r_{\|}^{\sfrac12}\lambda_{q+1} \norm{R_q}_{L^{\infty}}^{\sfrac12}
\end{align*}
Thus, to ensure the error is small we will require
\[r_{\perp} r_{\|}^{\sfrac12}\ll \lambda_{q+1}\,.\]
This condition, together with the condition \eqref{e:ell_lambda_heur}, rules out geometric growth of the frequency $\lambda_q$. Indeed for the purpose of proving non-uniqueness of the Navier-Stokes equations let $\lambda_q$ be of the form
\[\lambda_q=a^{(b^q)}\quad\mbox{for}\quad a,b\gg 1\,.\]

Now consider the estimate \eqref{e:perturbation_scaling_NS}. Na\"ively estimating, the principle perturbation, we have
\begin{align*}
\norm{w_{q+1}^{(p)}}_{L^2}&\les \sum_{\xi} \norm{a_{\xi}(\mathring R_q )}_{L^{\infty}}\norm{W_{\xi}}_{L^2} \les  \norm{\mathring R_q }_{L^{\infty}}^{\sfrac12}\,.
\end{align*}
We do not however inductively propagate good estimates on the $L^{\infty}$ norm of $\mathring R_q$ and as such, the above na\"ive estimate is not suitable in order to obtain \eqref{e:perturbation_scaling_NS}. To obtain a better estimate, we will utilize the following observation: given a function $f$ with frequency contained in a ball of radius $\kappa$ and a $\lambda^{-1}$-periodic function $g$, if $\lambda \gg \kappa$ then
\begin{align}\label{e:tricky_tricky}
\norm{f  g}_{L^p}\lesssim \norm{f}_{L^p}\norm{g}_{L^p}\,.
\end{align}
Hence using that $R_q$ is of frequency roughly $\lambda_q$ we obtain
\begin{align*}
\norm{w_{q+1}^{(p)}}_{L^2}&\les \sum_{\xi} \norm{a_{\xi}(\mathring R_q )}_{L^{2}}\norm{W_{\xi}}_{L^2}\\
&\les  \norm{\mathring R_q }_{L^{1}}^{\sfrac12}\\
&\les  \lambda_{q+1}^{-\beta}\,,
\end{align*}
where we have used \eqref{e:Reynolds_scaling_NS}.

In comparison to Beltrami waves, or Mikado waves used for the Euler constructions, the intermittent building blocks used in \cite{BV,BCV18} introduce addition difficulties in handling the resulting oscillation error. For the intermittent jets of \cite{BCV18} we have
\begin{equation}\label{e:cancellation}
\div \left(w_{q+1}^{(p)}\otimes w_{q+1}^{(p)}+{\mathring { R}}_{q}\right)\sim
\sum_{\xi}2 a_\xi^2 W_{\xi}\cdot \nabla W_{\xi}  + \mbox{(high frequency error)}
\end{equation}
Similar to how the Nash error for the Euler equations was dealt with, the high frequency error experiences a gain when one inverts the divergence equation. In order to take care of the main term in \eqref{e:cancellation}, the intermittent jets are carefully designed (cf.~\eqref{eq:useful:2}) to oscillate in time such that the term can be written as a temporal derivative:
\[\sum_{\xi}2W_{\xi}\cdot \nabla W_{\xi}= \frac{1}{\mu} \partial_t \left(\sum_\xi   |W_{\xi}|^2 \xi \right)\]
for some large parameter $\mu$. This error can absorbed by introducing a temporal corrector $w_{q+1}^{(t)}$
\[w_{q+1}^{(t)}:= -  \frac{1}{\mu} \mathbb P_H\mathbb P_{\neq 0} \left(\sum_\xi a_{ \xi }^2 |W_{ \xi }|^2 \xi\right)\,,\]
where ${\mathbb P}_{H}$ is the Helmholtz projection, and ${\mathbb P}_{\neq 0}$ is the projection onto functions with mean zero. Thus pairing the oscillation error with the time derivative of the temporal corrector, we obtain
\[\div \left(w_{q+1}^{(p)}\otimes w_{q+1}^{(p)}+{\mathring { R}}_{q}\right)+\partial_t w_{q+1}^{(t)}\sim ~ \mbox{(pressure gradient) }+\mbox{ (high frequency error)}\,.
\]
Finally, analogous to the Euler case, we define a divergence corrector $w_{q+1}^{(c)}$ to corrector for the fact that $w_{q+1}^{(p)}$ is not, as defined, divergence free. The perturbation $w_{q+1}$ is then defined to be
\[w_{q+1}:=w_{q+1}^{(p)}+w_{q+1}^{(t)}+w_{q+1}^{(c)}\,.\]
An important point to keep in mind is that the temporal oscillation in the definition of the intermittent jets will introduce an error arising from the term $\partial_t w_{q+1}^{(p)}$ which is proportional  to $\mu$. The oscillation error is inversely proportional to $\mu$, and thus $\mu$ will be required to be chosen carefully to optimize the two errors.

More recently, the intermittent convex integration construction introduced in~\cite{BV}, combined with additional new ideas, has been successfully applied in related contexts. Using {\em intermittent Mikado flows}, Modena and Sz\'ekeyhidi Jr.~and have adapted these methods to establish the existence of non-renormalized solutions to the  transport and continuity equations with Sobolev vector fields~\cite{ModenaSZ17,ModenaSZ18}.
In~\cite{Dai18}, Dai demonstrated that these  methods can be adapted to prove non-uniqueness of Leray-Hopf weak solutions for the 3D Hall-MHD system. T.~Luo and Titi~\cite{LuoTiti18} demonstrated that these methods are applicable also to the fractional Navier-Stokes equations with dissipation $(-\Delta)^\alpha$, and $\alpha < \sfrac 54$ (the Lions criticality threshold~\cite{Lions59}).  X.~Luo~\cite{Luo18} demonstrated the existence of non-trivial stationary solutions to the 4D Navier-Stokes equations. The extra dimension allowed Luo to avoid adding temporal oscillations to the intermittent building block used in the construction (compare this to the oscillations introduced in Section~\ref{s:intermittent} and parametrized by $\mu$). Very recently, Cheskidov and X.~Luo~\cite{CheskidovLuo19} have further improved this construction by introducing new building blocks called {\em viscous eddies}, which allowed them to treat the 3D stationary case.

\section{Euler: the existence of wild continuous weak solutions}
\label{sec:Euler:C0}

We consider zero mean weak solutions of the  the Euler equations~\eqref{eq:Euler} (cf.~Definition~\ref{eq:def:weak:Euler}).  
In~\cite[Theorem 1.1]{DLSZ13}, De Lellis and Sz\'ekelyhidi gave the first proof for the existence of a $C^0_{x,t}$ weak solution of the 3D Euler equations which is non-conservative. The main result of this work is as follows:

\begin{theorem}[Theorem 1.1,~\cite{DLSZ13}]
\label{thm:Euler:C0:DLSZ}
Assume $e\colon [0,1] \to (0,\infty)$ is a  smooth function. Then there is a continuous vector field $v\colon \T^3 \times [0,1] \to \R^3$  and a continuous scalar field $p\colon \T^3 \times [0,1] \to \R$ which solve the incompressible Euler equations \eqref{eq:Euler} in the sense of distributions, and such that
\begin{align*}
e(t) = \int_{\T^3} |v(x,t)|^2 dx 
\end{align*} for all $t\in [0,1]$.
\end{theorem}
In order to simplify the presentation, we only give here the details of an Euler $C^{0+}$ convex integration scheme, without attempting to attain a given energy profile (this would require adding one  more inductive estimate to the the list in \eqref{eq:inductive:Euler} below, see equation (7) in \cite{DLSZ13}).  The main result of this section is the existence of a continuous weak solution which is not conservative:
\begin{theorem}
\label{thm:Euler:C0}
There exists $\beta>0$ such that the following holds. There a weak solution $v \in C^{0}([0,1];C^\beta)$ of the Euler equations \eqref{eq:Euler} such that $\norm{v(\cdot,1)}_{L^2} \geq 2 \norm{v(\cdot,0)}_{L^2}$.
\end{theorem}
In Section~\ref{sec:Euler:C1/3} we present the necessary ideas required to obtain a $C^{\sfrac 13 - }$ solution, and discuss the necessary ingredients required to fix the energy profile.

\subsection{Inductive estimates and iteration proposition}

Below we assume $(v_q,\RR_q)$ is a given solution of the Euler-Reynolds system \eqref{e:euler_reynolds}.
We consider an increasing sequence $\{ \lambda_q\}_{q\in \N} \in 2\pi \N$ which diverges to $\infty$, and $\{\delta_q\}_{q\in \N} \in (0,1)$  a  sequence which is decreasing towards $0$, and such that $\delta_q^{\sfrac 12} \lambda_q$ is monotone increasing. It is convenient to specify these sequences, modulo some free parameters. For this purpose, we introduce $a\in \N$,   $\beta \in (0,1)$, and let 
\begin{align*}
\lambda_q &= 2\pi a^{(2^q)} \\
\delta_q &= \lambda_q^{-2\beta} \,.
\end{align*}
The parameter $\beta$ will be chosen sufficiently small, as specified in Proposition~\ref{prop:iteration:Euler}. The parameter $a$ will be  chosen as a sufficiently large multiple of a geometric constant $n_* \in  \N$ (which is fixed in Proposition~\ref{p:split}).

By induction on $q$ we will assume that the following bounds hold for the solution $(v_q,\RR_q)$ of \eqref{e:euler_reynolds}:
\begin{subequations}
\label{eq:inductive:Euler}
\begin{align}
\norm{v_q}_{C^0} &\leq 1  - \delta_{q}^{\sfrac12}
\label{eq:vq:C0:inductive}
\\
\norm{v_q}_{C^{1}_{x,t}} &\leq \delta_{q}^{\sfrac12}\lambda_q
\label{eq:vq:C1:inductive}
\\
\norm{\RR_{q}}_{C^0} &\leq c_R \delta_{q+1}
\label{eq:Rq:C0:inductive}
\end{align}
\end{subequations}
where $c_R >0$ is a sufficiently small universal constant (see estimates \eqref{eq:w:q+1:p:L:infinity} and and \eqref{eq:BBB:2} below).
Condition~\eqref{eq:vq:C0:inductive} is   not necessary  for a $C^0$-convex integration scheme, but it is convenient to propagate it.

The following proposition summarizes the iteration procedure which goes from level $q$ to $q+1$. 
\begin{proposition}[Main iteration]
\label{prop:iteration:Euler}
There exists a sufficiently small parameter $\beta   \in (0,1)$, such that the following holds.  There exists a sufficiently large constant $a_0 = a_0(c_R,\beta)$ such that for any $a\geq a_0$ which is a multiple of the geometric constant $n_*$, there exist functions $(v_{q+1},\RR_{q+1})$ which solve \eqref{e:euler_reynolds} and obey  \eqref{eq:inductive:Euler}  at level $q+1$, such that 
\begin{align}
\norm{v_{q+1} - v_q }_{C^0} \leq   \delta_{q+1}^{\sfrac 12} 
\label{eq:increment:Linfty}
\end{align}
holds.
\end{proposition}

\begin{remark} \label{rem:b:beta:precise} 
Inspecting the proof of Proposition~\ref{prop:iteration:Euler},we remark that it is sufficient to take  $\beta = \sfrac{1}{100}$.
\end{remark}

\subsection{Proof of Theorem~\ref{thm:Euler:C0}}
Fix the parameter  $ c>0$ as in Lemma~\ref{l:linear_algebra} below, and the parameters $\beta$ and  $a_0$ from Proposition~\ref{prop:iteration:Euler}. By possibly enlarging the value of $a\geq a_0$, we may ensure that $\delta_{0} \leq \sfrac 14$.

We define an incompressible, zero mean vector field $v_0$ by
\[
v_0(x,t) = \frac{t}{(2\pi)^{\sfrac 32}} ( \sin(\lambda_0^{\sfrac 12} x_3), 0, 0)  \, .
\]
Note that by construction we have 
$\sup_{t\in [0,1]} \norm{v_0(\cdot,t)}_{C^0} \leq \norm{v_0(\cdot,1)}_{C^0} =  (2\pi)^{-\sfrac 32}  \leq 1- \delta_0^{\sfrac 12}$,  so that \eqref{eq:vq:C0:inductive} is automatically satisfied. Moreover, $\norm{v_0}_{C^{1}_{x,t}} \leq \lambda_0^{\sfrac 12} \leq   \lambda_0 \delta_0^{\sfrac 12}=   \lambda_0^{1-\beta}$. This inequality holds because $\beta$ is strictly smaller than $\sfrac 12$, and $\lambda_0 = 2\pi a \geq 1$ may be chosen sufficiently large, depending on $\beta$. Thus, \eqref{eq:vq:C1:inductive} also holds at level $q=0$.

The vector field $v_0$ defined above is a shear flow, and thus $v_0 \cdot \nabla v_0 = 0$. Thus, it obeys \eqref{e:euler_reynolds} at $q=0$, with stress $\RR_0$ defined by 
\begin{align}
 \RR_0 =  \frac{1}{\lambda_0^{\sfrac 12} (2\pi)^{\sfrac 32}}  \left( {\begin{array}{ccc}
   0 & 0 & - \cos(\lambda_0^{\sfrac 12} x_3) \\
   0 & 0 & 0 \\
   - \cos(\lambda_0^{\sfrac 12} x_3) & 0 & 0   
  \end{array} } \right) \, . 
 \label{eq:R0:Euler:def}
\end{align}
Therefore, we have
\begin{align*}
 \norm{\RR_0}_{C^0} \leq  \lambda_0^{- \sfrac 12}  \leq c_R \delta_1 \, .
\end{align*}
The last inequality above uses that $\lambda_0^{\sfrac 12} \delta_1 = (2\pi)^{\sfrac 12 - 2\beta} a^{\sfrac 12  -4\beta} \geq a^{\sfrac 14} \geq c_R^{-1}$. This inequality holds because we may ensure $\beta  \leq \sfrac{1}{16}$ (see Remark~\ref{rem:b:beta:precise} above), and $a$ can be taken to be sufficiently large, in terms of the   universal constant $c_R$.  Thus, condition \eqref{eq:Rq:C0:inductive} is also obeyed for $q=0$. 
 
We may thus start the iteration Proposition~\ref{prop:iteration:Euler} with the pair $(v_0,\RR_0)$ and obtain a sequence of solutions $(v_q,\RR_q)$. By  \eqref{eq:inductive:Euler}, \eqref{eq:increment:Linfty} and interpolation we have that for any $\beta' \in (0,\beta)$, the following series is   summable
\begin{align*}
\sum_{q\geq 0} \norm{v_{q+1}-v_q}_{C^{\beta'}} 
\les \sum_{q\geq 0} \norm{v_{q+1}-v_q}_{C^0}^{1-\beta'}  \norm{v_{q+1}-v_q}_{C^1}^{\beta'} 
&\les    \sum_{q\geq 0} \delta_{q+1}^{\sfrac{1}{2} } \lambda_{q+1}^{\beta'}  \les   \sum_{q\geq 0} \lambda_{q+1}^{\beta' - \beta }  \les 1 
\end{align*}
where the implicit constant is universal. Thus, we may define a limiting function $v = \lim_{q\to \infty} v_q$ which lies in $C^0([0,1];C^{\beta'})$. Moreover, $v$ is a weak solution of the Euler equation \eqref{eq:Euler}, since by \eqref{eq:Rq:C0:inductive} we have that $\lim_{q\to \infty}\RR_q = 0$ in $C^0([0,1];C^0)$.  The regularity of the weak solution claimed in Theorem~\ref{thm:Euler:C0} then holds with $\beta$ replaced by $\beta'$.

It remains to show that $\norm{v(\cdot,1)}_{L^2} \geq 2 \norm{v(\cdot,0)}_{L^2}$. For this purpose note that since $2^{q+1} \geq 2 (q+1)$, we have
\begin{align*}
\norm{v - v_0}_{C^0} \leq \sum_{q\geq 0} \norm{v_{q+1} - v_q}_{C^0} 
\leq   \sum_{q\geq 0} \delta_{q+1}^{\sfrac 12} 
\leq   \sum_{q\geq 0} a^{-\beta (2^{q+1})} \leq   \sum_{q\geq 0} (a^{-\beta b})^{q+1}  = \frac{  a^{-2\beta }}{1-a^{-2\beta}} \leq \frac{1}{6 (2\pi)^{\sfrac 32}}
\end{align*}
once we choose $a$ sufficiently large. Since by construction $\norm{v_0(\cdot,0)}_{L^2} = 0$, and $\norm{v_0(\cdot,1)}_{L^2} = \sfrac{1}{\sqrt{2}}$, we obtain that 
\begin{align*}
2 \norm{v(\cdot,0)}_{L^2} &\leq 2 \norm{v_0(\cdot,0)}_{L^2} + 2 (2\pi)^{\sfrac 32} \norm{v(\cdot,0) - v_0(\cdot,0)}_{C^0}  \notag\\
&\leq   \frac13 \leq \frac{1}{\sqrt{2}} - \frac{1}{6} \leq \norm{v_0(\cdot,1)}_{L^2} - \norm{v(\cdot,1) - v_0(\cdot,1)}_{L^2} \leq \norm{v(\cdot,1)}_{L^2}
\end{align*}
holds. This concludes the proof of Theorem~\ref{thm:Euler:C0}.

\subsection{Mollification}
\label{sec:mollify}
In order to avoid a loss of derivative, we replace $v_q$ by a mollified velocity field $v_\ell$. For this purpose we choose a small parameter $\ell \in (0,1)$ which lies between $\lambda_{q}^{-1}$ and $\lambda_{q+1}^{-1} \approx \lambda_q^{-2}$ as 
\begin{align}
\ell  =  \lambda_{q}^{- \sfrac 32}.
\label{eq:ell:def}
\end{align}

Let $\{ \phi_{\eps} \}_{\eps>0}$ be a family of standard Friedrichs mollifiers (of compact support of radius $2$) on $\R^3$ (space), and $\{\varphi_{\eps}\}_{\eps>0}$ be a family of standard Friedrichs mollifiers (of compact support of width $2$) on $\R$ (time). We define a mollification of $v_q$ and $\RR_q$ in space and time, at length scale and time scale $\ell$   by
\begin{subequations}
\label{eq:v:R:ell:def}
\begin{align}
v_{\ell} = (v_q \ast_x \phi_{\ell}) \ast_t \varphi_{\ell}  \, ,\\
\RR_{\ell} = (\RR_q \ast_x \phi_{\ell}) \ast_t \varphi_{\ell} \,.
\end{align}
\end{subequations}
Then using \eqref{e:euler_reynolds} we obtain that $(v_\ell,\RR_\ell)$ obey
\begin{subequations}
\label{e:euler_reynolds_ell}
\begin{align}
\partial_t  v_\ell + \div(v_\ell \otimes v_\ell) + \nabla p_\ell 
&= \div \Big(  \RR_\ell +   R_{\rm commutator} \Big)\, ,  \\
\div v_\ell &= 0 \, ,
\end{align}
\end{subequations}
where  traceless symmetric commutator stress $  R_{\rm commutator}$ is given by
\begin{align}
  R_{\rm commutator} &= (v_\ell \mathring \otimes v_\ell) - ((v_q \mathring \otimes v_q)\ast_x \phi_{\ell}) \ast_t \varphi_{\ell} \, .\label{e:euler_reynolds_ell:a} 
\end{align}
Using a standard mollification estimate we obtain
\begin{align}
\norm{R_{\rm commutator}}_{C^0} \les \ell  \norm{v_q}_{C^1_{x,t}} \norm{v_q}_{C^0_{x,t}} \les \ell  \delta_{q}^{\sfrac 12} \lambda_q = \lambda_q^{-\beta -\sfrac 12}  \ll \delta_{q+2} \, .
 \label{eq:Rc:bound}
\end{align}
In the last estimate above we have used that $\beta$ may be chosen to be sufficiently small and $a$   sufficiently large. 
Moreover, with $\ell$ small as above we have
\begin{align}
 \norm{v_q - v_\ell}_{C^0}& \les  \ell\norm{v_q}_{C^1}\les \ell \lambda_q \delta_q^{\sfrac 12}  \ll  \delta_{q+1}^{\sfrac 12}
 \label{eq:V_ell_est}
 \end{align}
while for $N\geq 1$ we obtain from standard mollification estimates that
\begin{align}
  \norm{v_\ell}_{C^N_{x,t}} \les \ell^{-N+1} \norm{v_q}_{C^1_{x,t}} \les \lambda_q \delta_q^{\sfrac 12} \ell^{-N+1} \les  \ell^{-N}.
   \label{eq:v:ell:CN}
\end{align}
For $N=0$ we simply use that the mollifier has mass $1$ to obtain
\begin{align}
\norm{v_\ell}_{C^0} \leq \norm{v_q}_{C^0}  \leq 1- \delta_q^{\sfrac 12}.
   \label{eq:v:ell:C0}
\end{align}

\subsection{Beltrami waves}
\label{sec:Beltrami}
Given $\xi\in {\mathbb S}^2\cap {\mathbb Q}^3$ let $A_{\xi}\in {\mathbb S}^2 \cap {\mathbb Q}^3$ obey
$$
A_{\xi}\cdot {\xi}=0, \quad A_{-{\xi}}=A_{\xi} \,.
$$
We define the complex vector
$$
B_{\xi}= \tfrac{1}{\sqrt{2}} \left(A_{\xi}+i\xi\times A_{\xi}\right) \, .
$$
By construction, the vector $B_\xi$ has the properties 
$$
|B_\xi| = 1, \quad B_{\xi} \cdot \xi = 0, \quad i \xi \times B_\xi = B_\xi, \quad B_{-\xi} = \overline{B_\xi}\, .
$$
This implies that for any $\lambda \in   \Z$, such that $\lambda \xi \in   \Z^3$, the function
\begin{align}
W_{(\xi)}(x) := W_{\xi,\lambda}(x) :=  B_\xi e^{i \lambda \xi \cdot x}
\label{eq:Betrami:plane}
\end{align}
is $\T^3$ periodic, divergence free, and is an eigenfunction of the $\curl$ operator with eigenvalue $\lambda$. That is, $W_{(\xi)}$ is a complex Beltrami plane wave. The following lemma states a useful property for linear combinations of complex Beltrami plane waves.

\begin{proposition}[Proposition 3.1 in~\cite{DLSZ13}]
\label{p:Beltrami}
Let $\Lambda$ be a given finite subset of $\mathbb S^2\cap\mathbb Q^3$ such that $-\Lambda = \Lambda$, and let $\lambda\in  \mathbb Z$ be such that $\lambda\Lambda\subset  \mathbb Z^3$. Then for any choice of coefficients $a_{\xi}\in\CC$ with $\overline{a}_{\xi} = a_{-\xi}$ the vector field
\begin{equation}\label{e:Beltrami}
W(x) = \sum_{\xi\in \Lambda}
a_{\xi} B_{\xi}e^{i\lambda  \xi\cdot x}
\end{equation}
is a real-valued, divergence-free Beltrami vector field $\curl W = \lambda W $,
and thus it is a stationary solution of the Euler equations
\begin{equation}\label{e:Bequation}
\div (W\otimes W)=\nabla\frac{|W|^2}{2}.
\end{equation}
Furthermore, since $B_\xi \otimes B_{-\xi} + B_{-\xi} \otimes B_{\xi} = 2 \Re(B_\xi \otimes B_{-\xi}) = \Id - \xi \otimes \xi$, we have
\begin{equation}\label{e:av_of_Bel}
\fint_{\T^3} W\otimes W\,dx = \frac{1}{2} \sum_{\xi\in\Lambda} |a_{\xi}|^2 \left( \Id - \xi \otimes\xi\right)\, .  
\end{equation}
\end{proposition}

\begin{figure}[h!]
\begin{center}
\includegraphics[width=0.5\textwidth]{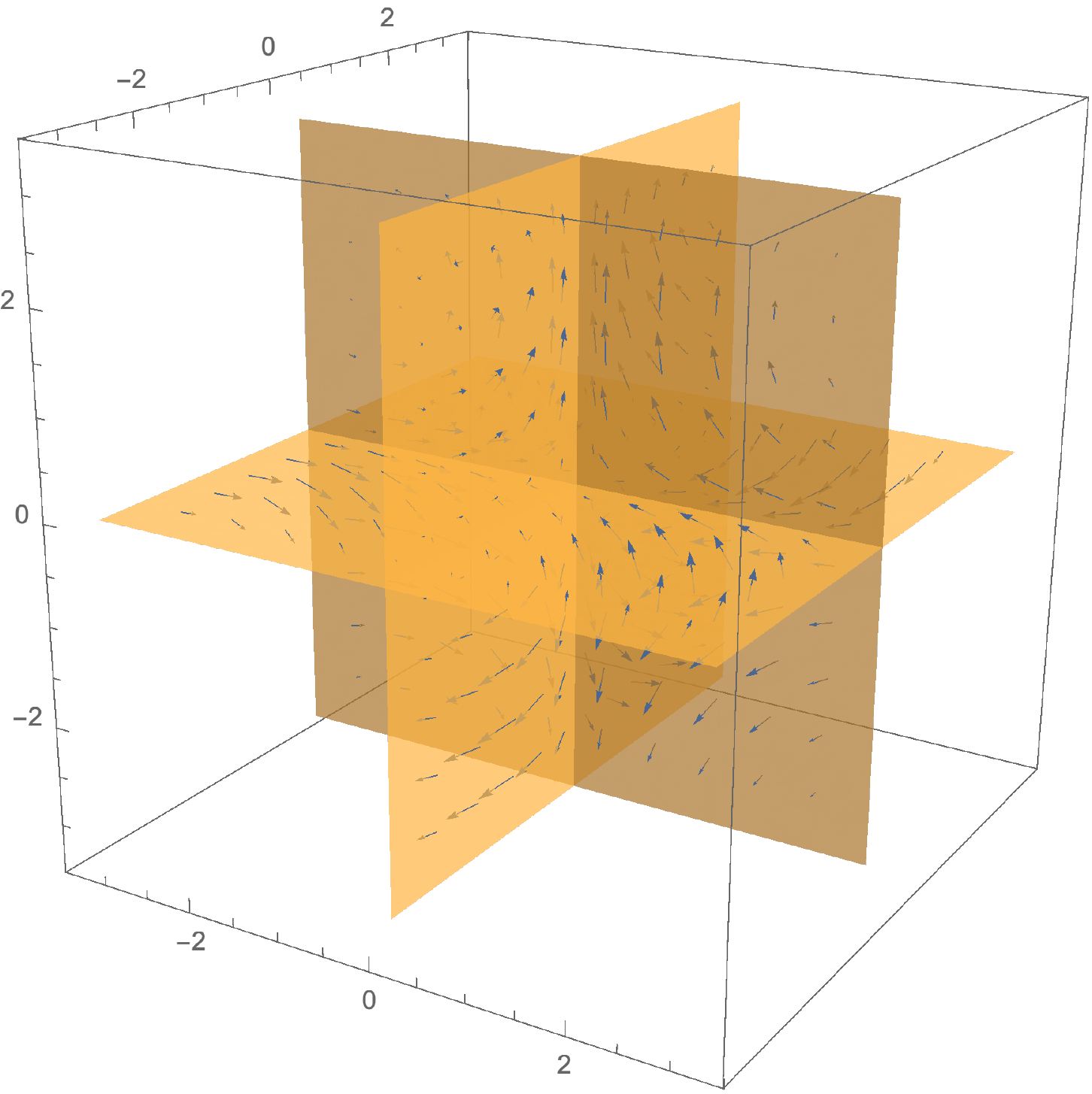}
\end{center}
\caption{{\small Example of a Beltrami flow $W(x)$ as defined in \eqref{e:Beltrami}.}}
\end{figure}

\begin{proposition}[Lemma 3.2 in~\cite{DLSZ13}]
\label{p:split}
There exists a sufficiently small $c_*>0$ with the following property. Let $B_{c_*}(\Id)$ denote the closed ball of symmetric $3\times 3$ matrices, centered at $\Id$, of radius $c_*$. Then, there exist pairwise disjoint subsets 
$$
\Lambda_{\alpha}\subset\mathbb S^2 \cap {\mathbb Q}^3\qquad \alpha\in \{0, 1\} \, ,
$$
and smooth positive functions  \[
\gamma_{ {\xi}}^{ {(\alpha)}}\in C^{\infty}\left(B_{c} (\Id)\right) \qquad \alpha\in \{0, 1\}, \, \xi\in\Lambda_{\alpha} \,,
\]
such that the following hold. For every $\xi\in \Lambda_{\alpha}$ we have $-\xi\in \Lambda_{\alpha}$ and $\gamma_{\xi}^{ {(\alpha)}} = \gamma_{-\xi}^{{(\alpha)}}$.  For each $R\in B_{c_*} (\Id)$ we have the identity
\begin{equation}\label{e:split}
R = \frac{1}{2} \sum_{\xi\in\Lambda_{\alpha}} \left(\gamma_{\xi}^{ {(\alpha)}}(R)\right)^2 \left(\Id - \xi\otimes \xi\right) .
\end{equation}
We label by $n_*$  the smallest natural number such that $n_* \Lambda_\alpha \subset \Z^3$ for all $\alpha \in \{1,2\}$.
\end{proposition}
It is sufficient to consider  index sets $\Lambda_0$ and $\Lambda_1$ in Proposition~\ref{p:split} to have $12$ elements.  Moreover, by abuse of notation, for $j \in \Z$ we denote $\Lambda_j = \Lambda_{j \, {\rm mod} \, 2}$. Also, it is convenient to denote by $M$ a {\em geometric constant} such that
\begin{align}
\sum_{\xi \in \Lambda_\alpha} \norm{\gamma_\xi^{(\alpha)}}_{C^1(B_{c_*}(\Id))} \leq M
 \label{eq:M:universal}
\end{align}
holds for  $\alpha \in \{0,1\}$ and $\xi \in \Lambda_\alpha$. This parameter  is universal.

\subsection{The perturbation}

\subsubsection{Flow maps and time cutoffs}
\label{sec:flow:maps:1}
In order to have an acceptable {\em transport error}, the perturbation $w_{q+1}$ needs to be transported by the flow of the vector field $\partial_t + v_\ell \cdot \nabla$, at least to leading order. The natural way to achieve this, is to replace the linear phase $\xi\cdot x$ in the definition of the Beltrami wave $W_{\xi,\lambda}$, with the nonlinear phase $\xi \cdot \Phi(x,t)$, where $\Phi$ is transported by the aforementioned vector field.

We subdivide $[0,1]$ into time intervals of size $\ell$, and solve transport equations on these intervals.\footnote{Standard ODE arguments show  that the time $\tau$ such that the flow of $\partial_t + v_\ell \cdot \nabla$ remains close to its initial datum on $[-\tau,\tau]$, should obey $\tau \norm{v_{\ell}}_{C^1} \ll 1$. This is the same as the CFL condition~\cite{CFL28}. Since in this exposition we do not aim for the most optimized possible convex-integration scheme, instead of introducing a new parameter  for the CFL-time, which is then to be optimized later, we work with the already available parameter $\ell$. Indeed, \eqref{eq:ell:def} shows that $\ell \lambda_q \delta_q^{\sfrac 12}= \lambda_{q}^{-\sfrac 12 - \beta} \ll 1$, which in view of \eqref{eq:v:ell:CN} shows that $\ell \norm{v_{\ell}}_{C^1} \ll 1$ holds, as desired.} For $j \in \{0,\ldots,\lceil \ell^{-1} \rceil\}$,\footnote{Here we use $\lceil x \rceil$ to denote the smallest integer $n\geq x$.} we define the map $\Phi_j: \R^3\times [0,1]\to \R^3$ as the $\T^3$ periodic solution of 
\begin{subequations}\label{e:Phi-transport}
\begin{align}
\left(\partial_t  + v_{\ell}\cdot  \nabla\right) \Phi_j =0 \, \\ 
\Phi_j (x, j\ell)=x .
\end{align}
\end{subequations}
This map $\Phi_j$ obeys the expected estimates 
\begin{subequations}
\label{eq:Phi:j:bounds}
\begin{align}
\sup_{t\in [(j-1) \ell, (j+1)\ell]} \norm{\nabla \Phi_j(t) - \Id}_{C^0} 
&\les \ell \norm{v_{\ell}}_{C^1} \les \ell \lambda_q \delta_q^{\sfrac 12} \ll 1\,
\label{eq:Phi:near:Id} \\
\sup_{t\in [(j-1) \ell, (j+1)\ell]} \norm{\nabla \Phi_j(t)}_{\dot{C}^{1}_{x,t}}
&\les \lambda_q \delta_q^{\sfrac 12}
\label{eq:Phi:C2} \\
\sup_{t\in [(j-1) \ell, (j+1)\ell]} \norm{\nabla \Phi_j(t)}_{C^{n}} 
&\les \ell^{1-n} \lambda_q \delta_q^{\sfrac 12} \les \ell^{-n}
\label{eq:Phi:CN}
\end{align}
\end{subequations}
for $n\geq 1$,  which are a consequence of the Gr\"onwall inequality for derivatives of \eqref{e:Phi-transport} (see, e.g.~\cite[Proposition D.1]{BDLISZ15}).  We also let $\chi$ be a non-negative  bump function supported in $(-1,1)$ which is identically $1$ on $(-\sfrac 14, \sfrac 14)$, and such that the shifted bump functions
\begin{align}
\chi_j(t) :=\chi\Bigl(\ell^{-1} t- j \Bigr),\label{e:chi_l}
\end{align}
form a partition of unity in time once they are squared
\begin{align}
\sum_{j} \chi_j^2(t) = 1
\label{eq:chi:partition}
\end{align}
for all $t\in [0,1]$. Note that the sum over $j$ is finite, $j \in \{0, 1,\ldots,\lceil \ell^{-1} \rceil\}$, and that at each time $t$ at most two cutoffs are nontrivial.

\subsubsection{Amplitudes}
In view of Proposition~\ref{p:split}, we introduce the amplitude functions
\begin{align}
a_{(\xi)}(x,t):=a_{q+1,j,\xi}(x,t)&:= c_R^{\sfrac 14} \delta_{q+1}^{\sfrac12} \, \chi_j(t) \,  \gamma_{\xi} \left(\Id - \frac{\RR_{\ell}(x,t)}{c_R^{\sfrac 12} \delta_{q+1}}\right).\label{e:a_xi}
\end{align}
The division of $\RR_\ell$ by the parameter $c_R^{\sfrac 12}  \delta_{q+1}$ ensures via \eqref{eq:Rq:C0:inductive} and the fact that the mollifier has mass $1$, that 
\begin{align*}
\norm{ \frac{\RR_{\ell}(x,t)}{c_R^{\sfrac 12} \delta_{q+1}} }_{C^0}\leq c_R^{\sfrac 12} \leq c_*.
\end{align*}
Therefore,  $\Id - \RR_\ell c_R^{-\sfrac 12} \delta_{q+1}^{-1}$ lies in the domain of the functions $\gamma_{(\xi)}$ and we deduce from \eqref{e:split} and \eqref{eq:chi:partition} that
\begin{align}
c_R^{\sfrac 12} \delta_{q+1} \Id - \RR_\ell = \frac{1}{2} \sum_j \sum_{\xi\in\Lambda_j} a_{(\xi)}^2 \left(\Id - \xi\otimes \xi\right) 
\label{eq:amplitudes:are:good}
\end{align}
holds pointwise, for any $\alpha \in \{0,1\}$. For a given $j$, we write $\Lambda_j = \Lambda_0$  if $j$ is even and $\Lambda_j = \Lambda_1$ is $j$ is odd. This justifies definition~\eqref{e:a_xi} of the amplitudes. 

\subsubsection{Principal part of the corrector}
\label{eq:compose:with:flow:map:and:get:a:coffee}
Using the notation from \eqref{eq:Betrami:plane}, \eqref{e:a_xi}, and \eqref{e:Phi-transport} we let
\begin{align}
w_{(\xi)}(x,t) := w_{q+1,j,\xi}(x,t)& = a_{q+1,j,\xi}(x,t)\,W_{\xi,\lambda_{q+1}} (\Phi_j(x,t)) =  a_{q+1,j,\xi}(x,t)\,B_\xi e^{i \lambda_{q+1} \xi \cdot \Phi_j(x,t)} .\label{e:w_xi}
\end{align}
and define the \emph{principal part} $w_{q+1}^{(p)} $of the perturbation $w_{q+1}$ as
\begin{align}\label{e:def_wo}
w_{q+1}^{(p)} (x,t) := \sum_{j} \sum_{\xi \in \Lambda_{j}}  w_{(\xi)} (x,t) \, .
\end{align}
From \eqref{eq:M:universal}, \eqref{e:a_xi}, and the fact that $\chi_j^2$ form a partition of unity, it follows that 
\begin{align}
 \norm{w_{q+1}^{(p)}}_{C^0} \leq M  c_R^{\sfrac 14} \delta_{q+1}^{\sfrac 12} \leq \frac{\delta_{q+1}^{\sfrac 12}}{2}
 \label{eq:w:q+1:p:L:infinity}
\end{align}
where we have used that $c_R$ may be taken sufficiently small, in terms of the universal constant $M$.

\subsubsection{Incompressibility correction}
In order to define the incompressibility correction it is useful to introduce the scalar \emph{phase} function
\begin{equation}\label{e:phi_xi}
\phi_{(\xi)}(x,t):= \phi_{q+1,j,\xi}(x,t):=e^{i\lambda_{q+1}\xi\cdot (\Phi_j(x,t)-x)} \,.
\end{equation}
In view of \eqref{eq:Phi:near:Id}, we think of $\phi_{(\xi)}$ as oscillating at a frequency $\ll \lambda_{q+1}$. 
Also, with this notation, \eqref{e:w_xi} reads as
\[
w_{(\xi)}(x,t) = a_{(\xi)}(x,t) \phi_{(\xi)}(x,t) W_{(\xi)}(x) = a_{(\xi)}(x,t) \phi_{(\xi)}(x,t) B_\xi e^{i \lambda_{q+1}\xi\cdot x} \, ,
\]
and the term $W_{(\xi)}$ oscillates the fastest (at frequency $\lambda_{q+1}$). We will add a corrector to $w_{(\xi)}$ such that the resulting function is a perfect $\curl$, making it thus divergence free. For this purpose recall that $\curl W_{(\xi)} = \lambda_{q+1} W_{(\xi)}$ and therefore, since $a_{(\xi)}$ and $\phi_{(\xi)}$ are scalar functions, we have
\begin{align*}
a_{(\xi)} \phi_{(\xi)} W_{(\xi)} = \frac{1}{\lambda_{q+1}} \curl\left( a_{(\xi)} \phi_{(\xi)} W_{(\xi)} \right) - \frac{1}{\lambda_{q+1}} \nabla \left(a_{(\xi)} \phi_{(\xi)}\right) \times W_{(\xi)}.
\end{align*}
We therefore define 
\begin{align}
 w_{(\xi)}^{(c)}(x,t) &:= \frac{1}{\lambda_{q+1}} \nabla\left(a_{(\xi)}\phi_{{(\xi)}} \right) \times B_\xi e^{i\lambda_{q+1}\xi \cdot x} \notag\\
 &= \frac{1}{\lambda_{q+1}}\left(\nabla a_{(\xi)}  + a_{(\xi)} i \lambda_{q+1} (\nabla \Phi_j(x,t) - \Id) \xi \right) \times B_\xi e^{i\lambda_{q+1}\xi \cdot \Phi_j(x,t)} \notag\\
 &= \left(\frac{\nabla a_{(\xi)}}{\lambda_{q+1}}  + i \, a_{(\xi)}    (\nabla \Phi_j(x,t) - \Id) \xi \right) \times W_{(\xi)}(\Phi_j(x,t)).
\label{eq:w:c:small}
\end{align}

The {\em incompressibility correction} $w_{q+1}^{(c)}$ is then defined as
 \begin{align}
w_{q+1}^{(c)}(x,t)
&:= \sum_{j} \sum_{\xi \in \Lambda_{j}}   w_{(\xi)}^{(c)}(x,t)\, ,
\label{e:corrector}
\end{align}
so that setting
\[w_{q+1}:=w_{q+1}^{(p)}+w_{q+1}^{(c)}\]
we obtain from the above computations that 
\begin{equation}\label{e:w_compact_form}
w_{q+1} = \frac{1}{\lambda_{q+1}} \sum_{j} \sum_{\xi \in \Lambda_{j}}  \curl \left(a_{(\xi)}\,  W_{(\xi)} \circ \Phi_j  \right)
\end{equation}
and so clearly $w_{q+1}$ is divergence and mean free.

Note that by \eqref{eq:Rq:C0:inductive}, \eqref{eq:chi:partition}, standard mollification estimates, \eqref{eq:M:universal} and \eqref{eq:Phi:near:Id} imply that 
\begin{align}
\norm{w_{q+1}^{(c)}}_{C^0} 
&\leq 2 \sup_{j} \sum_{\xi \in \Lambda_{j}}  \frac{\norm{\nabla a_{(\xi)}}_{C^0}}{\lambda_{q+1}} + \norm{a_{(\xi)}}_{C^0} \norm{\nabla \Phi_j - \Id}_{C^0(\supp \chi_j)} \notag\\
&\les  \delta_{q+1}^{\sfrac 12} \left( \frac{1}{\ell \lambda_{q+1}} + \ell \lambda_q \delta_{q}^{\sfrac 12}\right)  
\ll \delta_{q+1}^{\sfrac 12} \, .
 \label{eq:w:q+1:c:L:infinity}
\end{align}

\subsubsection{The velocity inductive estimates}
We define the velocity field at level $q+1$ as
\begin{align}
v_{q+1}:= v_\ell + w_{q+1}.
\label{eq:v_q+1:def}
\end{align}
At this stage we verify that \eqref{eq:vq:C0:inductive} and \eqref{eq:vq:C1:inductive} hold at level $q+1$.

First, we note that \eqref{eq:w:q+1:p:L:infinity} and \eqref{eq:w:q+1:c:L:infinity} give that 
\begin{align*}
\norm{w_{q+1}}_{C^0} \leq \frac{3}{4} \delta_{q+1}^{\sfrac 12}
\end{align*}
which combined with \eqref{eq:V_ell_est} gives the proof of \eqref{eq:increment:Linfty}.  Moreover, combining the above estimate with \eqref{eq:v:ell:C0} yields 
\begin{align*}
\norm{v_{q+1}}_{C^0} \leq 1 - \delta_q^{\sfrac 12} +   \delta_{q+1}^{\sfrac 12} \leq 1 - \delta_{q+1}^{\sfrac 12}
\end{align*}
since $4 \delta_{q+1} \leq \delta_q$ holds upon choosing $a$ sufficiently large. This proves \eqref{eq:vq:C0:inductive}.

A short calculation shows that upon applying a spatial or a temporal derivative to \eqref{e:def_wo} and \eqref{e:corrector}, similarly to \eqref{eq:w:q+1:p:L:infinity} and \eqref{eq:w:q+1:c:L:infinity} we obtain that for some $q$-independent constant $C$ we have
\begin{align}
\norm{w_{q+1}}_{C^1_{x,t}} &\leq c_R^{\sfrac 14} M \lambda_{q+1} \delta_{q+1}^{\sfrac 12} + C \ell^{-1} \delta_{q+1}^{\sfrac 12} + C \ell^{-2} \lambda_{q+1}^{-1} \delta_{q+1}^{\sfrac 12} + C  \lambda_q \delta_q^{\sfrac 12} \delta_{q+1}^{\sfrac 12}  \leq   \lambda_{q+1} \delta_{q+1}^{\sfrac 12} \,.
\label{eq:BBB:2}
\end{align}
In the above inequality above we have used that  $\ell = \lambda_q^{-\sfrac 32} \approx \lambda_{q+1}^{- \sfrac 34}$, have taken $c_R$ sufficiently small in terms of $M$, $\beta$ sufficiently small  and  $a$  sufficiently large.
This proves \eqref{eq:vq:C1:inductive} at level $q+1$.

\subsection{Reynolds Stress}

\subsubsection{Inverse divergence operator and stationary phase bounds}
\label{sec:RSZ:bounds}
We recall~\cite[Definition 4.2]{DLSZ13} the operator $\RSZ$ which acts on  vector fields $v$ with $\int_{\T^3} v dx = 0$ as
\begin{align}
(\RSZ v)^{k\ell} = (\partial_k \Delta^{-1} v^{\ell} + \partial_\ell \Delta^{-1} v^k)  - \frac{1}{2} \left( \delta_{k \ell} + \partial_k \partial_\ell \Delta^{-1}\right)\div \Delta^{-1} v  
\label{eq:RSZ}
\end{align} 
for $k,\ell \in \{1,2,3\}$. The above inverse divergence operator has the property that  $\RSZ v(x)$ is a symmetric trace-free matrix for each $x \in \T^3$, and $\RSZ$ is an right inverse of the $\div$ operator, i.e.\ $\div (\RSZ v) = v$. When $v$ does not obey $\int_{\T^3} v dx = 0$, we overload notation and denote $\RSZ v := \RSZ( v -\int_{\T^3} v dx)$. Note that  $\nabla \RSZ$ is a Calder\'on-Zygmund operator.  

The following lemma makes rigorous the fact that $\RSZ$ obeys the same elliptic regularity estimates as $|\nabla|^{-1}$. We recall the following {\em stationary phase lemma} (see for example~\cite[Lemma 2.2]{DSZ17}), adapted to our setting. 
\begin{lemma}
\label{lem:stationary:phase}
Let $\lambda \xi \in \Z^3$, $\alpha \in (0,1)$, and $m\geq 1$. Assume that $a \in C^{m,\alpha}(\T^3)$ and $\Phi\in C^{m,\alpha}(\T^3;\R^3)$ are smooth functions such that the phase function $\Phi$ obeys
\begin{equation*}
{C}^{-1}\leq \abs{\nabla \Phi} \leq {C}
\end{equation*}
 on $\T^3$, for some constant $C\geq 1$. Then, with the inverse divergence operator $\RSZ$ defined in \eqref{eq:RSZ} we have
\begin{align*}
\norm{\RSZ \left( a(x) e^{i \lambda \xi \cdot \Phi(x)} \right) }_{C^\alpha} 
\les \frac{  \norm{a}_{C^0}}{\lambda^{1-\alpha}} +   \frac{ \norm{a}_{C^{m,\alpha}}+\norm{a}_{C^0}\norm{\nabla \Phi}_{C^{m,\alpha}}}{\lambda^{m-\alpha}} \, ,
\end{align*}
where the implicit constant depends on $C$, $\alpha$ and $m$ (in particular, not on the frequency $\lambda$).
\end{lemma}

The above lemma is used when estimating the $C^0$ norm of the new stress. Indeed, for a fixed $t \in [(j-1) \ell, (j+1)\ell]$, in view of the bounds \eqref{eq:Phi:j:bounds},  assuming that $\ell \lambda_q \delta_q^{\sfrac 12} \ll 1$ we have that $\frac 12 \leq |\nabla \Phi_j(\cdot,t)| \leq 2$ on $\T^3$. 
Thus, by Lemma~\ref{lem:stationary:phase} we obtain that if $a$ is a smooth periodic function such that 
\begin{align}
\norm{a}_{C^n} \les C_a \ell^{-n}
\label{eq:stationary:phase:in}
\end{align}
holds for some constant $C_a >0$, for all $0 \leq n\leq m+1$, where the implicit constant only depends on $m$, and if $m+1  \geq \sfrac{1}{\alpha} $, then
\begin{align}
\norm{\RSZ \left( a \, W_{(\xi)}\circ \Phi_j \right) }_{C^\alpha} 
\les \frac{C_a}{\lambda_{q+1}^{1-\alpha}} \left( 1 +   \frac{\ell^{-m-1}}{\lambda_{q+1}^{m-1}}  \right) \les \frac{C_a}{\lambda_{q+1}^{1-\alpha}} \, .
\label{eq:stationary:phase:out}
\end{align} 
The implicit constant depends only on $\alpha$ and $m$. In the second inequality above we have used that  $\ell^{-1} \leq \lambda_{q+1}^{\sfrac 34}$, and thus  $\ell^{-m-1} \lambda_{q+1}^{1-m}  \leq 1$ as soon as $m \geq 7$.

The same proof that was used in~\cite{DSZ17} to prove Lemma~\ref{lem:stationary:phase} gives another useful estimate.  Let $\xi \in \Lambda_{j}$ and $\xi' \in \Lambda_{j'}$ for $|j-j'|\leq 1$ be such that $\xi + \xi'\neq 0$.  Then we have that $|\xi+\xi'|\geq c_*> 0$, for some universal constant $c_* \in(0,1)$. By appealing to the estimates \eqref{eq:Phi:j:bounds}, one may show that for a smooth periodic function $a(x)$ which obeys \eqref{eq:stationary:phase:in} for some constant $C_a>0$, we have that 
\begin{align}
\norm{\RSZ \left( a \left( W_{(\xi)}\circ \Phi_j \otimes W_{\xi'} \circ \Phi_{j'}\right) \right) }_{C^\alpha} 
\les \frac{C_a}{\lambda_{q+1}^{1-\alpha}} \left( 1 +   \frac{\ell^{-m-1}}{\lambda_{q+1}^{m-1}}  \right) \les \frac{C_a}{\lambda_{q+1}^{1-\alpha}} \, .
\label{eq:stationary:phase:double}
\end{align} 
The implicit constant depends only on $\alpha$ and $m$. 

\subsubsection{Decomposition of the new Reynolds stress}

Our goal is to show that the vector field  $v_{q+1}$ defined in \eqref{eq:v_q+1:def} obeys \eqref{e:euler_reynolds} at level $q+1$, for a Reynolds stress $\RR_{q+1}$ and pressure scalar $p_{q+1}$ which we are computing next. Upon subtracting  \eqref{e:euler_reynolds} at level $q+1$ the system \eqref{e:euler_reynolds_ell} we obtain that 
\begin{align}
\div \RR_{q+1} - \nabla p_{q+1}
&=\underbrace{(\partial_t +v_{\ell}\cdot \nabla) w^{(p)}_{q+1}}_{\div (R_{\rm transport})}+ \underbrace{\div(w_{q+1}^{(p)} \otimes w_{q+1}^{(p)} + \RR_{\ell})}_{\div (R_{\rm oscillation}) + \nabla p_{\rm oscillation}} \notag \\
& \qquad+ \underbrace{w_{q+1}\cdot  \nabla v_{\ell}}_{\div( R_{\rm Nash})}  +\underbrace{ (\partial_t + v_\ell \cdot\nabla) w^{(c)}_{q+1}+ \div\left(w_{q+1}^{(c)} \otimes w_{q+1}+ w_{q+1}^{(p)} \otimes w_{q+1}^{(c)}\right) }_{\div (R_{\rm corrector})+ \nabla p_{\rm corrector}} \notag \\
&\qquad+ \div\left(R_{\rm commutator}\right) - \nabla  p_\ell \,.
\label{eq:R:q+1:A}
\end{align}
Here, $R_{\rm commutator}$ is as defined in \eqref{e:euler_reynolds_ell:a}, and we have used the inverse divergence operator from \eqref{eq:RSZ} to define
\begin{align}
R_{\rm transport} &:= \RSZ\left( (\partial_t +v_{\ell}\cdot \nabla) w^{(p)}_{q+1} \right) \label{eq:Euler:R:transport}\\
R_{\rm Nash} &:= \RSZ\left( w_{q+1} \cdot \nabla v_\ell \right)\label{eq:Euler:R:Nash} \\
R_{\rm corrector} &:= \RSZ\left( (\partial_t + v_\ell \cdot\nabla) w^{(c)}_{q+1} \right)+  \left(w_{q+1}^{(c)} \mathring \otimes w_{q+1}^{(c)} + w_{q+1}^{(c)} \mathring \otimes w_{q+1}^{(p)} + w_{q+1}^{(p)} \mathring\otimes w_{q+1}^{(c)} \right) \label{eq:Euler:R:corrector} 
\end{align}
while $p_{\rm corrector} := 2 w_{q+1}^{(c)} \cdot w_{q+1}^{(p)} + |w_{q+1}^{(c)}|^2$.
The remaining stress $R_{\rm oscillation}$ and corresponding pressure $p_{\rm oscillation}$ are defined as follows.

First, note that for $j, j'$ such that $|j-j'|\geq 2$, we have $\chi_j(t) \chi_{j'}(t) =0$. Second, for $|j-j'|=1$, we have that $\Lambda_{j} \cap \Lambda_{j'} = \emptyset$. And third, we note that similarly to \eqref{e:Bequation} we have $\div ( W_{(\xi)} \otimes W_{(\xi')} + W_{(\xi')} \otimes W_{(\xi)} ) = \nabla (W_{(\xi)} \cdot W_{(\xi')})$, which follows from the identity 
$(B_\xi \otimes B_{\xi'} + B_{\xi'} \otimes B_\xi) \cdot (\xi+\xi') = (B_\xi \cdot B_{\xi'}) (\xi+\xi')$.  Therefore, we may use \eqref{e:av_of_Bel}, \eqref{eq:amplitudes:are:good}, \eqref{eq:chi:partition}--\eqref{e:def_wo}, and the notation \eqref{e:phi_xi} to write
\begin{align*}
&\div(w_{q+1}^{(p)} \otimes w_{q+1}^{(p)} + \RR_{\ell})  = \div\left( \sum_{j,j',\xi,\xi'}   w_{(\xi)} \otimes w_{(\xi')} +\RR_{\ell}\right)\\
&= \div\left( \sum_{j,\xi}   w_{(\xi)} \otimes w_{(-\xi)} +\RR_{\ell}\right) +   \sum_{j,j',\xi+\xi' \neq 0}  \div \left( w_{(\xi)} \otimes w_{(\xi')} \right)  \\
&= \div\left(\frac 12 \sum_{j,\xi}     a_{(\xi)}^2 \left(\Id - \xi \otimes \xi\right) +\RR_{\ell}\right) 
+  \sum_{j,j',\xi+\xi' \neq 0}  \div \left( a_{(\xi)} a_{(\xi')} \phi_{(\xi)} \phi_{(\xi')} W_{(\xi)} \otimes W_{(\xi')} \right)  \\
&=  \frac 12 \sum_{j,j',\xi+\xi'\neq 0}   a_{(\xi)} a_{(\xi')} \phi_{(\xi)} \phi_{(\xi')} \div \left( W_{(\xi)} \otimes W_{(\xi')} + W_{(\xi')} \otimes W_{(\xi)} \right) \notag\\
&\qquad + \sum_{j,j',\xi+\xi'\neq 0}  
\left(W_{(\xi)} \otimes W_{({\xi'})}\right) \nabla \left( a_{(\xi)}a_{(\xi')}\phi_{(\xi)}\phi_{(\xi')} \right)   \\
&=  \frac 12 \sum_{j,j',\xi+\xi'\neq 0}   a_{(\xi)} a_{(\xi')} \phi_{(\xi)} \phi_{(\xi')} \nabla \left( W_{(\xi)} \cdot W_{(\xi')}  \right)  + \sum_{j,j',\xi+\xi'\neq 0} 
\left(W_{(\xi)} \otimes W_{({\xi'})}\right) \nabla \left( a_{(\xi)}a_{(\xi')}\phi_{(\xi)}\phi_{(\xi')} \right)   \\
&=\nabla p_{\rm oscillation} +  \div(R_{\rm oscillation}) \,.
\end{align*}
Above, we   have   denoted
\[  p_{\rm oscillation} :=  \frac 12 \sum_{j,j',\xi+\xi'\neq 0}  a_{(\xi)} a_{(\xi')} \phi_{(\xi)} \phi_{(\xi')}   \left( W_{(\xi)} \cdot W_{(\xi')}  \right) \, , \]
and 
\begin{align}
R_{\rm oscillation} &:= \sum_{j,j',\xi+\xi'\neq 0}  \RSZ \left(  \left(W_{(\xi)} \otimes W_{({\xi'})} - \frac{W_{(\xi)} \cdot W_{(\xi')}}{2} \Id \right) \nabla\left(a_{(\xi)}a_{(\xi')}\phi_{(\xi)}\phi_{(\xi')}\right)  \right)\, .
\label{eq:Euler:R:oscillation}
\end{align}

At this point all the terms in \eqref{eq:R:q+1:A} are well defined. We have $p_{q+1} = p_\ell - p_{\rm oscillation}  - p_{\rm corrector}$ and  
\begin{align}
\RR_{q+1} = R_{\rm transport} + R_{\rm oscillation} + R_{\rm Nash} + R_{\rm corrector} + R_{\rm commutator} \,
\label{eq:R:q+1:B}
\end{align}
which are stresses defined in \eqref{eq:Euler:R:transport}, \eqref{eq:Euler:R:oscillation}, \eqref{eq:Euler:R:Nash}, \eqref{eq:Euler:R:corrector}, and respectively \eqref{e:euler_reynolds_ell:a}.

\subsubsection{Estimates for the new Reynolds stress}
\label{sec:Reynolds:stress:C0}
To conclude the proof of the inductive lemma, we need to show that the stress defined in \eqref{eq:R:q+1:B} obeys the bound \eqref{eq:Rq:C0:inductive} at level $q+1$. Recall that the commutator stress was bounded in \eqref{eq:Rc:bound}, and that it obeys a suitable bound if $\ell$ is sufficiently small. The main terms are the transport error and the oscillation error, which we bound first.

{\bf Transport error. \,}
Inspecting the definition of $w^{(p)}_{q+1}$ from \eqref{e:w_xi} and \eqref{e:def_wo}, we notice that the material derivative cannot land on the highest frequency term, namely $W_{(\xi)} \circ \Phi_j$, as this term is perfectly transported by $v_\ell$. Therefore, we have 
\begin{align*}
(\partial_t +v_{\ell}\cdot \nabla)w_{q+1}^{(p)}  = \sum_{j, \xi}  (\partial_t + v_\ell \cdot \nabla) a_{(\xi)} \, W_{(\xi)} \circ \Phi_j \, .
\end{align*}
Returning to the definition of $a_{(\xi)}$ in \eqref{e:a_xi} we may show using standard mollification estimates, that the bounds \eqref{eq:Rq:C0:inductive}  and \eqref{eq:v:ell:C0}, imply 
\begin{align}
\norm{ (\partial_t + v_\ell \cdot \nabla) \RR_\ell }_{C^0} 
\les \norm{\partial_t \RR_\ell}_{C^0} + \norm{v_\ell}_{C^0} \norm{\RR_\ell}_{C^1}
\les \ell^{-1} \norm{\RR_q}_{C^0} 
\les \ell^{-1} \delta_{q+1}.
\label{eq:Dt:R:ell:C0}
\end{align}
In a similar spirit to the above estimate, and taking into account \eqref{e:chi_l},  we may in fact show that 
\begin{align}
\label{eq:a:Cn}
\norm{a_{(\xi)}}_{C^n} \les \delta_{q+1}^{\sfrac 12} \ell^{-n}, \qquad \mbox{and} \qquad \norm{(\partial_t + v_\ell \cdot \nabla) a_{(\xi)}}_{C^n} \les \delta_{q+1}^{\sfrac 12} \ell^{-1-n}
\end{align}
holds for all $n \geq 0$. Thus both $a_{(\xi)}$ and the material derivative of $a_{(\xi)}$ obey the bound \eqref{eq:stationary:phase:in}, with $C_a = \delta_{q+1}^{\sfrac 12}$,  respectively $C_a = \delta_{q+1}^{\sfrac 12} \ell^{-1}$. We thus conclude from \eqref{eq:stationary:phase:in} that for a sufficiently large universal $m$, we have
\begin{align*}
\norm{R_{\rm transport}}_{C^0} \les \frac{\ell^{-1} \delta_{q+1}^{\sfrac 12}}{\lambda_{q+1}^{1-\alpha}} \left( 1 + \frac{\ell^{-m-1}}{\lambda_{q+1}^{m-1}} \right) \les \frac{\ell^{-1} \delta_{q+1}^{\sfrac 12}}{\lambda_{q+1}^{1-\alpha}} \les \lambda_{q+1}^{\alpha - \sfrac 14}  \ll \delta_{q+2} \, .
\end{align*}
Here we have taken $\alpha$ and $\beta$ sufficiently small, and $a$ sufficiently large.

{\bf Oscillation error. \,} For the oscillation error we apply the \eqref{eq:stationary:phase:double} version of the stationary phase estimate. 
First we use \eqref{eq:Euler:R:oscillation} to rewrite
\begin{align*}
R_{\rm oscillation}
&= \sum_{j,j',\xi+\xi'\neq 0}   \RSZ \Bigg( \left(W_{(\xi)}\circ \Phi_j \otimes W_{({\xi'})}\circ\Phi_{j'} - \frac{W_{(\xi)}\circ\Phi_{j} \cdot W_{(\xi')}\circ\Phi_{j'}}{2} \Id \right) \notag\\
&\qquad \qquad   \Big( \nabla(a_{(\xi)}a_{(\xi')}) + i \lambda_{q+1} a_{(\xi)} a_{(\xi')} \left( (\nabla \Phi_j - \Id) \cdot \xi +  (\nabla \Phi_{j'} - \Id) \cdot \xi'  \right)\Big)  \Bigg) \, .
\end{align*}
Then, similarly to \eqref{eq:a:Cn} we have 
\begin{align*}
 \norm{\nabla(a_{(\xi)}a_{(\xi')})}_{C^n}  
 \les  \delta_{q+1} \ell^{-n-1} 
\end{align*}
and by also appealing to \eqref{eq:Phi:j:bounds} we also obtain
\begin{align*}
\lambda_{q+1} \norm{a_{(\xi)} a_{(\xi')} \left( (\nabla \Phi_j - \Id) \cdot \xi +  (\nabla \Phi_{j'} - \Id) \cdot \xi'  \right)}_{C^n} 
\les \lambda_{q+1} \delta_{q+1} (\ell \lambda_q \delta_q^{\sfrac 12}) \ell^{-n} 
\end{align*}
for   $n \geq 0$. Using \eqref{eq:ell:def} and \eqref{eq:stationary:phase:double} we thus obtain from that the oscillation stress as defined in \eqref{eq:Euler:R:oscillation} obeys
\begin{align*}
\norm{R_{\rm oscillation}}_{C^0}
&\les \frac{\delta_{q+1} \ell^{-1} + \lambda_{q+1} \delta_{q+1} (\ell \lambda_q \delta_{q}^{\sfrac 12})}{\lambda_{q+1}^{1-\alpha}}
\les \lambda_{q+1}^{\alpha-\sfrac 14}  + \lambda_{q}^{2\alpha-\sfrac 12}    \ll \delta_{q+2} \, 
\end{align*}
as desired.

{\bf Nash error. \,}
Using that 
\begin{align*}
 w_{q+1} \cdot \nabla v_\ell = \sum_{j,\xi}   a_{(\xi)} \, W_{(\xi)}\circ \Phi_j \cdot \nabla v_\ell \, ,
\end{align*}
and the available estimate
\begin{align*}
 \norm{a_{(\xi)} \, \nabla v_\ell}_{C^n} \les \delta_{q+1}^{\sfrac 12} \lambda_q \delta_q^{\sfrac 12} \ell^{-n} \, ,
\end{align*}
for all $n\geq 0$, allows us to appeal to \eqref{eq:stationary:phase:out} and conclude that
\begin{align*}
 \norm{R_{\rm Nash}}_{C^0} \les \frac{\delta_{q+1}^{\sfrac 12} \lambda_q \delta_q^{\sfrac 12}}{\lambda_{q+1}^{1-\alpha}} \les \lambda_{q}^{2\alpha -1}    \ll \delta_{q+2} \, .
\end{align*}

{\bf Corrector error. \,}  The corrector error has two pieces, the transport derivative of $w_{q+1}^{(c)}$ by the flow of $v_\ell$, and the residual contributions from the nonlinear term, which are easier to estimate due to \eqref{eq:w:q+1:p:L:infinity} and \eqref{eq:w:q+1:c:L:infinity}:
\begin{align*}
\norm{w_{q+1}^{(c)} \otimes w_{q+1}^{(c)} + w_{q+1}^{(c)} \otimes w_{q+1}^{(p)}  + w_{q+1}^{(p)} \otimes w_{q+1}^{(c)}}_{C^0}
\les \delta_{q+1} \left( \frac{1}{\ell \lambda_{q+1}} + \ell \lambda_q \delta_{q}^{\sfrac 12}\right) 
  \ll \delta_{q+2}\, .
\end{align*}
Inspecting the definition of $w^{(c)}_{q+1}$ from \eqref{eq:w:c:small} and \eqref{e:corrector}, we notice that the material derivative cannot land on the highest frequency term, namely $W_{(\xi)} \circ \Phi_j$, as this term is perfectly transported by $v_\ell$. We thus have
\begin{align*}
(\partial_t  + v_{\ell} \cdot  \nabla) w_{q+1}^{(c)} 
&=   \sum_{j,\xi}  \left( (\partial_t  + v_{\ell} \cdot  \nabla) \left(\frac{\nabla a_{(\xi)}}{\lambda_{q+1}}  + i \, a_{(\xi)}    (\nabla \Phi_j(x,t) - \Id) \xi \right) \right)
\times W_{(\xi)}(\Phi_j(x,t)).
\end{align*}
The available estimates for $a_{(\xi)}$ and $\nabla \Phi_j -\Id$ yield 
\begin{align*}
 \norm{a_{(\xi)}    (\nabla \Phi_j(x,t) - \Id)  }_{C^n} \les \delta_{q+1}^{\sfrac 12} (\ell\lambda_q\delta_{q}^{\sfrac 12}) \ell^{-n}
\end{align*}
and using \eqref{eq:Dt:R:ell:C0} we also obtain  
\begin{align*}
\frac{1}{\lambda_{q+1}} \norm{(\partial_t + v_\ell \cdot \nabla) \nabla a_{(\xi)}  }_{C^n} \les \frac{\ell^{-2} \delta_{q+1}^{\sfrac 12}}{\lambda_{q+1}}  \ell^{-n}
\end{align*}
for all $n \geq 0$. From \eqref{eq:stationary:phase:out} we thus obtain that 
\begin{align*}
\norm{\mathcal R\left(\partial_t w^{(c)}_{q+1}+ v_{\ell} \cdot  \nabla w_{q+1}^{(c)} \right)}_{L^{\infty}}\les \frac{\ell^{-2} \delta_{q+1}^{\sfrac 12}}{\lambda_{q+1}^{2-\alpha}} 
\les \lambda_{q}^{2\alpha-1}   
\ll \delta_{q+2} \, .
\end{align*}

\subsubsection{Proof of \eqref{eq:Rq:C0:inductive} at level $q+1$} Summarizing the estimates obtained for the five stresses in \eqref{eq:R:q+1:B}, we obtain that 
\begin{align}
\norm{\RR_{q+1}}_{C^0} \leq c_R \delta_{q+2}.
\end{align}
By taking $a$ sufficiently large and $\beta$ sufficiently small, it follows that the constant $c_R$ may be taken arbitrarily small. This concludes the proof of Proposition~\ref{prop:iteration:Euler}.

\section{Euler: the full flexible part of the Onsager conjecture}
\label{sec:Euler:C1/3}
The result of the previous section gives us the existence of H\"older continuous weak solutions of the 3D Euler equations which are not conservative (more generally, which can attain any given smooth energy profile). In this section our goal is to describe the H\"older $1/3^-$ scheme of~\cite{Isett16,BDLSV17}.

\begin{remark}[The H\"older $1/5^-$ scheme]
In order to achieve a H\"older exponent $< \sfrac 15$, in the proof of Theorem~\ref{thm:Euler:C0} one has to carefully take into account estimates for the {\em material derivative of the Reynolds stress} $(\partial_t + v_\ell \cdot \nabla) \RR_q$ (cf.~\cite{IsettThesis,Buckmaster14,BDLISZ15} for details). In principle, material derivatives should cost less than regular spatial or temporal derivatives. Indeed, by scaling, one expects material derivatives to cost a factor roughly proportional to the Lipschitz norm of $v_q$. Taking advantage of this observation, one can improve on the estimate \eqref{eq:Dt:R:ell:C0}. As it stands, the estimate \eqref{eq:Dt:R:ell:C0} scales particularly badly and is the principal reason the proof of Theorem~\ref{thm:Euler:C0} given above requires significant super-exponential growth in frequency. 
\end{remark}

\begin{remark}[Almost everywhere in time H\"older $1/3^-$ scheme]
As noted in Section~\ref{sec:Reynolds:stress:C0} the principal errors are transport error and the oscillation error. Note that the transport error is concentrated on the subset of times where the temporal cutoffs $\chi_i$ and $\chi_{i+1}$ overlap. In~\cite{Buckmaster15} it was noted that one can obtain any H\"older exponent $<\sfrac 13$ {\em almost everywhere in time}, by designing a scheme that concentrates such errors on a zero measure set of times. By taking advantage of this idea and using a delicate bookkeeping scheme, in~\cite{BDLSZ16} non-conservative solutions were constructed in the space $L^1_t C^{\sfrac 13-}$. 
\end{remark}

The flexible side of the Onsager conjecture was resolved by Isett in~\cite{Isett16}, who proved the existence of non-conservative weak solutions of 3D Euler in the regularity class $C^\beta_{x,t}$, for any $\beta < \sfrac 13$. The proof of Isett builds upon the ideas in the above mentioned works, and utilizes two new key ingredients. The first, is the usage of Mikado flows which were introduced earlier by Daneri and Sz\'ekelyhidi~\cite{DSZ17}. These are a rich family of pressure-less stationary solutions of the 3D Euler equation (straight pipe flows), which have a better (when compared to Beltrami flows) self-interaction behavior in the oscillation error  as they are advected by a mean flow. We discuss Mikado flows in Section~\ref{sec:Mikado} below. The second key ingredient is due to Isett, and may be viewed as the principal main idea in his proof: prior to adding the  convex integration perturbation $w_{q+1}$, it is very useful to  replace the approximate solution $(v_q,\RR_q)$ with another pair $(\overline v_q, \Rbar_q)$, which has the property that $\overline v_q$ is close to $v_q$, but more importantly, $\Rbar_q(t)$ vanishes on every other interval of size $\approx \norm{\nabla v_q}_{C^0}^{-1}$ within $[0,T]$. We discuss this gluing procedure in Section~\ref{sec:gluing} below. The main result of~\cite{Isett16} is Theorem~\ref{thm:Onsager:critical:Isett} above. This result was subsequently extended in~\cite{BDLSV17} to the class of dissipative weak solutions, which in particular can attain any given energy profile, cf.~Theorem~\ref{thm:Onsager:critical:dissipative} above.

In this section, we will present an outline of the arguments employed in \cite{BDLSV17} in order to prove Theorem~\ref{thm:Onsager:critical:dissipative}. The principal new idea in \cite{BDLSV17} that allows for the construction of dissipative solutions involves adding kinetic energy at each iterative step that \emph{wiggles} through space-time (this is explained in detail in Section \ref{sec:Mikado:w}).  Additionally, the gluing procedure is implemented differently in~\cite{BDLSV17} which leads to a more efficient proof when compared to \cite{Isett16}. Since the proof outlined below will follow \cite{BDLSV17}  very closely, some details will be omitted and we refer the reader to \cite{BDLSV17} for the complete proof. We note however that compared to \cite{BDLSV17}, in Section \ref{sec:Mikado}, we have adopted a slightly different presentation of Mikado flows in order that their construction be more directly comparable to the intermittent jets of Section~\ref{sec:NSE:L2}.

Finally, we note that in \cite{Isett17}, Isett showed that one can further optimize the schemes of  \cite{Isett16,BDLSV17} in order to construct non-conservative weak solutions to the Euler equations that lie in the intersection of all H\"older spaces $C^{\beta}$ for $\beta<\sfrac13$. It is an open problem to determine whether non-conservative weak solutions to the Euler equations exist that have H\"older exponent exactly $\sfrac13$. Or, in view of~\cite{CCFS08}, it would be desirable to determine whether there exist non-conservative weak solutions which lie in the space $L^3_t B^{\sfrac 13}_{3,\infty,x}$.

\subsection{Iteration lemma and proof of Theorem~\ref{thm:Onsager:critical:dissipative}}
Define the frequency parameter $\lambda_q$ and the amplitude parameter $\delta_q$ by
\begin{subequations}
\begin{align}
\lambda_q&=  \lceil a^{(b^q)} \rceil \label{e:freq_def}\\
\delta_q&=\lambda_q^{-2\beta} \label{e:size_def}
\end{align}
\end{subequations}
where $a>1$ is a  large parameter, $b>1$ is close to $1$ and $0<\beta<\sfrac13$. Note that with this notation we have $1\leq \lambda_q a^{-(b^q)} \leq 2$. To simplify matters, we will assume the prescribed energy profile satisfies the   bound
\begin{equation}\label{e:energy_unif_bnd}
\sup_t e'(t)\leq 1\,.
\end{equation}
In the proof of Theorem~\ref{thm:Onsager:critical:dissipative} below, we will demonstrate how this additional constraint may be removed.

The inductive estimates on the approximate solution $(v_q,\RR_q)$ of the Euler-Reynolds system \eqref{e:euler_reynolds} are nearly identical to \eqref{eq:inductive:Euler}, except that now we also include the estimate for the energy iterate:
\begin{subequations}
\label{eq:inductive:C1/3}
\begin{align}
\norm{\mathring R_q}_{C^0}&\leq  \delta_{q+1}\lambda_q^{-3\alpha}\label{e:R_q_inductive_est}\\
\norm{v_q}_{C^1}&\leq M \delta_q^{\sfrac12}\lambda_q\label{e:v_q_inductive_est}\\
\norm{v_q}_{C^0} & \leq 1- \delta_q^{\sfrac12}\label{e:v_q_0}\\
\delta_{q+1}\lambda_q^{-\alpha} &\leq e(t)-\int_{\T^3}\abs{v_q}^2\,dx\leq \delta_{q+1}\label{e:energy_inductive_assumption} \, .
\end{align}
\end{subequations}
Here $0 < \alpha  < 1$ is a small parameter to be chosen suitably in terms of $\beta$ and $b$,  and $M$ is a universal geometric constant which is fixed throughout the iteration scheme.  The wiggle room given by the $\lambda_q^{-3\alpha}$ factor in \eqref{e:R_q_inductive_est} (when compared to \eqref{eq:Rq:C0:inductive}), is useful in the gluing step $\RR_q \mapsto \Rbar_q$, and in bounding errors arising from the fact that Calder\'on-Zygmund operators are not bounded on $C^0$. The iterative proposition used to go from step $q$ to step $q+1$ is as follows. 

\begin{proposition}\label{p:main} Assume $0<\beta<\sfrac13$ and $0 < b - 1 < \sfrac{(1-3\beta)}{2\beta}$. 
Then there exists an $\alpha_0 = \alpha_0(\beta,b) >0$, such that for any $0<\alpha<\alpha_0$ there exists an $a_0 = a_0(\beta,b,\alpha,M)>0$, such that for any $a\geq a_0$ the following holds: Given $(v_q,\mathring R_q)$ which solve \eqref{e:euler_reynolds} and obey the estimates \eqref{eq:inductive:C1/3}, there exists a solution  $(v_{q+1}, \mathring R_{q+1})$ to \eqref{e:euler_reynolds} which satisfy  \eqref{eq:inductive:C1/3} with $q$ replaced by $q+1$, and moreover, we have 
\begin{equation}
\norm{v_{q+1}-v_q}_{C^0}
+ \lambda_{q+1}^{-1} \norm{v_{q+1}-v_q}_{C^1} 
\leq M\delta_{q+1}^{\sfrac12}\label{e:v_diff_prop_est}.
\end{equation}
\end{proposition}

\subsubsection{Proof of Theorem~\ref{thm:Onsager:critical:dissipative}}

We first observe, that without loss of generalization, we may restrict ourselves to considering \emph{normalized} energy profiles that satisfy the following estimates
\begin{equation}\label{e:normalized_energy}
\inf_t e(t) \geq \delta_{1}\lambda_0^{-\alpha},\qquad\sup_t e(t)\leq \delta_{1}, \quad\mbox{and}\quad \sup_t e'(t)\leq 1\, ,
\end{equation}
Indeed, the reduction to normalized energy profiles follows as a consequence of the scaling scale invariance $v(x,t)\mapsto \Gamma v(x,\Gamma t)$ of the Euler equation (see \cite{BDLSV17} for details).

Observe that by setting $(v_0,R_0)=(0,0)$, the pair $(v_0, R_0)$ trivially satisfying the assumptions  \eqref{e:R_q_inductive_est}--\eqref{e:v_q_0}; moreover, \eqref{e:energy_inductive_assumption} and \eqref{e:energy_unif_bnd} are implied by \eqref{e:normalized_energy}.  Applying Proposition \ref{p:main} iteratively, we obtain a sequence of velocities $v_q$ converging uniformly to a weak solution $v$ of the Navier-Stokes equations \eqref{eq:NSE}. Utilizing the estimate \eqref{e:v_diff_prop_est} yields
\begin{align*}
\sum_{q=0}^{\infty} \norm{v_{q+1}-v_q}_{C^{\beta'}} &
\lesssim \; \sum_{q=0}^{\infty} \norm{v_{q+1}-v_q}_{C^0}^{1-\beta'}\norm{v_{q+1}-v_q}_{C^1}^{\beta'}
\lesssim \; \sum_{q=0}^{\infty} \delta_{q+1}^{\frac{1-\beta'}2}\left(\delta_{q+1}^{\sfrac12}\lambda_q\right)^{\beta'} \notag
\lesssim  \sum_{q=0}^{\infty} \lambda_q^{\beta'-\beta}
\end{align*}
Hence we obtain that $v\in C^0_tC^{\beta'}_x$ for all $\beta'<\beta$. The time regularity follows as a consequence of the work \cite{Isett13} (see also \cite{BDLSV17} for a simplified argument).

\begin{remark}[Vanishing viscosity of smooth solutions to the forced Navier-Stokes equations]
\label{rem:vanishing:viscosity}
A direct consequence of Proposition~\ref{p:main} is that for any $\beta < \sfrac 13$ one may prove the following statement, which is motivated by the discussion in Section~\ref{sec:physics}. For every $\nu>0$, one may find a $C^\infty$ smooth solution $v^\nu$ of the Navier-Stokes equations~\eqref{eq:NSE} with a $C^\infty$ smooth forcing term $f^\nu = \div F^\nu$, such that as $\nu \to 0$ we have:   $\norm{F^\nu}_{C^0} \les \nu^{\frac{2\beta}{(1+\beta)}}$ and thus $F^\nu \to 0$ as $\nu\to 0$;  $\norm{v^\nu}_{C^{\beta}} \les 1$, independently of $\nu$; and such that $v^\nu$ converges as $\nu\to 0$ to a dissipative weak solution $v$ of the of Euler equation. However, inspecting the energy balance for this forced Navier-Stokes solution, we see that the work of force term on the right side of \eqref{eq:NSE:global:energy:balance} has size $\approx \nu^{\frac{3\beta-1}{1+\beta}}$, which does not remain bounded as $\nu \to 0$ since $\beta < \sfrac 13$. In light of this remark, it would be interesting to see whether a convex-integration scheme could be designed to reach the endpoint exponent $\beta = \sfrac 13$. See Problem~\ref{p:Ronaldo} below.
\end{remark}

\subsection{Mollification}

Similarly to Section~\ref{sec:mollify}, here we replace the pair $(v_q,\RR_q)$ with a mollified pair $(v_\ell,\RR_\ell)$ which obey the Euler-Reynolds system. Notationally, unlike the definition of $\RR_\ell$ in Section~\ref{sec:mollify}, here, $\RR_\ell$ will denote the sum of the mollified stress and the commutator stress (compare to \eqref{eq:v:R:ell:def} and \eqref{e:euler_reynolds_ell:a}). For   $\ell >0$ defined as
\begin{align}
\label{e:ell_def}
\ell =\frac{\delta_{q+1}^{\sfrac 12}}{\delta_q^{\sfrac12 }\lambda_q^{1+\sfrac{3\alpha}{2}}}\, ,
\end{align}
we let
\begin{align*}
v_{\ell}:=& v_q \ast_x \phi_\ell\\
\RR_{\ell}:=& \mathring R_q \ast_x  \phi_\ell  -(v_q\mathring\otimes v_q)\ast_x \phi_\ell  + v_{\ell}\mathring\otimes v_{\ell} 
\end{align*}
which obey the Euler-Reynolds system~\eqref{e:euler_reynolds} for a suitable pressure $p_\ell$ which has zero mean. 

In order to bound $\RR_\ell$ we use  from~\cite[Lemma 1]{CDLSJ12} the following generalization of the Constantin-E-Titi~\cite{ConstantinETiti94} commutator estimate:
\begin{proposition}\label{p:CET}
Let $f,g\in C^{\infty}(\T^3\times[0,1])$ and $\psi$ a standard radial smooth and compactly supported kernel. For any $r\geq 0$ and $\theta \in (0,1]$ we have the estimate
\[
\Bigl\|(f*\psi_ \ell)( g*\psi_\ell)-(fg)*\psi_ \ell\Bigr\|_{C^r} \les \ell^{2\theta -r}  \|f\|_{C^\theta} \|g\|_{C^\theta} \, ,
\]
where the implicit constant  depends only on $r$ and $\psi$.
\end{proposition}
With the choice of $\ell$ given above in \eqref{e:ell_def}, using the usual mollification bounds and Proposition~\ref{p:CET}, we obtain (cf.~\cite[Proposition 2.2]{BDLSV17})
\begin{subequations}
\label{eq:v_ell:all:bounds}
\begin{align}
\norm{v_{\ell}-v_q}_{C^0}&\lesssim \delta_{q+1}^{\sfrac12}\lambda_q^{-\alpha}\,,\label{e:v:ell:0}\\
\norm{v_{\ell}}_{C^{N+1}} &\lesssim \delta_q^{\sfrac 12}\lambda_q\ell^{-N}  \,,  \label{e:v:ell:k}\\
\norm{\mathring{R}_{\ell}}_{C^{N+\alpha}}&\lesssim  \delta_{q+1}\ell^{-N+\alpha} \,, \label{e:R:ell}\\
\abs{\int_{\T^3}\abs{v_q}^2-\abs{v_{\ell}}^2\,dx} &\lesssim \delta_{q+1}\ell^{\alpha}\,,
\label{e:vq_vell_energy_diff}
\end{align}
\end{subequations}
for all $N\geq 0$ where the implicit constant may depend on $N$ and $\alpha$.

\subsection{Gluing}\label{sec:gluing}
The gluing step is the fundamental new idea introduced in~\cite{Isett16}. 
The idea is to replace $v_\ell$ with $\overline v_q$, which is obtained by gluing in time exact solutions $v_i$ of the three dimensional Euler equations whose initial datum matches $v_\ell$ at certain times $t_i$. The resulting glued stress $\Rbar_q$ will then be localized on pairwise disjoint time intervals $I_i$ of length proportional to $\tau_q = |t_{i} - t_{i-1}|$, and   that the spacing between this intervals is also proportional to  $ \tau_q$. The stress $\RR_\ell$ will then be corrected in a convex integration step involving Mikado flows (cf.\ Section~\ref{sec:Mikado:w}). As mentioned previously, the Mikado flows need to be advected by the mean flow $v_{\ell}$ causing them to be twisted and deformed. However, if the intervals $I_i$ and $I_j$, with $i\neq j$, were to intersect, then this could result in cross-interactions between the deformed Mikado flows, leading to unacceptable Reynolds errors. In this subsection we describe an efficient  implementation of the gluing procedure, as detailed in~\cite{BDLSV17}. 

\subsubsection{Exact solutions of the Euler equations and their stability.}
To make the above idea precise, we introduce the parameter
\begin{equation}\label{e:tau_def}
\tau_q= \frac{\ell^{2\alpha}}{\delta_q^{\sfrac 12}\lambda_q}.
\end{equation} 
and define 
\begin{align}
t_i = i \tau_q, \qquad I_i = [t_i + \sfrac{\tau_q}{3}, t_i + \sfrac{2\tau_q}{3}] \cap [0,T], \qquad J_i = (t_i - \sfrac{\tau_q}{3}, t_i+ \sfrac{\tau_q}{3}) \cap [0,T].
\label{eq:ti:Ii}
\end{align}
In order to justify the choice \eqref{e:tau_def}, we recall that the Euler equations with $C^{1+\alpha}$ initial datum may be solved uniquely locally in time on a time interval with length which is inversely proportional to the $C^{1+\alpha}$ norm of the initial datum. From \eqref{e:v:ell:k} and \eqref{e:tau_def} we note that indeed $\tau_q$ obeys the CFL-like condition:
\begin{equation}
\tau_q \norm{v_{\ell}}_{1+\alpha} \lesssim  \tau_q\delta_q^{\sfrac{1}{2}} \lambda_q \ell^{-\alpha}  \lesssim \ell^{\alpha} \ll 1 \label{e:CFL}
\end{equation}
as long as $a$ is sufficiently large. 

With the above considerations, for each $i$ we solve on $[t_i - \tau_q, t_i+\tau_q]$ the Euler equations exactly, with initial data given by $v_\ell(\cdot,t_i)$:
\begin{subequations}
\label{eq:Euler:i}
\begin{align}
\partial_t v_i+ v_i \cdot \nabla v_i +\nabla p_i &=0\\
\div v_i &= 0\\
v_i(\cdot,t_i)&=v_\ell(\cdot,t_i)
\end{align}
\end{subequations}
By \eqref{e:CFL} and the classical local existence theorem~\cite{Hoelder33} for \eqref{eq:Euler:i} we have that $v_i$ is a uniquely defined $C^{1+\alpha}$ smooth function on  $\T^3 \times [t_i - \tau_q, t_i+\tau_q]$. Moreover, using bounds similar to \eqref{eq:Phi:j:bounds}, we have that 
\begin{align}
\norm{v_i(t)}_{C^{N+\alpha}} \les \delta_q^{\sfrac 12} \lambda_q \ell^{1-N-\alpha} \les \tau_q^{-1} \ell^{1-N+\alpha}
\label{eq:v_i:C^N}
\end{align}
holds for all $t\in [t_i-\tau_q,t_i+\tau_q]$ and $N\geq 1$.

In fact we can say more: the solution $v_i(\cdot,t)$ is in fact very close to $v_\ell(\cdot,t)$ for all $|t-t_i|\leq \tau_q$, in the sense that the following estimates hold (cf.~\cite[Proposition 3.3]{BDLSV17})
\begin{subequations}
\label{eq:vi:minus:v_ell}
\begin{align}
 \norm{v_i-v_{\ell}}_{C^{N+\alpha}} &\lesssim  \tau_q\delta_{q+1}\ell^{-N-1+\alpha}\,, \label{e:z_diff_k}\\
\norm{\nabla(p_{i} - p_{\ell})}_{C^{N+\alpha}} &\lesssim \delta_{q+1}\ell^{-N-1+\alpha}\,, \label{e:pressure_1}\\
\norm{(\partial_t + v_\ell \cdot \nabla) (v_i-v_{\ell})}_{C^{N+\alpha}} &\lesssim  \delta_{q+1}\ell^{-N-1+\alpha}\,, \label{e:z_D_t}
\end{align}
\end{subequations}
for all $N\geq 0$, uniformly on $[t_i-\tau_q,t_i+\tau_q]$. The bounds \eqref{eq:vi:minus:v_ell} hold upon noting that the incompressible vector field $v_i - v_\ell$ obeys
\begin{align}
\left(\partial_t + v_{\ell} \cdot \nabla\right) (v_{\ell} - v_i)  
= (v_i - v_{\ell}) \cdot \nabla v_i - \nabla (p_{\ell} - p_i)+\div \mathring{R}_{\ell}.
\label{e:z_diff_evo}
\end{align}
with $0$ initial condition at $t = t_i$, whereas $\nabla(p_i - p_\ell)$ is given explicitly as
\begin{align}
\nabla (p_{\ell} - p_i) = \nabla \Delta^{-1} \div\bigl( \nabla v_\ell(v_i-v_\ell) + \nabla v_i{(v_i-v_\ell)} + \div\mathring{R}_{\ell} \bigr) \, .
\label{e:eqnpi}
\end{align}
Equations \eqref{e:z_diff_evo} and \eqref{e:eqnpi} imply the desired estimates \eqref{eq:vi:minus:v_ell} using standard estimates for the transport equation (see also~\eqref{eq:Phi:j:bounds}), mollification and interpolation bounds for H\"older spaces.

Lastly, as we shall see in \eqref{eq:bar:R_q:def} below,  we require sharp estimates on $v_i- v_\ell$ in negative order spaces. For this purpose we introduce vector potentials associated to the exact solutions $v_i$
\begin{align*}
 z_i = {\mathcal B} v_i = (-\Delta)^{-1} \curl v_i \, .
\end{align*}
Note that $\div z_i = 0$ and $\curl z_i = v_i - \int_{\T^3} v_i dx$. Similarly we use the Biot-Savart operator ${\mathcal B}$ to define $z_\ell = {\mathcal B} v_\ell$. With this notation we then have (cf.~\cite[Proposition 3.4]{BDLSV17})
\begin{subequations}
\label{eq:vector:potential:stability}
 \begin{align}
\norm{z_i-z_\ell}_{C^{N+\alpha}} &\lesssim  \tau_q\delta_{q+1}\ell^{-N+\alpha}\,,   \label{e:z_diff} \\
\norm{(\partial_t + v_\ell \cdot \nabla) (z_i-z_\ell)}_{C^{N+\alpha}} &\lesssim  \delta_{q+1}\ell^{-N+\alpha}\,, \label{e:z_diff_Dt}
\end{align}
\end{subequations}
for all $t \in [t_i-\tau_q,t_i+\tau_q]$ and all $N\geq 0$. Proving \eqref{eq:vector:potential:stability} for $N=0$ already contains all the necessary ideas required for $N\geq 0$. Here the main idea is to explicitly compute the commutator between the advective derivative operator $\partial_t + v_\ell\cdot\nabla$ and the Biot-Savart operator ${\mathcal B}$. Using a delicate but explicit computation, one arrives at
\begin{align*}
 (\partial_t + v_\ell \cdot\nabla) (z_i - z_\ell) 
 &= \Delta^{-1}\curl \div \left(\RR_\ell \right) + \Delta^{-1} \nabla \div \left( ((z_i - z_\ell) \cdot \nabla) v_\ell \right) \notag\\
 &\qquad + \Delta^{-1} \curl \div \left( ((z_i - z_\ell) \times \nabla)v_\ell + ((z_i - z_\ell) \times \nabla) v_i^T \right) \, .
\end{align*}
Since $ \Delta^{-1}\curl \div$, $\Delta^{-1} \nabla \div$, and $\Delta^{-1} \curl \div$ are bounded operators on $C^\alpha$, the bounds \eqref{eq:vector:potential:stability} follow from the above identity by using standard estimates for the transport equation (Gr\"onwall's inequality).

\subsubsection{The glued velocity and the glued stress}
Having constructed the exact solutions $(v_i,p_i) $, we glue them together using time-dependent cut-off functions $\chi_i$ as follows
\begin{align}
\label{eq:bar:v_q:def}
\overline v_q(x,t):= \sum_i \chi_i(t) v_i(x,t)\,, \qquad \mbox{and} \qquad \overline p_q^{(1)}(x,t) =\sum_i \chi_i(t) p_i(x,t)
 \, .
\end{align}
Here, the $\{\chi_i\}_i$ form a partition of unity in time for $[0,T]$ with the property that $\supp \chi_i\cap \supp \chi_{i+2}=\emptyset$ and moreover
\begin{subequations}\label{e:chi:time:width}
\begin{align}
\supp \chi_i&\subset [t_i-\sfrac{2\tau_q}{3},t_i+\sfrac{2\tau_q}{3}]  = I_{i-1} \cup J_i \cup I_i \, ,\\
\chi_i&=1\quad\textrm{ on } (t_i - \sfrac{\tau_q}{3}, t_i+ \sfrac{\tau_q}{3}) = J_i \, ,  \\
\norm{\partial_t^N \chi_i}_{C^0} &\lesssim \tau_q^{-N} \label{e:dt:chi}\,.
\end{align}
\end{subequations} 
Since the cutoffs $\chi_i$ only depend on time, we notice that the vector field $\overline v_q$ defined in in \eqref{eq:bar:v_q:def} is divergence-free. Moreover, by the definition of the cutoff functions, on every $J_i$ interval we have $\overline v_q = v_i$ and $\overline p_q^{(1)} = p_i$. Therefore, on $\cup_i J_i$ we have that $(\overline v_q, \overline p_q^{(1)} )$ is an exact solution of the Euler equations. On the other hand, on every $I_i$ interval we have
\begin{align*}
\overline v_q = \chi_i v_i + (1-\chi_i) v_{i+1} \, \qquad \mbox{and} \qquad \overline p_q^{(1)} = \chi_i p_i + (1-\chi_i) p_{i+1} \,,
\end{align*} 
which leads to 
\begin{align*}
\partial_t\overline v_q+\div(\overline v_q\otimes \overline v_q)+\nabla\overline p_q^{(1)} = 
\partial_t\chi_i(v_i-v_{i+1})-\chi_i(1-\chi_i)\div\left((v_i-v_{i+1})\otimes (v_i-v_{i+1})\right) \,.
\end{align*}
Using the elliptic inverse divergence operator $\RSZ$ from \eqref{eq:RSZ}, for all $t \in \cup_i I_i$ we thus define
\begin{subequations}
 \begin{align}
\mathring{\overline{R}}_q&=\partial_t\chi_i\mathcal{R}(v_i-v_{i+1})-\chi_i(1-\chi_i)(v_i-v_{i+1})\mathring{\otimes} (v_i-v_{i+1}) \, , \label{eq:bar:R_q:def}\\
\overline{p}_q&=\overline{p}_q^{(1)} -\chi_i(1-\chi_i)\left(|v_i-v_{i+1}|^2-\int_{\T^3}|v_i-v_{i+1}|^2\,dx\right) \, .
\end{align}
\end{subequations}
As discussed above, on $[0,T] \setminus \cup_i I_i = \cup_i J_i$ we have $\mathring{\overline{R}}_q = 0$ and let $\overline{p}_q=\overline{p}_q^{(1)}$. 
By construction, $\mathring{\overline{R}}_q$ is traceless symmetric, supported on the union of the $I_i's$, and we have that $(\overline v_q, \mathring{\overline{R}}_q, \overline{p}_q)$ solve the Euler-Reynolds system \eqref{e:euler_reynolds} on $\T^3 \times [0,T]$. 

Thus, it just remains to estimate the glued velocity field and Reynolds stress defined in \eqref{eq:bar:v_q:def} and \eqref{eq:bar:R_q:def}. For this purpose we note that upon rewriting $v_i - v_{i+1} = (v_i - v_\ell) - (v_{i+1} - v_\ell)$ we obtain from \eqref{eq:v_ell:all:bounds} and \eqref{eq:vi:minus:v_ell} that (cf.~\cite[Proposition 4.3]{BDLSV17})
\begin{subequations}
\label{eq:bar:v_q:all:bounds}
\begin{align}
 \norm{\bar v_q - v_{\ell}}_{C^\alpha} &\lesssim \delta_{q+1}^{\sfrac12}\ell^{\alpha} \label{e:vq:vell} \\
\norm{\overline{v}_q-v_\ell}_{C^{N+\alpha}} &\lesssim \tau_q\delta_{q+1}\ell^{-1-N+\alpha} \label{e:vq:vell:additional} \\
\norm{\bar v_q}_{C^{N+1}} &\lesssim \delta_{q}^{\sfrac12} \lambda_q \ell^{-N}\label{e:vq:1} \\
\left| \int_{\T^3} |\bar v_q|^2 - |v_\ell |^2 dx \right| &\lesssim \delta_{q+1}\ell^\alpha\label{e:voverline_vell_energy_diff} 
\end{align}
\end{subequations} 
holds for all $N\geq 0$. 
On the other hand, for the glued stress $\Rbar_q$ we obtain from \eqref{eq:vi:minus:v_ell} and \eqref{eq:vector:potential:stability}, upon rewriting $v_i - v_\ell$ as $\curl(z_i - z_\ell)$ and using that $\RSZ \curl$ is bounded on $C^\alpha$, that (cf.~\cite[Proposition 4.4]{BDLSV17})
\begin{subequations}
\label{eq:bar:R_q:all:bounds}
\begin{align}
 \norm{\mathring{\overline R_q}}_{C^{N+\alpha}} &\lesssim \delta_{q+1}\ell^{-N+\alpha} \label{e:Rq:1}\\
\norm{(\partial_t + \overline v_q\cdot \nabla) \mathring{\overline R_q}}_{C^{N+\alpha}} &\lesssim \delta_{q+1}\delta_q^{\sfrac12}\lambda_q\ell^{-N-\alpha}. \label{e:Rq:Dt}
\end{align}
\end{subequations}
For the estimate \eqref{e:Rq:Dt} above, we additionally need to commute the first order operator $v_\ell \cdot \nabla$ past the zero order operator $\RSZ \curl$. The bound \eqref{e:Rq:Dt} we use that this commutator is bounded on $C^\alpha$, and the operator norm of the commutator is bounded by a constant times $\norm{v_\ell}_{C^{1+\alpha}}$. A similar estimate in $C^{N+\alpha}$ also holds after applying the Leibniz rule.

\subsection{Mikado flows} 
\label{sec:Mikado}
In this section we recall the definition and the main properties of the Mikado flows constructed in~\cite{DSZ17}. Here we give a slightly different presentation of the  construction of~\cite{DSZ17}, which is consistent with our definition of intermittent jets in Section~\ref{s:intermittent}. The below geometric lemma is a variation of~\cite[Lemma 1]{Nash54} (cf.~\cite[Lemma 2.4]{DSZ17} and Proposition~\ref{p:split} above).
\begin{lemma}\label{l:linear_algebra} Denote by $\overline{B_{\sfrac 12}}(\Id)$ the closed ball of radius $1/2$ around the identity matrix, in the space of symmetric $3\times 3$ matrices.
 There exist mutually disjoint sets  $\{\Lambda_i \}_{i=0,1}  \subset\mathbb S^2\cap \mathbb Q^3$ such that for each $\xi \in \Lambda_i$ there exist $C^\infty$ smooth  functions $\gamma_{\xi}: B_{\sfrac 12}(\Id)\rightarrow \mathbb R$ which obey
\begin{align*}
R=\sum_{\xi\in \Lambda_{i}} \gamma^2_{\xi}(R)(\xi\otimes \xi)
\end{align*} 
for every symmetric matrix $R$ satisfying $\abs{R-\Id}\leq \sfrac12$, and for each $i \in \{0,1\}$.  

For a sufficiently large geometric constant $C_\Lambda \geq 1$, to be chosen precisely in Section~\ref{sec:Onsager:velocity:inductive} below, we define the constant
\begin{align}
M = C_\Lambda \sup_{\xi \in \Lambda_i} \left( \norm{\gamma_\xi}_{C^0} + \norm{\nabla \gamma_\xi}_{C^0} \right) \, ,
\label{eq:Onsager:M:def}
\end{align}
which appears in \eqref{e:v_q_inductive_est}. Moreover, for $i \in \{0,1\}$, and each $\xi\in \Lambda_i $, let use define $A_{\xi}\in\mathbb S^2\cap \mathbb Q^3$ to be an orthogonal vector to $\xi$. Then for each $\xi\in \Lambda_i$, we have that $\{\xi,A_{\xi},\xi\times A_{\xi}\}\subset \mathbb S^2\cap \mathbb Q^3$ form an orthonormal basis for $\mathbb R^3$. Furthermore, similarly to the constant $n_*$ of Proposition~\ref{p:split}, we label by $n_*$  the smallest natural such that 
\begin{equation}\label{e:Mignolet}
\left\{n_*\, \xi,~n_* \, A_{\xi},~n_* \, \xi\times A_{\xi}\right\}\subset  \Z^3\,
\end{equation}
for every $\xi \in \Lambda_i$ and for every $i \in \{0,1\}$. That is, $n_*$ is the $\rm{l.c.m.}$ of the denominators of the rational numbers $\xi, A_\xi$, and $\xi \times A_\xi$. 
\end{lemma}

For $\eps_\Lambda > 0 $, to be chosen later in terms of the set $\Lambda_i$, let $\Psi:\mathbb R^2\rightarrow \mathbb R^2$ be a $C^\infty$ smooth function with support contained in a ball of radius $\eps_\Lambda$ around the origin. We normalize $\Psi$  such that $\phi=- \Delta  \Psi$  obeys
\begin{equation}\label{eq:phi_normalize}
\int_{\R^2} \phi^2(x_1,x_2)\,dx_1dx_2= 4 \pi^2\,.
\end{equation}
Moreover, since $\supp \Psi, \phi \subset  \T^2$, we abuse notation and still denote by $\Psi, \phi$ the $\T^2$-periodized versions of $\Psi$ and $\phi$.   
Then, for any large  $\lambda \in \N$ and every $\xi \in \Lambda_i$, we introduce the functions
\begin{subequations}
\label{eq:phi:rotate}
\begin{align}
\Psi_{(\xi)}(x)&:=\Psi_{\xi, \lambda}(x):=\Psi(n_* \lambda (x-\alpha_{\xi})\cdot A_\xi,n_* \lambda (x-\alpha_\xi) \cdot (\xi\times A_{\xi}))\, ,\\
\phi_{(\xi)}(x)&:=\phi_{\xi,\lambda}(x):=\phi(n_* \lambda (x-\alpha_{\xi})\cdot A_\xi,n_* \lambda (x-\alpha_\xi)\cdot (\xi\times A_{\xi}))\,, 
\end{align}
\end{subequations}
where $\alpha_\xi \in \R^3$ are {\em shifts} whose purpose is to ensure that the functions $\{\Psi_{(\xi)} \}_{\xi \in \Lambda_i}$ have mutually disjoint support. Note that since $n_*  A_\xi$ and $n_* \xi \times A_\xi \in \Z^3$, and $\lambda \in \N$, the functions $\Psi_{(\xi)}$ and $\phi_{(\xi)}$ are $(\sfrac{\T}{\lambda})^3$-periodic. By construction we have that $\{\xi, A_\xi,\xi \times A_\xi\}$ are an orthonormal basis or $\R^3$, and hence $\xi \cdot \nabla \Psi_{(\xi)}(x) = \xi \cdot \nabla \phi_{(\xi)}(x) = 0$.  From the normalization of $\phi$ we have that $\fint_{\T^3} \phi_{(\xi)}^2 dx = 1$ and $\phi_{(\xi)}$ has zero mean on $(\sfrac{\T}{\lambda})^3$.  Since $\phi=-\Delta  \Psi$  we have that $(n_*\lambda)^2 \phi_{(\xi)} = - \Delta \Psi_{(\xi)}$.  Last, we emphasize that the existence of the shifts $\alpha_\xi$, which ensure that the supports of $\Psi_{(\xi)}$ are mutually disjoint for $\xi \in \Lambda_i$, is guaranteed by choosing $\eps_\Lambda$  sufficiently small solely in terms of the set $\Lambda_i$. Indeed, we can always ensure that the rational direction vectors in $\Lambda_i$ give (periodized) straight lines which do not intersect, when shifted by suitably chosen vectors $\alpha_\xi$.  

With this notation, the {\em Mikado flows} $W_{(\xi)} \colon \T^3 \to \R^3$ are defined as
\begin{align}
W_{(\xi)}(x) := W_{\xi,\lambda}(x) := \xi \, \phi_{(\xi)}(x) \, .
\label{eq:Mikado:def}
\end{align}
Since $\xi \cdot \nabla \phi_{(\xi)} = 0$, we immediately deduce that 
\begin{align}
 \div W_{(\xi)} = 0 \, \qquad \mbox{and} \qquad \div \left(W_{(\xi)} \otimes W_{(\xi)} \right) = 0 \,.
 \label{eq:Mikado:1}
\end{align}
Therefore, the Mikado flows are exact, smooth, pressureless solutions of the stationary 3D Euler equations. By construction, the functions $W_{(\xi)}$ have zero mean on $\T^3$ and are in fact $(\sfrac{\T}{\lambda})^3$-periodic. Moreover, by our choice of $\alpha_\xi$ we have that 
\begin{align}
W_{(\xi)} \otimes W_{(\xi')} \equiv 0 \qquad \mbox{whenever} \qquad \xi \neq \xi' \in   \Lambda_i \,,
\label{eq:Mikado:2}
\end{align}
for $i\in \{0,1\}$, and our normalization of $\phi_{(\xi)}$ ensures that 
\begin{align}
\fint_{\T^3} W_{(\xi)}(x) \otimes W_{(\xi)}(x)\,dx=\xi\otimes \xi\,. 
\label{eq:Mikado:3}
\end{align}
Lastly, using~\eqref{eq:Mikado:3}, the definition of the functions $\gamma_\xi$ in Lemma~\ref{l:linear_algebra} and the $L^2$ normalization of the functions $\phi_{(\xi)}$ we have that 
\begin{align}
 \sum_{\xi \in \Lambda_i} \gamma_{\xi}^2(R) \fint_{\T^3} W_{(\xi)}(x) \otimes W_{(\xi)}(x) dx = R \, ,
 \label{eq:Mikado:4}
\end{align}
for every $i \in \{0,1\}$ and any symmetric matrix $R \in \overline B_{\sfrac 12}(\Id)$.

We summarize the above properties \eqref{eq:Mikado:1}--\eqref{eq:Mikado:4} of the Mikado building blocks defined in \eqref{eq:Mikado:def} in the following result:
\begin{lemma}[Lemma 2.3,\,\cite{DSZ17}]
Given a symmetric matrix $R \in \overline{B}_{\sfrac 12}(\Id)$ and $\lambda \in \N$, the Mikado flow
\begin{align*}
{\mathcal W}(R,x) = \sum_{\xi \in \Lambda_i} \gamma_\xi(R) \, W_{\xi,\lambda}(x)
\end{align*}
obeys
\begin{align*}
\div {\mathcal W} = 0, \quad \div ( {\mathcal W} \otimes {\mathcal W} ) = 0, \quad \int_{\T^3} {\mathcal W} \, dx = 0, \quad \fint_{\T^3} {\mathcal W} \otimes {\mathcal W} \, dx = R.
\end{align*}
That is, ${\mathcal W}$ is a zero mean, presureless, solution of the stationary 3D Euler equations, which may be used to cancel the stress $R$.
\end{lemma}

\begin{figure}[h!]
\begin{center}
\includegraphics[width=0.4\textwidth]{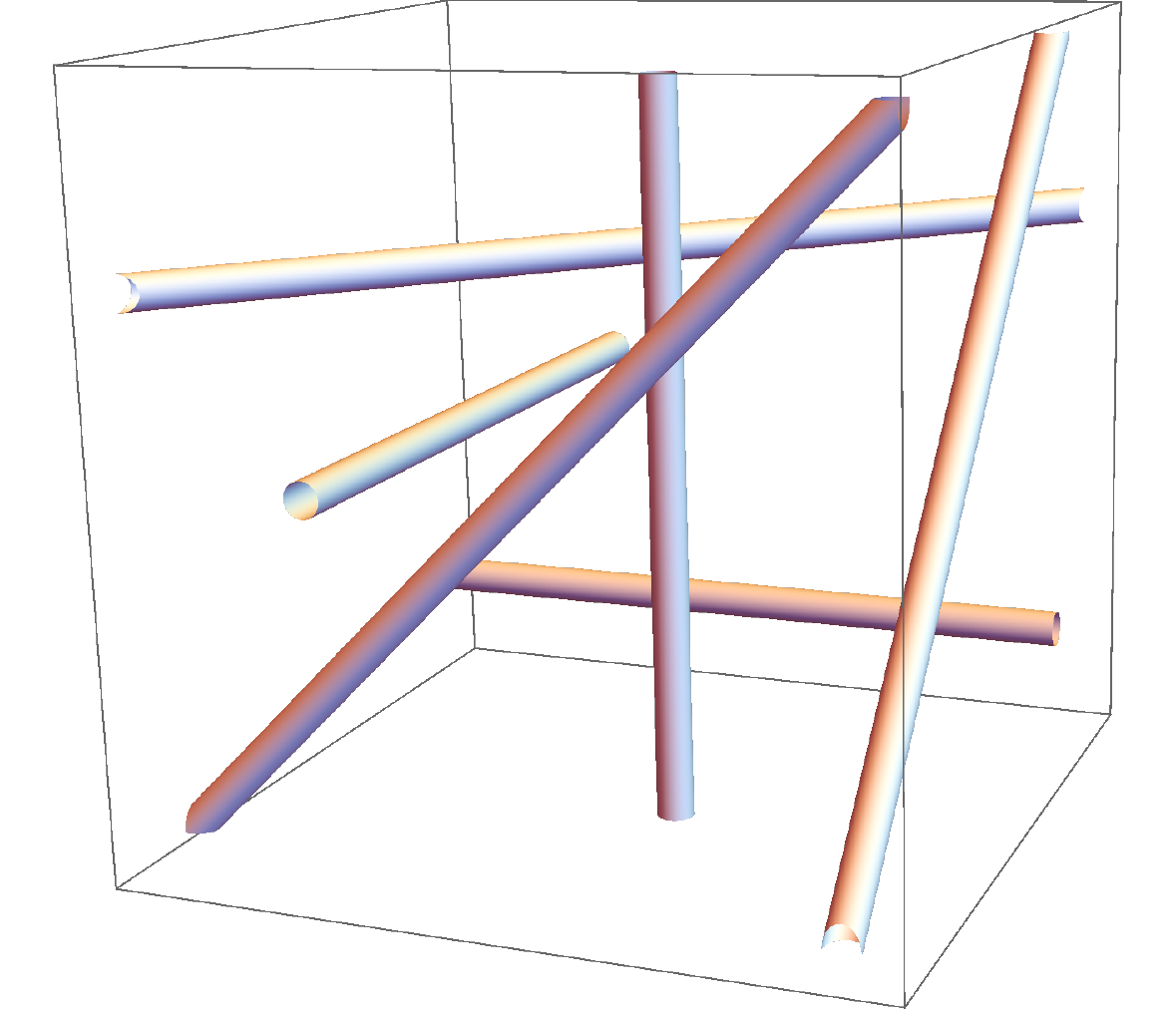}
\end{center}
\caption{{\small Example of a Mikado flow ${\mathcal W}$ restricted to one of the $(\sfrac{\T}{\lambda})^3$ periodic boxes.}}
\end{figure}

To conclude this section we note that $W_{(\xi)}$ may be written as the $\curl$ of a vector field, a fact which is useful in defining the incompressibility corrector in Section~\ref{sec:Onsager:principal:corrector}. Indeed, since $\xi \cdot \nabla \Phi_{(\xi)} = 0$, and since by definition we have that $-\frac{1}{(n_* \lambda)^2} \Delta \Phi_{(\xi)} = \phi_{(\xi)}$ we obtain
\begin{align}
\curl\left( \frac{1}{(n_* \lambda)^2} \nabla \Psi_{(\xi)} \times \xi \right)  = \curl\left( \frac{1}{(n_* \lambda)^2} \curl (\xi \Psi_{(\xi)}) \right) =    - \xi \left( \frac{1}{(n_* \lambda)^2} \Delta \Psi_{(\xi)} \right) = W_{(\xi)} \,.
\label{eq:Mikado:curl}
\end{align}
For notational simplicity, we   define 
\begin{align}
V_{(\xi)} =  \frac{1}{(n_* \lambda)^2} \nabla \Psi_{(\xi)} \times \xi 
\label{eq:Mikado:curl:2}
\end{align}
so that $\curl V_{(\xi)} = W_{(\xi)}$. With this notation we have the bounds
\begin{align}
\norm{W_{(\xi)}}_{C^N} + \lambda_{q+1} \norm{V_{(\xi)}}_{C^N} \les \lambda_{q+1}^{N}
\label{eq:Mikado:bounds}
\end{align}
for $N\geq 0$.

\subsection{The perturbation} 
\label{sec:Mikado:w}
In order to define the perturbation $w_{q+1} = v_{q+1} - \overline v_q$, we need to introduce a few objects.

\subsubsection{Cutoffs} 
\begin{figure}[h!]
\centering
  \begin{subfigure}[b]{0.45\textwidth}
    \includegraphics[width=0.97\textwidth]{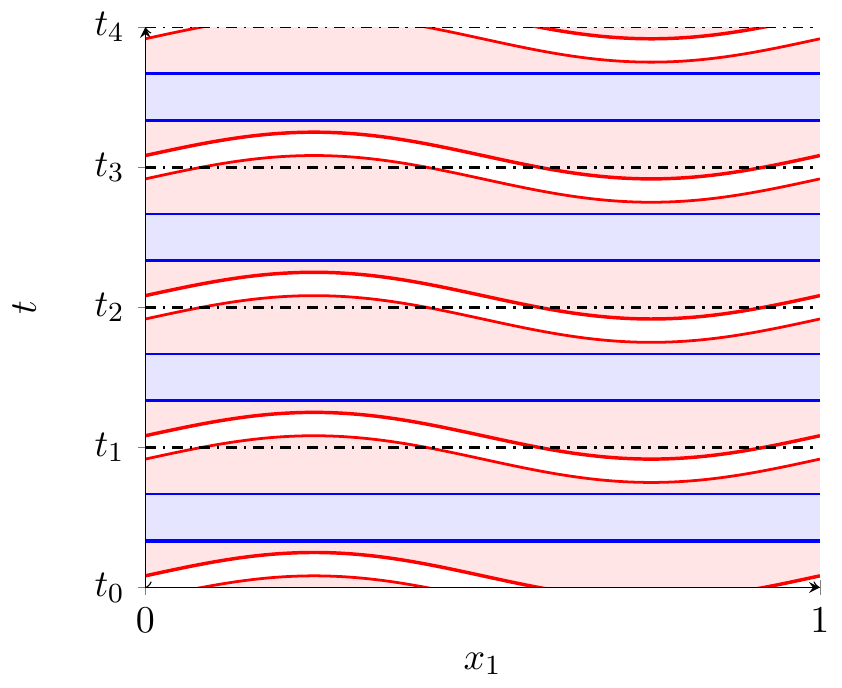}
    \caption{{\small The support of $\Rbar_q$ is given by the blue regions. The support of the cut-off functions $\eta_i$, which marks the region where the convex integration perturbation is supported, is given by the region between two consecutive red squiggling stripes.}}
    \label{fig:We}
  \end{subfigure}
    \begin{subfigure}[b]{0.05\textwidth}
    \,
  \end{subfigure}
  \begin{subfigure}[b]{0.45\textwidth}
    \includegraphics[width=1.07\textwidth]{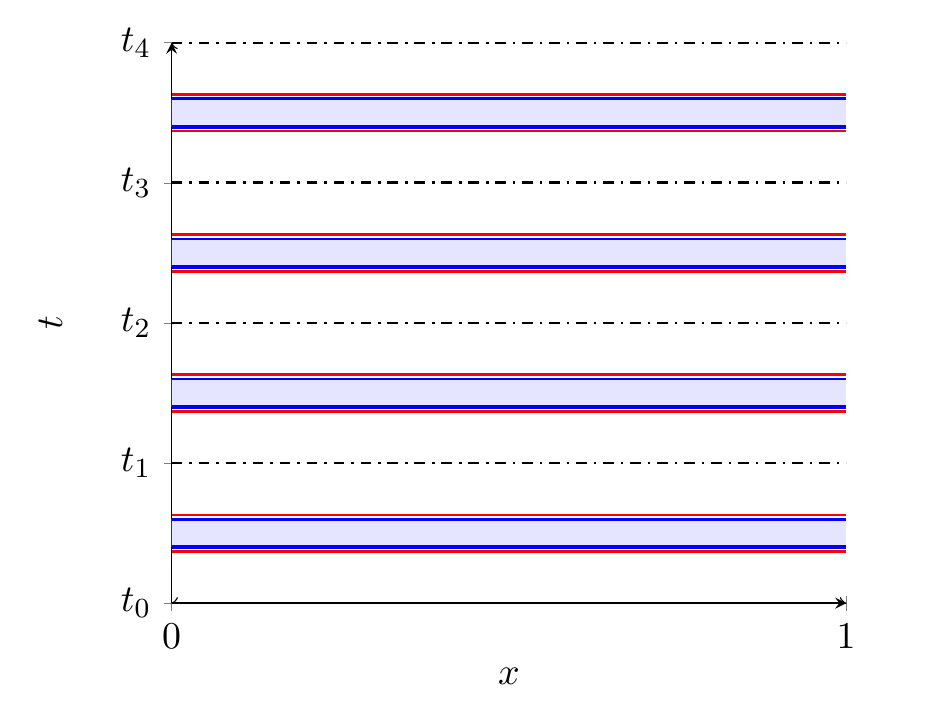}
    \caption{{\small In contrast, in the construction of~\cite{Isett16}, the support of the convex integration perturbation is (nearly) the same as the support of $\Rbar_q$. Consequently, on the time intervals between two consecutive blue regions, no energy is added to the solution.}}
    \label{fig:Phil}
  \end{subfigure}
  \begin{subfigure}[b]{0.1\textwidth}
  \end{subfigure}
  \caption{{\small The support of $\Rbar_q$ and $w_{q+1}$.}}
\end{figure}

Recall that $\Rbar_q$ has support in $\T^3 \times \cup_i I_i$, where $I_i$ is as defined in \eqref{eq:ti:Ii}. Accordingly, we define a family of cutoff functions $\{ \eta_i \}$ with the following properties 
\begin{enumerate}[(i)]
\item $\eta_i \in C^{\infty}_{x,t}$, $0 \leq \eta_i \leq 1$, and $\eta_i \, \eta_j \equiv 0$ for every  $i \neq j$ , 
\item $\eta_i \equiv 1$,  on $\T^3 \times I_i$ ,  
\item $\supp (\eta_i) \subset  \T^3 \times \tilde{I_i}$, where $\tilde I_i := J_i \cup I_i \cup J_{i+1} $, 
\item \label{eq:eta:i:special} $c_\eta \leq \sum_i \int_{\T^3} \eta_i^2(x,t) dx \leq 2 (2\pi)^3$,  for all $t\in [0,T]$,  where $c_\eta>0$ is a universal constant,  
\item $\norm{\partial_t^n \eta_i}_{C^m} \les_{n,m} \tau_q^{-n}$,  for all $n,m\geq 0$.
\end{enumerate}
The construction of such a sequence of cutoff functions is elementary (see, e.g.~\cite[Lemma 5.3]{BDLSV17}), and Figure~\ref{fig:We} shows how the supports of the $\eta_i$ relate to the support of $\Rbar_q$.  We emphasize that condition \eqref{eq:eta:i:special} is the one which allows us to alter the energy profile of $v_{q+1}$ even on time intervals where $\Rbar_q =0$.

We introduce the function $\rho_q(t)$,  which measures the remaining energy profile error after the gluing step, and after leaving ourselves room for adding a future velocity increment 
\begin{align*}
 \rho_q(t) = \frac 13 \left( e(t) - \frac{\delta_{q+2}}{2} - \int_{\T^3} |\overline v_q(x,t)|^2 dx \right)\,.
\end{align*}
By  \eqref{e:vq_vell_energy_diff} and \eqref{e:voverline_vell_energy_diff} we have 
$\abs{ \int_{\T^3} |v_q|^2  - |\overline v_q|^2} \les \delta_{q+1} \ell^\alpha $, which may be combined with \eqref{e:energy_inductive_assumption} and the choice of $\ell \leq \lambda_q^{-1 - \sfrac{3\alpha}{2}}$ to deduce that 
\begin{align}
\frac{\delta_{q+1}}{8 \lambda_q^\alpha} \leq  \rho_q(t)  \leq \delta_{q+1} 
\label{eq:rho:q:upper:lower}
\end{align}
for all $t \in [0,T]$. Similarly, using the assumed bound $\norm{\partial_t e}_{C^0} \les 1$, the energy inequality for the Euler-Reynolds system obeyed by $(\overline v_q,\Rbar_q)$, and the bounds \eqref{e:vq:1} and \eqref{e:Rq:1}, we obtain
\begin{align*}
\norm{\partial_t \rho_q}_{C^0} \les \delta_{q+1} \delta_q^{\sfrac 12} \lambda_q \, .
\end{align*}

The last cutoff function combines $\eta_i$ and $\rho_q$, and is defined by
\begin{align}
 \rho_{q,i}(x,t) = \rho_q(t) \, \frac{\eta_i^2(x,t)}{\sum_j \int_{\T^3} \eta_j^2 (y,t) dy}\,.
 \label{eq:rho:q:i:def}
\end{align}
By the normalization of $\rho_{q,i}$ we have that $\sum_i \int_{\T^3} \rho_{q,i}(x,t) = \rho_q(t)$ for all $t\in [0,T]$, and tracing back the properties of $\rho_q$ and $\eta_i$ we may verify that the following estimates hold
\begin{align}
\norm{\rho_{q,i}}_{C^0} \leq \frac{\delta_{q+1}}{c_\eta}, \quad \norm{\rho_{q,i}}_{C^N} \leq \delta_{q+1}, \quad \norm{\partial_t \rho_{q,i}}_{C^N} \les \delta_{q+1} \tau_q^{-1} \, ,
\label{eq:rho:q:i:bnd}
\end{align}
for all $N\geq 0$, where $c_\eta$ is the universal constant from property \eqref{eq:eta:i:special}.

\subsubsection{Flow maps}
Similarly to Section~\ref{sec:flow:maps:1}, we define the (backward) flow maps $\Phi_i$ for the velocity field $\overline v_{q}$ as the solution of the transport equation
\begin{subequations}
\label{eq:Phi:i:def}
\begin{align}
(\partial_t + \overline v_q  \cdot \nabla) \Phi_i &=0 \, ,\\  
\Phi_i\left(x,t_i\right) &= x \, ,
\end{align}
\end{subequations}
for all $t \in \supp (\eta_i) \subset \T^3 \times \tilde I_i$. For the reminder of this section, it is convenient to denote the material derivative as $D_{t,q}$, that is
\begin{align*}
 D_{t,q} = \partial_t + \overline v_q  \cdot \nabla_x \, .
\end{align*}
Since for every $t \in \tilde I_i$ we have $|t-t_i| \leq 2\tau_q$, by the definition of $\tau_q$ in \eqref{e:tau_def}, and using the estimates \eqref{e:vq:1}, we have that $\tau_q \norm{\nabla \overline v_q}_{C^0} \les \ell^{2\alpha} \ll 1$, and thus the CFL-condition is obeyed on $\tilde I_i$. From standard estimates for the transport equation (cf.~\cite[Proposition D.1]{BDLISZ15}), similarly to the bounds \eqref{eq:Phi:j:bounds} discussed earlier, we have 
\begin{subequations}
\label{eq:Phi:i:bnd}
\begin{align}
 \norm{\nabla \Phi_i(t) - \Id}_{C^0} &\les \ell^{2\alpha} \leq \frac 12
 \label{eq:Phi:i:bnd:a} \\
 \norm{\nabla \Phi_i}_{C^N} + \norm{(\nabla \Phi_i)^{-1}}_{C^N} &\les \ell^{-N}
 \label{eq:Phi:i:bnd:b} \\
 \norm{D_{t,q} \nabla \Phi_i}_{C^N} &\les \delta_q^{\sfrac 12} \lambda_q \ell^{-N}
 \label{eq:Phi:i:bnd:c}
\end{align}
\end{subequations}
for all $t\in \tilde I_i$ and $N\geq 0$. In order to establish these bounds it is useful to recall that after applying a gradient to \eqref{eq:Phi:i:def} we obtain the identity 
\begin{align*}
D_{t,q} \nabla \Phi_i = - \nabla \Phi_i \, D \overline v_q  \, .
\end{align*}
Also, we note from \eqref{eq:Phi:i:bnd:a} that $(\nabla \Phi_i)^{-1}$ is a well-defined object on $\tilde I_i$.

\subsubsection{Amplitudes}
Since $\eta_i \equiv 1$ on $\T^3 \times I_i$, $\eta_i \eta_j \equiv 0$ for $i\neq j$, and since $\supp (\Rbar_q) \subset \T^3 \times \cup_i I_i$, we have that 
\begin{align}
\sum_i \eta_i^2 \Rbar_q = \Rbar_q \, .
\label{eq:eta:partition}
\end{align}
Moreover, the cutoff functions $\eta_i$ already incorporate in them a temporal cutoff (recall that $\supp(\eta_i) \subset \T^3 \times \tilde I_i$), and thus it is convenient to introduce 
\begin{align*}
R_{q,i}   = \rho_{q,i}  \Id - \eta_i^2 \, \Rbar_q \,
\end{align*}
which is a stress supported in $\supp(\eta_i)$, and which obeys $\sum_i \RR_{q,i} = - \Rbar_q$.

For reasons which will become apparent only later (cf.~\eqref{eq:w:q+1:is:good}), we also define the symmetric tensor
\begin{align}
\tilde R_{q,i} = \frac{1}{\rho_{q,i}} \nabla \Phi_i \,  R_{q,i} \, \nabla \Phi_i^T = \Id + \left(\nabla \Phi_i \, \nabla \Phi_i^T - \Id\right) - \nabla \Phi_i\, \frac{\eta_i^2 \Rbar_q}{\rho_{q,i}} \, \nabla \Phi_i^T 
\label{eq:tilde:R:q:i:def}
\end{align}
for all $(x,t) \in \supp (\eta_i)$. By the above identity, estimate \eqref{e:Rq:1} with $N=0$, the property \eqref{eq:eta:i:special}, the bounds \eqref{eq:rho:q:upper:lower} and \eqref{eq:Phi:i:bnd:a}, we have that 
\begin{align*}
\norm{\tilde R_{q,i}(\cdot,t) - \Id }_{C^0} \les \ell^{\alpha} \leq \frac 12 \, \qquad \mbox{for all} \qquad t \in \tilde I_i \, ,
\end{align*}
once we ensure $a$ is taken to be sufficiently large.  Furthermore, using the estimates \eqref{eq:bar:R_q:all:bounds}, the properties of the $\eta_i$, \eqref{eq:rho:q:i:bnd} and \eqref{eq:Phi:i:bnd}, one may show that 
\begin{align}
 \norm{\tilde R_{q,i}}_{C^N} + \tau_q \norm{D_{t,q} \tilde R_{q,i}}_{C^N} \les \ell^{-N} \, \qquad \mbox{on} \qquad \supp(\eta_i) \, ,
 \label{eq:tilde:R:q:i:bnd}
\end{align}
for all $N\geq 0$.  One last important property of the stress $\tilde R_{q,i}$ is the identity
\begin{align}
 \sum_i \rho_{q,i} (\nabla \Phi_i)^{-1} \tilde R_{q,i} (\nabla \Phi_i)^{-T} = \left( \sum_i \rho_{q,i} \right) \Id - \Rbar_q \, ,
  \label{eq:tilde:R:q:i:identity}
\end{align}
which is useful in cancelling the glued stress. Here we have again appealed to \eqref{eq:eta:partition}. 

Thus, since $ \tilde R_{q,i}$ obeys the conditions of Lemma~\ref{l:linear_algebra} on $\supp (\eta_i)$, and since $\rho_{q,i}^{\sfrac 12}$ is a  multiple of $\eta_i$, we may define the amplitude functions 
\begin{align}
 a_{(\xi,i)}(x,t) = \rho_{q,i}(x,t)^{\sfrac 12} \,  \gamma_{\xi} (\tilde R_{q,i}) 
 \label{eq:a:xi:i:def}
\end{align}
where the $\gamma_\xi$ are the functions from Lemma~\ref{l:linear_algebra}. Note importantly that the amplitude functions already include a temporal cutoff, which shows that $\supp(a_{(\xi,i)}) \subset \supp(\eta_i)$. The amplitude functions $a_{(\xi)}$ inherit the expected $C^N$ bounds and material derivative bounds from \eqref{eq:rho:q:i:bnd}, \eqref{eq:tilde:R:q:i:bnd}, the product at the chain rules
\begin{align}
\norm{a_{(\xi,i)}}_{C^N} + \tau_q \norm{D_{t,q} a_{(\xi,i)}}_{C^N} \les \delta_{q+1}^{\sfrac 12} \ell^{-N}
\label{eq:Onsager:a:xi:CN}
\end{align}
for $N\geq 0$.

\subsubsection{Principal part of the velocity increment and the incompressibility corrector}
\label{sec:Onsager:principal:corrector}
For the remainder of the paper we consider Mikado building blocks as defined in \eqref{eq:Mikado:def} with $\lambda = \lambda_{q+1}$, i.e. 
\[ W_{(\xi)}(x)= W_{\xi,\lambda_{q+1}}(x) \, . \] 
Recall: for the index sets $\Lambda_i$ of Lemma~\ref{l:linear_algebra}, we overload notation and write $\Lambda_i = \Lambda_0$ for $i$ even, and $\Lambda_i = \Lambda_1$ for $i$ odd. With this notation, we now define the principal part of the velocity increment as 
\begin{align}
w_{q+1}^{(p)}(x,t) = \sum_{i} \sum_{\xi \in \Lambda_i} a_{(\xi,i)}(x,t) (\nabla \Phi_i(x,t))^{-1} W_{(\xi)}(\Phi_i(x,t)) \, .
\label{eq:Onsager:w:q+1:p}
\end{align}
When compared to the ansatz we made earlier for the $C^{0+}$ result (see~\eqref{e:w_xi}--\eqref{e:def_wo}), we notice the presence of $(\nabla \Phi_i)^{-1}$. The reason for this modification is as follows. At time $t=t_i$, we have $\Phi_i(x,t_i) = x$, $\nabla \Phi_i = \Id$, and by \eqref{eq:Mikado:1} we have that the vector field 
\begin{align*}
U_{i,\xi} = (\nabla \Phi_i)^{-1} W_{(\xi)}(\Phi_i)
\end{align*}
is incompressible at $t=t_i$. We then notice that $U_{i,\xi}$ is {\em Lie-advected} by the flow of the incompressible vector field $\overline v_q$, in the sense that 
\begin{align}
D_{t,q} U_{i,\xi} = (U_{i,\xi} \cdot \nabla) \overline v_q = (\nabla \overline v_q)^T U_{i,\xi} \, .
\label{eq:Onsager:Lie:advect}
\end{align}
This implies directly that $D_{t,q}( \div U_{i,\xi} ) = 0$, and thus the divergence free nature of $U_{i,\xi}$ is carried from $t=t_i$ to all $t$ close to $t_i$. This shows that the function $w_{q+1}^{(p)}$ defined in \eqref{eq:Onsager:w:q+1:p} is to leading order in $\lambda_{q+1}$ divergence-free (i.e.~the incompressibility corrector will turn out to be small).

At this stage we may also explain  why $R_{q,i}$ was not just normalized by $\rho_{q,i}$ but also conjugated with $\nabla \Phi_i$ respectively $(\nabla \Phi_i)^T$, in order to obtain $\tilde R_{q,i}$ (cf.~\eqref{eq:tilde:R:q:i:def}). Using the spanning property of the Mikado building blocks \eqref{eq:Mikado:4}, the fact that they have mutually disjoint support~\eqref{eq:Mikado:2}, identity \eqref{eq:tilde:R:q:i:identity} above, and the fact that the $\eta_i$ have mutually disjoint supports, we obtain
\begin{align}
&w_{q+1}^{(p)} \otimes w_{q+1}^{(p)} 
= \sum_i \sum_{\xi \in \Lambda_i} a_{(\xi,i)}^2 (\nabla \Phi_i)^{-1} \left( (W_{(\xi)}\circ \Phi_i) \otimes (W_{(\xi)}\circ \Phi_i) \right) (\nabla \Phi_i)^{-T} \notag\\
&= \sum_i \rho_{q,i}  (\nabla \Phi_i)^{-1}   \left( \sum_{\xi \in \Lambda_i}  \gamma_\xi^2(\tilde R_{q,i})  \left( ( W_{(\xi)} \otimes W_{(\xi)})\circ \Phi_i \right) \right) (\nabla \Phi_i)^{-T} \notag\\
&= \sum_i \rho_{q,i}  (\nabla \Phi_i)^{-1}  \tilde R_{q,i} (\nabla \Phi_i)^{-T}  + \sum_i \sum_{\xi \in \Lambda_i} a_{(\xi,i)}^2 (\nabla \Phi_i)^{-1} \left( \left( \Proj_{\neq 0}(W_{(\xi)} \otimes W_{(\xi)}) \right)\circ \Phi_i \right)  (\nabla \Phi_i)^{-T}   \notag\\
&= \left( \sum_i \rho_{q,i} \right) \Id - \Rbar_q + \sum_i \sum_{\xi \in \Lambda_i} a_{(\xi,i)}^2 (\nabla \Phi_i)^{-1} \left( \left(\Proj_{\geq \sfrac{\lambda_{q+1}}{2}}(W_{(\xi)} \otimes W_{(\xi)}) \right)\circ \Phi_i \right)  (\nabla \Phi_i)^{-T}   
\label{eq:w:q+1:is:good}
\end{align}
where we have denoted by $\Proj_{\neq 0} f(x) = f(x) - \fint_{\T^3} f(y) dy$, the projection of $f$ onto its nonzero frequencies. We have also used that since  $W_{(\xi)}\otimes W_{(\xi)}$ is $(\sfrac{\T}{\lambda_{q+1}})^3$-periodic, the identity $\Proj_{\neq 0} (W_{(\xi)} \otimes W_{(\xi)}) = \Proj_{\geq \sfrac{\lambda_{q+1}}{2}} (W_{(\xi)} \otimes W_{(\xi)})$ holds. The calculation \eqref{eq:w:q+1:is:good} shows that by design, the low frequency part of $w_{q+1}^{(p)} \otimes w_{q+1}^{(p)}$ cancels the glued stress $\Rbar_q$, modulo a multiple of the identity, which is then used to correct the energy profile and which contributes a pressure term to the equation.

Based on the definition \eqref{eq:Onsager:w:q+1:p}  of the principal part of the velocity increment, we construct an incompressibility corrector. As was observed in~\cite{DSZ17}, for any smooth vector field $V$, we have the identity
\begin{align*}
(\nabla \Phi_i)^{-1} \left( (\curl V)\circ \Phi_i \right)= \curl\left( (\nabla \Phi_i)^T (V\circ \Phi_i)\right).
\end{align*}
Recalling identity \eqref{eq:Mikado:curl} and the definition \eqref{eq:Mikado:curl:2}, we may write $W_{(\xi)} = \curl V_{(\xi)}$ and thus the above identity shows that 
\begin{align*}
(\nabla \Phi_i)^{-1} (W_{(\xi)} \circ \Phi_i) = \curl\left( (\nabla \Phi_i)^T (V_{(\xi)}\circ \Phi_i)\right) \,.
\end{align*}
From the above identity and \eqref{eq:Onsager:w:q+1:p}, it follows that if we define the incompressibility corrector as
\begin{align}
w_{q+1}^{(c)}(x,t) = \sum_i \sum_{\xi \in \Lambda_i} \nabla a_{(\xi,i)}(x,t) \times \left( (\nabla \Phi_i(x,t))^T (V_{(\xi)} (\Phi_i(x,t) )  \right)
\label{eq:Onsager:w:q+1:c}
\end{align}
then the total velocity increment $w_{q+1}$ obeys
\begin{align}
w_{q+1} = w_{q+1}^{(p)} + w_{q+1}^{(c)} = \curl \left( \sum_i \sum_{\xi \in \Lambda_i} a_{(\xi,i)} \, (\nabla \Phi_i)^T (V_{(\xi)}\circ \Phi_i) \right) 
\label{eq:Onsager:w:q+1}
\end{align}
so that it is automatically incompressible.

\subsubsection{Velocity inductive estimates}
\label{sec:Onsager:velocity:inductive}
The velocity field at level $q+1$ is constructed as
\begin{align}
v_{q+1} = \overline v_q + w_{q+1} = v_q + (v_\ell - v_q) + (\overline v_q - v_\ell) + w_{q+1} \,.
\label{eq:Onsager:v:q+1:def}
\end{align}
From \eqref{eq:Onsager:M:def}, \eqref{eq:rho:q:i:bnd}, \eqref{eq:Phi:i:bnd:a}, \eqref{eq:a:xi:i:def}, and \eqref{eq:Onsager:w:q+1:p}, and the fact that the $\eta_i$ have disjoint supports, once $a$ is sufficiently large we obtain that 
\begin{subequations}
\begin{align}
\norm{w_{q+1}^{(p)}}_{C^0} &\leq  \frac{2|\Lambda_i| \norm{\phi}_{C^0} }{c_\eta^{\sfrac 12} C_\Lambda}  M \delta_{q+1}^{\sfrac 12} \leq \frac{M}{8} \delta_{q+1}^{\sfrac 12} 
\label{eq:Liverpool:1}\\
 \norm{w_{q+1}^{(p)}}_{C^1} &\leq  \frac{4 |\Lambda_i| n_* \norm{\phi}_{C^1}}{c_\eta^{\sfrac 12}C_\Lambda}  M \delta_{q+1}^{\sfrac 12} \lambda_{q+1} \leq \frac{M}{8} \delta_{q+1}^{\sfrac 12} \lambda_{q+1}
\end{align}
\end{subequations}
by choosing the parameter $C_\Lambda$ from \eqref{eq:Onsager:M:def} to be large enough. Note that $C_\Lambda$ only depends on the cardinality of $\Lambda_i$, on the universal constant $c_\eta$, the geometric integer $n_*$, and on the $C^1$ norm of the function $\phi$, which in turn depends solely on the geometric constant $\eps_\Lambda$.

For the incompressibility corrector we lose a factor of $\ell^{-1}$ from the gradient landing on $a_{(\xi,i)}$, but we gain a factor of $\lambda_{q+1}$ because we have $V_{(\xi)}$ instead of $W_{(\xi)}$ (recall \eqref{eq:Mikado:curl:2}). Therefore, we may show that 
\begin{align}
\norm{w_{q+1}^{(c)}}_{C^0} + \frac{1}{\lambda_{q+1}} \norm{w_{q+1}^{(c)}}_{C^1} \les \delta_{q+1}^{\sfrac 12} \frac{\ell^{-1}}{\lambda_{q+1}} \,.
\label{eq:Liverpool:2}
\end{align}
We note that by choosing $\alpha$ to be sufficiently small in therms of $b$ and $\beta$, we have
\begin{align}
\frac{\ell^{-1}}{\lambda_{q+1}} = \frac{\delta_q^{\sfrac 12} \lambda_q^{1+ \sfrac{3\alpha}{2}}}{\delta_{q+1}^{\sfrac 12} \lambda_{q+1}} = \frac{\lambda_q^{1-\beta + \sfrac{3\alpha}{2}}}{\lambda_{q+1}^{1-\beta}} \leq 2 \lambda_{q}^{\sfrac{3\alpha}{2} - (b-1)(1-\beta)} \leq \lambda_{q}^{- \sfrac{(b-1)(1-\beta)}{2}}  \ll 1 \, ,
\label{eq:Onsager:ell:gap}
\end{align}
and thus by choosing $a$ sufficiently large we may ensure that the velocity increment defined in \eqref{eq:Onsager:w:q+1:c} obeys
\begin{align*}
\norm{w_{q+1}}_{C^0} + \frac{1}{\lambda_{q+1}} \norm{w_{q+1}}_{C^1} \leq \frac{M}{2} \delta_{q+1}^{\sfrac 12} \,.
\end{align*} 
By combining the above estimate with \eqref{eq:v_ell:all:bounds}, \eqref{eq:bar:v_q:all:bounds}, and \eqref{eq:Onsager:v:q+1:def}, after choosing $a$ sufficiently large, we deduce that \eqref{e:v_diff_prop_est} is satisfied, and moreover that the bounds \eqref{e:v_q_inductive_est} and \eqref{e:v_q_0} hold with $q$ replaced with $q+1$.

\subsection{Reynolds stress} 
Recall that the pair $(\overline v_q,\Rbar_q)$ solves the Euler-Reynolds system \eqref{e:euler_reynolds}, and that $v_{q+1}$ is defined in \eqref{eq:Onsager:v:q+1:def}. In this subsection we define the new Reynolds stress $\RR_{q+1}$, and show that it obeys the estimate
\begin{align}
\norm{\RR_{q+1}}_{C^\alpha} \les \frac{\delta_{q+1}^{\sfrac 12} \delta_{q}^{\sfrac 12} \lambda_q}{\lambda_{q+1}^{1-4\alpha}} \, .
\label{eq:Onsager:R:q+1:to:do}
\end{align}
The above bound immediately implies the desired estimate \eqref{e:R_q_inductive_est} at level $q+1$, upon noting that the following parameter inequality holds (after taking $\alpha$ sufficiently small and $a$ sufficiently large)
\begin{align}
\frac{\delta_{q+1}^{\sfrac 12} \delta_{q}^{\sfrac 12} \lambda_q}{\lambda_{q+1}^{1-4\alpha}} \leq \frac{\delta_{q+2}}{\lambda_{q+1}^{4\alpha}} \,.
\label{eq:parameter:stuffing}
\end{align}
The remaining power of  $\lambda_{q+1}^{-\alpha}$ is used to absorb the implicit constant in \eqref{eq:Onsager:R:q+1:to:do}. 

In order to define $\RR_{q+1}$, similarly to \eqref{eq:R:q+1:A} we write
\begin{align}
\div \RR_{q+1} - \nabla p_{q+1}
&=\underbrace{ D_{t,q} w^{(p)}_{q+1}}_{\div (R_{\rm transport})}+ \underbrace{\div(w_{q+1}^{(p)} \otimes w_{q+1}^{(p)} + \Rbar_q)}_{\div (R_{\rm oscillation}) + \nabla p_{\rm oscillation}}+ \underbrace{w_{q+1}\cdot  \nabla \overline v_{q}}_{\div( R_{\rm Nash})}   \notag \\
& \qquad +\underbrace{ D_{t,q} w^{(c)}_{q+1}+ \div\left(w_{q+1}^{(c)} \otimes w_{q+1}+ w_{q+1}^{(p)} \otimes w_{q+1}^{(c)}\right) }_{\div (R_{\rm corrector})+ \nabla p_{\rm corrector}}   - \nabla  \overline p_q \, .
\label{eq:Onsager:R:q+1}
\end{align}
The various traceless symmetric stresses present implicitly in \eqref{eq:Onsager:R:q+1} are defined using the inverse divergence operator $\RSZ$, and by recalling the identity \eqref{eq:w:q+1:is:good} (for the oscillation error) as
\begin{subequations}
\begin{align}
R_{\rm transport} &= \RSZ \left( D_{t,q} w^{(p)}_{q+1} \right) \label{eq:Onsager:R:transport}\\
R_{\rm oscillation} &= \sum_i \sum_{\xi \in \Lambda_i} \RSZ \div \left( a_{(\xi,i)}^2 (\nabla \Phi_i)^{-1} \left( \left(\Proj_{\geq \sfrac{\lambda_{q+1}}{2}}(W_{(\xi)} \otimes W_{(\xi)}) \right)\circ \Phi_i \right)  (\nabla \Phi_i)^{-T}   \right) \label{eq:Onsager:R:oscillation} \\
R_{\rm Nash} &= \RSZ \left( w_{q+1}\cdot  \nabla \overline v_{q} \right) \label{eq:Onsager:R:Nash}\\
R_{\rm corrector} &= \RSZ\left( D_{t,q} w^{(c)}_{q+1} \right)+  \left(w_{q+1}^{(c)} \mathring \otimes w_{q+1}^{(c)} + w_{q+1}^{(c)} \mathring \otimes w_{q+1}^{(p)} + w_{q+1}^{(p)} \mathring\otimes w_{q+1}^{(c)} \right) \label{eq:Onsager:R:corrector}
\end{align}
\end{subequations}
while the pressure terms are given by $p_{\rm oscillation} = \sum_i \rho_{q,i}$ and $p_{\rm corrector}  =  2 w_{q+1}^{(c)} \cdot w_{q+1}^{(p)} + |w_{q+1}^{(c)}|^2$. With this notation we have $p_{q+1} = \overline p_q - p_{\rm oscillation} - p_{\rm corrector}$ and  
\begin{align}
\RR_{q+1} = R_{\rm transport} + R_{\rm oscillation} + R_{\rm Nash} + R_{\rm corrector}   \, .
\label{eq:R:Onsager:q+1:B}
\end{align}
Prior to estimating the above stresses, it is convenient to adapt the stationary phase bounds of Section~\ref{sec:RSZ:bounds} from Beltrami flows to Mikado flows.

\subsubsection{Inverse divergence and stationary phase bounds}
In order to apply Lemma~\ref{lem:stationary:phase}, and obtain  bounds similar to \eqref{eq:stationary:phase:out} and \eqref{eq:stationary:phase:double}, but for Mikado flows instead of the Beltrami flows, we decompose the function $\phi_{(\xi)}$ which defined $W_{(\xi)}$ in \eqref{eq:Mikado:def} as a Fourier series. Recall that $\phi_{(\xi)}$ defined in \eqref{eq:phi:rotate} is $(\sfrac{\T}{\lambda_{q+1}})^3$ periodic and has zero mean. Additionally, the function $\phi$ is $C^\infty$ smooth. Therefore, we may decompose
\begin{align}
\phi_{(\xi)}(x) = \phi_{\xi,\lambda_{q+1}} (x) = \sum_{k \in \Z^3 \setminus \{0\}} f_\xi(k) e^{i \lambda_{q+1} k \cdot (x-\alpha_\xi)} \,
\label{eq:Mikado:Fourier}
\end{align}
where the complex numbers $f_{\xi}(k)$ are the Fourier series coefficients of the $C^\infty$ smooth, mean-zero $\T^3$ periodic function $z\mapsto \phi(n_* z\cdot A_\xi, n_* z\cdot (\xi \times A_\xi))$. The  shift $x\mapsto x-\alpha_\xi$ has no effect on the estimates. Moreover, the Fourier coefficients decay arbitrarily fast. For any $m\in \N$ we have $|f_\xi(k)| = |f_\xi(k) e^{i \lambda_{q+1} k \cdot  \alpha_\xi}| \leq C |k|^{-m}$, where the constant $C$ depends on $m$ and on geometric parameters of the construction, such as $n_*$, the sets $\Lambda_i$, the shifts $\alpha_\xi$, and norms of the bump function $\phi(x_1,x_2)$. Thus, $C$ is independent of $\lambda_{q+1}$, or any other $q$-dependent parameter. 

A similar Fourier series decomposition applies to the function $\frac{1}{n_*\lambda_{q+1}} \nabla \Psi_{(\xi)} = (\nabla \Psi)_{(\xi)}$ which is used in \eqref{eq:Mikado:curl:2} to define $V_{(\xi)}$. For this function we also obtain that its Fourier series coefficients decay arbitrarily fast, with constants that are bounded independently of $q$ (and hence $\lambda_{q+1}$). 

Therefore, for a smooth function $a(x,t)$, in order to estimate $\RSZ(a \, W_{(\xi)} \circ \Phi_i)$, we use identity \eqref{eq:Mikado:Fourier}, and apply Lemma~\ref{lem:stationary:phase} for each $k$ individually, and then sum in $k$ using the fast decay of the Fourier coefficients $f_\xi(k)$. Without giving all the details, we summarize this procedure as follows.  Let $a \in C^0([0,T];C^{m,\alpha}(\T^3) )$ be such that  $\supp (a) \subset \supp (\eta_i)$, which ensures that the phase $\Phi_i$ obeys the conditions of Lemma~\ref{lem:stationary:phase} by \eqref{eq:Phi:i:bnd:a}. Using \eqref{eq:Phi:i:bnd:b} we obtain from Lemma~\ref{lem:stationary:phase} that 
\begin{align}
\norm{\RSZ \left( a \; (W_{(\xi)} \circ \Phi_i ) \right) }_{C^\alpha}  + \lambda_{q+1} \norm{\RSZ \left( a \; (V_{(\xi)} \circ \Phi_i ) \right) }_{C^\alpha}
\les \frac{  \norm{a}_{C^0}}{\lambda_{q+1}^{1-\alpha}} +   \frac{ \norm{a}_{C^{m,\alpha}}+\norm{a}_{C^0} \ell^{-m-\alpha}}{\lambda_{q+1}^{m-\alpha}} \, ,
\label{eq:Mikado:stationary:phase:1}
\end{align}
where the implicit constant is independent of $q$. 

Recalling that $W_{(\xi)} \otimes W_{(\xi)} = (\xi \otimes \xi) \phi_{(\xi)}^2$, 
and using that the function $ \Proj_{\geq \sfrac{\lambda_{q+1}}{2}}\phi_{(\xi)}^2$ is also zero mean $(\sfrac{\T}{\lambda_{q+1}})^3$-periodic, a similar argument shows that 
\begin{align}
\norm{\RSZ \left( a \; \left( \Big( \Proj_{\geq \sfrac{\lambda_{q+1}}{2}}(W_{(\xi)} \otimes W_{(\xi)} ) \Big) \circ\Phi_i \right) \right) }_{C^\alpha}   
\les \frac{  \norm{a}_{C^0}}{\lambda_{q+1}^{1-\alpha}} +   \frac{ \norm{a}_{C^{m,\alpha}}+\norm{a}_{C^0} \ell^{-m-\alpha}}{\lambda_{q+1}^{m-\alpha}} 
\label{eq:Mikado:stationary:phase:2}
\end{align}
holds. The above estimate is useful for estimating the oscillation error.

\subsubsection{The estimate for $\RR_{q+1}$}
In this section we show that the  stresses defined in \eqref{eq:R:Onsager:q+1:B} obey \eqref{eq:Onsager:R:q+1:to:do}. 
The Nash error and the corrector error are in a sense lower order, and they can be treated similarly (or using similar bounds) to the transport and oscillation errors. Because of this, we omit the details for estimating $R_{\rm Nash}$ and $R_{\rm corrector}$.  

{\bf Transport error. \,}
Recalling the definition of $w_{q+1}^{(p)}$ in \eqref{eq:Onsager:w:q+1:p}, and the  Lie-advection identity \eqref{eq:Onsager:Lie:advect} we obtain that the transport stress in \eqref{eq:Onsager:R:transport} is given by  
\begin{align}
R_{\rm transport} &= \sum_{i} \sum_{\xi \in \Lambda_i} \RSZ \left( a_{(\xi,i)} (\nabla \overline v_q)^T (\nabla \Phi_i)^{-1} W_{(\xi)}(\Phi_i) \right)  +  \RSZ \left( \left( D_{t,q} a_{(\xi,i)} \right)  (\nabla \Phi_i)^{-1} W_{(\xi)}(\Phi_i) \right) \,.
\label{eq:Onsager:transport:1}
\end{align}
In order to bound the terms in \eqref{eq:Onsager:transport:1} we use \eqref{eq:Mikado:stationary:phase:1} to gain a factor of $\lambda_{q+1}^{-1+\alpha}$ from the operator $\RSZ$ acting on the highest frequency term $W_{(\xi)}\circ \Phi_i$. The derivatives of $a_{(\xi,i)}$, $\nabla \overline v_q$, and $(\nabla \Phi_i)^{-1}$ are estimated using \eqref{eq:Onsager:a:xi:CN}, \eqref{e:vq:1}, and \eqref{eq:Phi:i:bnd:b} respectively. These bounds show that each additional spacial derivatives costs a power of $\ell^{-1}$. We obtain from \eqref{eq:Mikado:stationary:phase:1} that
\begin{align*}
\norm{R_{\rm transport}}_{C^\alpha} &\les \frac{\delta_{q+1}^{\sfrac 12} \delta_q^{\sfrac 12} \lambda_q}{\lambda_{q+1}^{1-\alpha}} \left( 1 + \frac{\ell^{-m-\alpha}}{\lambda_{q+1}^{m-1}} \right)  +  \frac{\delta_{q+1}^{\sfrac 12} \tau_q^{-1}}{\lambda_{q+1}^{1-\alpha}} \left( 1 + \frac{\ell^{-m-\alpha}}{\lambda_{q+1}^{m-1}} \right)    \, .
\end{align*}
Recalling \eqref{eq:Onsager:ell:gap}, we have that $(\ell \lambda_{q+1})^{-1} \leq \lambda_q^{-\sfrac{(b-1)(1-\beta)}{2}}$, and thus upon taking the parameter $m$ in  to be sufficiently large (in terms of $\beta$ and $b$), we obtain that $R_{\rm transport}$ indeed is bounded by the right side of \eqref{eq:Onsager:R:q+1:to:do}, as desired. 

{\bf Oscillation error. \,}
For the oscillation error, which is defined in \eqref{eq:Onsager:R:oscillation}, the main observation is that when the $\div$ operator lands on the highest frequency term, namely $\left(\Proj_{\geq \sfrac{\lambda_{q+1}}{2}}(W_{(\xi)} \otimes W_{(\xi)}) \right)\circ \Phi_i$, due to certain cancellations this term vanishes. Since by construction we have $(\xi \cdot \nabla) \phi_{(\xi)} = 0$ it also follows that  $(\xi \cdot \nabla) \Proj_{\geq \sfrac{\lambda_{q+1}}{2}}( \phi_{(\xi)}^2 )  = 0$. Therefore,
\begin{align*}
& \div \left( a_{(\xi,i)}^2 (\nabla \Phi_i)^{-1} \left( \left(\Proj_{\geq \sfrac{\lambda_{q+1}}{2}}(W_{(\xi)} \otimes W_{(\xi)}) \right)\circ \Phi_i \right)  (\nabla \Phi_i)^{-T}   \right)  
\notag\\
&= \div \left( a_{(\xi,i)}^2 (\nabla \Phi_i)^{-1} (\xi \otimes \xi) (\nabla \Phi_i)^{-T} \left( \left(\Proj_{\geq \sfrac{\lambda_{q+1}}{2}}(\phi_{(\xi)}^2) \right)\circ \Phi_i \right)     \right)  \notag\\
&=\left( \left(\Proj_{\geq \sfrac{\lambda_{q+1}}{2}}(\phi_{(\xi)}^2) \right)\circ \Phi_i \right)    \div \left( a_{(\xi,i)}^2 (\nabla \Phi_i)^{-1}    (\xi \otimes \xi) (\nabla \Phi_i)^{-T}\right)    \notag\\
&\qquad +  \underbrace{a_{(\xi,i)}^2 (\nabla \Phi_i)^{-1} (\xi \otimes \xi) (\nabla \Phi_i)^{-T} \left((\nabla \Phi_i)^{T} \left(\nabla \Proj_{\geq \sfrac{\lambda_{q+1}}{2}}(\phi_{(\xi)}^2) \right)\circ \Phi_i \right)  }_{=0} .
\end{align*}
The above identity shows that 
\[
R_{\rm oscillation} = \sum_i \sum_{\xi \in \Lambda_i} \RSZ \left( \left( \left(\Proj_{\geq \sfrac{\lambda_{q+1}}{2}}(W_{(\xi)} \otimes W_{(\xi)}) \right)\circ \Phi_i \right) \div \left( a_{(\xi,i)}^2 (\nabla \Phi_i)^{-1} (\xi\otimes \xi) (\nabla \Phi_i)^{-T}  \right)  \right)    \, ,
\]
at which point we may appeal to the stationary phase estimate \eqref{eq:Mikado:stationary:phase:2} combined with the bounds \eqref{eq:Onsager:a:xi:CN}  and \eqref{eq:Phi:i:bnd:b} to obtain
\begin{align*}
\norm{R_{\rm oscillation}}_{C^\alpha} \les \frac{\delta_{q+1}\ell^{-1}}{\lambda_{q+1}^{1-\alpha}}  \left( 1 + \frac{\ell^{-m-\alpha}}{\lambda_{q+1}^{m-1}} \right) \les \frac{\delta_{q+1}^{\sfrac 12} \delta_q^{\sfrac 12}\lambda_q }{\lambda_{q+1}^{1-\sfrac{5\alpha}{2}}}    \, . 
\end{align*}
Here we have again taken $m$ sufficiently large, and have recalled the definition of $\ell$ in \eqref{e:ell_def}. Thus the oscillation error is also bounded by the right side of \eqref{eq:Onsager:R:q+1:to:do}, as claimed.

\subsubsection{Energy increment} 
To conclude the proof of Proposition~\ref{p:main}, it remains to show that \eqref{e:energy_inductive_assumption} holds with $q$ replaced by $q+1$. In order to prove this bound we show that 
\begin{align}
\left| e(t) - \int_{\T^3} |v_{q+1}(x,t)|^2 dx - \frac{\delta_{q+2}}{2} \right| \les \frac{\delta_{q+1}^{\sfrac 12} \delta_q^{\sfrac 12} \lambda_{q}^{1+2\alpha}}{\lambda_{q+1}} 
\label{eq:Onsager:energy:increment}
\end{align}
holds. Recalling the parameter estimate \eqref{eq:parameter:stuffing}, and taking $a$ sufficiently large to absorb all the implicit constants, it is clear that \eqref{eq:Onsager:energy:increment} implies the bound \eqref{e:energy_inductive_assumption} at level $q+1$. 

In order to prove \eqref{eq:Onsager:energy:increment}, the principal observation is the following. Taking the trace of \eqref{eq:w:q+1:is:good}, since $\Rbar_q$ is traceless we obtain
\begin{align*}
\int_{\T^3} |w_{q+1}|^2 dx 
&= 3 \sum_i \int_{\T^3} {\rho_{q,i}} dx \notag\\
&  + \sum_i \sum_{\xi \in \Lambda_i} \int_{\T^3} a_{(\xi,i)}^2 \tr\left( (\nabla \Phi_i)^{-1} (\xi \otimes \xi) (\nabla \Phi_i)^{-T} \right) \left( \left(\Proj_{\geq \sfrac{\lambda_{q+1}}{2}}(W_{(\xi)} \otimes W_{(\xi)}) \right)\circ \Phi_i \right)  dx \,.
\end{align*}
The second term in the above identity can be made arbitrarily small, since it is the $L^2$ inner product of a function whose oscillation frequency is $\les \ell^{-1}$ (cf.~\eqref{eq:Onsager:a:xi:CN}  and \eqref{eq:Phi:i:bnd:b}) and a function which is $\lambda_{q+1}$ periodic and zero mean. On the other hand, by the design of the functions $\rho_{q,i}$ (cf.~\eqref{eq:rho:q:i:def}) we have 
\begin{align*}
3 \sum_i \int_{\T^3} {\rho_{q,i}} dx = 3 \rho_q(t) = e(t) - \frac{\delta_{q+2}}{2} - \int_{\T^3} |\overline v_q(x,t)|^2 dx \, .
\end{align*}
Since $v_{q+1} = \overline v_q + w_{q+1}$, the above identity implies that 
\begin{align*}
e(t) - \int_{\T^3} |v_{q+1}(x,t)|^2 dx - \frac{\delta_{q+2}}{2} = - 2 \int_{\T^3} \overline v_q \cdot w_{q+1} dx - 2 \int_{\T^3} w_{q+1}^{(p)} \cdot w_{q+1}^{(c)} dx - \int_{\T^3} |w_{q+1}^{(c)}|dx \,.
\end{align*}
The corrector terms in the above give estimates consistent with \eqref{eq:Onsager:energy:increment} by appealing to \eqref{eq:Liverpool:1}, \eqref{eq:Liverpool:2}, and \eqref{eq:Onsager:ell:gap}. For the first term on the right side of the above we recall (cf.~\eqref{eq:Onsager:w:q+1}) that $w_{q+1}$ may be written as the $\curl$ of a vector field whose size is $\delta_{q+1}^{\sfrac 12} \lambda_{q+1}^{-1}$. Integrating by parts the $\curl$ and using \eqref{e:vq:1} with $N=0$ we conclude the proof of \eqref{eq:Onsager:energy:increment}, and hence of   Proposition~\ref{p:main}.

\section{Navier-Stokes: existence of weak solutions with finite energy}
\label{sec:NSE:L2}
We consider weak solutions of the Navier-Stokes equations \eqref{eq:NSE} (see Definition~\ref{d:weak_sol} above).
The viscosity parameter obeys $\nu \in (0,1]$, and we consider solutions $v$ which have zero mean on $\T^3$. The main result of~\cite{BV} is Theorem~\ref{thm:BV:main}.

In order to keep the exposition as simple as possible, without omitting any of the main ideas involved in proving Theorem~\ref{thm:BV:main}, in this section we prove a simpler result which states that there exists a weak solution in the aforementioned regularity class, whose kinetic energy is {\em not monotone decreasing}. Thus, this weak solution is not equal to any of the Leray-Hopf weak solutions arising from the same $L^2$ initial datum $v|_{t=0}$. 

\begin{theorem}
\label{thm:NSE:main}
There exists $\beta > 0$, such that the following holds.  There exists a sufficiently small $\nu \in (0,1]$ and a weak solution $v$
of the Navier-Stokes equations \eqref{eq:NSE}, which lies in $ C^0([0,1];H^\beta(\T^3)) \cap C^0([0,1];W^{1,1+\beta}(\T^3))$, and such that $\norm{v(\cdot,1)}_{L^2} \geq 2 \norm{v(\cdot,0)}_{L^2}$. 
\end{theorem}
Theorem~\ref{thm:NSE:main} is proven in Section~\ref{sec:NSE:proof} below. This proof can be used to also establish  Theorem~\ref{thm:BV:main} if one adds a few inductive estimates to the the list in \eqref{eq:inductive:NSE} below, cf.~(2.4)--(2.6) in \cite{BV}. For simplicity we omit these details here, and only prove this more restrictive result. 

\subsection{Inductive estimates}
Let $(v_q,\RR_q)$ be a given solution of the Navier-Stokes-Reynolds system~\eqref{e:Navier-Reynolds}. 
We consider the same parameters $\lambda_q \to + \infty$ and $\delta_q \to 0^+$  defined by
\begin{subequations}
\label{eq:parameters:def}
\begin{align}
\lambda_q &= 2\pi a^{(b^q)} \\
\delta_q &= \lambda_q^{-2\beta} \,.
\end{align}
\end{subequations}
The sufficiently large (universal) parameter $b$ is free, and so is the sufficiently small parameter $\beta = \beta(b)$. The parameter $a$ is   chosen to be a sufficiently large multiple of the geometric constant $n_*$ from~\eqref{e:Mignolet} above. 

For $q\geq 0$ we make the following inductive assumptions
\begin{subequations}
\label{eq:inductive:NSE}
\begin{align}
\norm{v_q}_{L^2} &\leq 1 - \delta_q^{\sfrac 12} 
\label{eq:vq:L2:NSE}\\
\norm{v_q}_{C^{1}_{x,t}} &\leq  \lambda_q^4
\label{eq:vq:C1:NSE}\\
\norm{\RR_{q}}_{L^{1}} &\leq c_R \delta_{q+1}
\label{eq:Rq:L1:NSE}
\end{align} 
\end{subequations}
where $c_R>0$ is a sufficiently small universal constant (determined in \eqref{eq:a:xi:L2} in terms of the  parameter $M$ from Lemma~\ref{l:linear_algebra} below).  
The inductive proposition is almost identical to~Proposition~\ref{prop:iteration:Euler}.

\begin{proposition}[Main iteration]
\label{prop:iteration:NSE}
There exists a sufficiently large universal parameter $b>1$, and a sufficiently small parameter $\beta = \beta(b) >0 $, such that the following holds.  There exists a sufficiently large constant $a_0 = a_0(c_R, \beta,b)$ such that for any $a\geq a_0$ which is a multiple of the geometric constant $n_*$, there exist functions $(v_{q+1},\RR_{q+1})$ which solve \eqref{e:Navier-Reynolds} and obey  \eqref{eq:inductive:NSE}  at level $q+1$, such that 
\begin{align}
\norm{v_{q+1} - v_q }_{L^2} \leq  \delta_{q+1}^{\sfrac 12} 
\label{eq:increment:L2}
\end{align}
holds. All parameters in the Proposition, $b$, $\beta$, and $a_0$, are independent of $\nu \in (0,1]$ .
\end{proposition}

\begin{remark}[Quantifying the parameters $\beta$ and $b$] 
\label{rem:b:beta:NSE}
For the purpose of specifying the parameters $b$ and $\beta$ it is convenient to first fix an auxiliary parameter $0 < \alpha \ll 1$.
It is sufficient to take $\alpha \leq \sfrac{1}{240}$, arbitrary. Then, inspecting the proof of Proposition~\ref{prop:iteration:NSE}, we may verify that $b \in \N$ may be chosen to be any multiple of $7$ which obeys $b > \sfrac{4}{\alpha}$. With such a value of $b$ fixed, $\beta$ may be chosen to obey $\beta < \sfrac{\alpha}{b}$. We did not try here to optimize these conditions.  The particular choices $\alpha = \sfrac{1}{240}$, $b = 1001$, and $\beta = \sfrac{1}{2^{18}}$ are permissible.
\end{remark}

\subsection{Proof of Theorem~\ref{thm:NSE:main}}
\label{sec:NSE:proof}
Fix the parameters $b, \beta, c_R$ and  $a_0$ from Proposition~\ref{prop:iteration:NSE}.
By possibly enlarging the value of $a\geq a_0$, we may ensure that $\delta_{0} \leq \sfrac 14$.

We define an incompressible, zero mean vector field $v_0$ by
\[
v_0(x,t) = \frac{t }{(2\pi)^{\sfrac 32}} ( \sin(\lambda_0 x_3), 0, 0) 
\]
and define the kinematic viscosity 
\[
\nu =  \lambda_0^{-2} \in (0,1] \, .
\]
Note that by construction we have 
$\sup_{t\in [0,1]} \norm{v_0(\cdot,t)}_{L^2} \leq \norm{v_0(\cdot,1)}_{L^2} = \sfrac{1}{\sqrt{2}} \leq 1- \delta_0^{\sfrac 12}$,  so that \eqref{eq:vq:L2:NSE} is automatically satisfied. Moreover, $\norm{v_0}_{C^{1}_{x,t}} \leq \lambda_0 \leq \lambda_0^4$ since $\lambda_0 = 2\pi a \geq 1$ by construction.

The vector field $v_0$ defined above is a shear flow, and thus $v_0 \cdot \nabla v_0 = 0$. Thus, it obeys \eqref{e:Navier-Reynolds} at $q=0$, with stress $\RR_0$ defined by 
\begin{align}
 \RR_0 =   - \frac{1 + \nu t \lambda_0^2}{\lambda_0 (2\pi)^{\sfrac 32}}  \left( {\begin{array}{ccc}
   0 & 0 &  \cos(\lambda_0^{\sfrac 12} x_3) \\
   0 & 0 & 0 \\
    \cos(\lambda_0^{\sfrac 12} x_3) & 0 & 0   
  \end{array} } \right)
 \, . 
 \label{eq:R0:NSE:def}
\end{align}
Recalling that $\nu = \lambda_0^{-2}$, we see that for some universal constant $C>0$ we have
\begin{align*}
 \norm{\RR_0}_{L^1} \leq \frac{C}{\lambda_0}  \leq \delta_1 \, .
\end{align*}
The last inequality above uses that $\lambda_0 \delta_1 = (2\pi)^{1-2\beta} a^{1-2\beta b} \geq a^{\sfrac 12} \geq C$. This inequality holds because $\beta b \leq \sfrac 14$ (see Remark~\ref{rem:b:beta:NSE} above), and $a$ can be taken to be larger than $C^2$, which is a universal constant.  Thus, condition \eqref{eq:Rq:L1:NSE} is also obeyed for $q=0$. 
 
We may thus use the iteration Proposition~\ref{prop:iteration:NSE} and obtain a sequence of solutions $(v_q,\RR_q)$ which obey \eqref{eq:inductive:NSE} and \eqref{eq:increment:L2}. By interpolation we have that for any $\beta' \in (0, \frac{\beta}{4+\beta})$, the following series is   summable
\begin{align*}
\sum_{q\geq 0} \norm{v_{q+1}-v_q}_{H^{\beta'}} 
\les \sum_{q\geq 0} \norm{v_{q+1}-v_q}_{L^2}^{1-\beta'}  \norm{v_{q+1}-v_q}_{H^{1}}^{\beta'}
\les  \sum_{q\geq 0} \delta_{q+1}^{\frac{1-\beta'}{2}} \lambda_{q+1}^{4 \beta'}  \les  \sum_{q\geq 0} \lambda_{q+1}^{-\beta (1-\beta') + 4\beta'}  \les 1
\end{align*}
where the implicit constant is universal. Thus, we may define a limiting function $v = \lim_{q\to \infty} v_q$ which lies in $C^0([0,1];H^{\beta'})$. Moreover, $v$ is a weak solution of the Navier-Stokes equation \eqref{eq:NSE}, since by \eqref{eq:Rq:L1:NSE} we have that $\lim_{q\to \infty}\RR_q = 0$ in $C^0([0,1];L^1)$. From the maximal regularity of the heat equation we also obtain that for some $\beta'' \in (0,\beta')$ we have $v \in C^0([0,1];W^{1,1+\beta''}(\T^3))$. The regularity of the weak solution claimed in Theorem~\ref{thm:NSE:main} then holds with $\beta$ replaced by $\beta''>0$.

It remains to show that $\norm{v(\cdot,1)}_{L^2} \geq 2 \norm{v(\cdot,0)}_{L^2}$. For this purpose note that since $b^{q+1} \geq b (q+1)$, we have
\begin{align*}
\norm{v - v_0}_{L^2} \leq \sum_{q\geq 0} \norm{v_{q+1} - v_q}_{L^2} 
\leq   \sum_{q\geq 0} \delta_{q+1}^{\sfrac 12}  
\leq   \sum_{q\geq 0} a^{-\beta (b^{q+1})} \leq   \sum_{q\geq 0} (a^{-\beta b})^{q+1}  = \frac{  a^{-\beta b}}{1-a^{-\beta b}} \leq \frac{1}{6}
\end{align*}
once we choose $a$ sufficiently large, in terms of $\beta$, and $b$. Using that by construction we have $\norm{v_0(\cdot,0)}_{L^2} = 0$, and $\norm{v_0(\cdot,1)}_{L^2} = \sfrac{1}{\sqrt{2}}$, we obtain that 
\begin{align*}
2 \norm{v(\cdot,0)}_{L^2} &\leq 2 \norm{v_0(\cdot,0)}_{L^2} + 2\norm{v(\cdot,0) - v_0(\cdot,0)}_{L^2}  \notag\\
&\leq \frac{1}{3}  \leq \frac{1}{\sqrt{2}} - \frac{1}{6} \leq \norm{v_0(\cdot,1)}_{L^2} - \norm{v(\cdot,1) - v_0(\cdot,1)}_{L^2} \leq \norm{v(\cdot,1)}_{L^2}
\end{align*}
holds.  This concludes the proof of Theorem~\ref{thm:NSE:main}.

\subsection{Mollification}
\label{sec:mollify:NSE}

Similarly to the Euler section~\ref{sec:mollify}, we replace the pair $(v_q,\RR_q)$ by a mollified pair $(v_\ell,\RR_\ell)$ defined exactly as in \eqref{eq:v:R:ell:def}  by
\begin{align*}
v_{\ell} = (v_q \ast_x \phi_{\ell}) \ast_t \varphi_{\ell} \, , \qquad \mbox{and} \qquad 
\RR_{\ell} = (\RR_q \ast_x \phi_{\ell}) \ast_t \varphi_{\ell} \, .
\end{align*}
similarly, to \eqref{e:euler_reynolds_ell}, we obtain that $(v_\ell,\RR_\ell)$ obey
\begin{subequations}
\label{e:NSE_reynolds_ell}
\begin{align}
\partial_t  v_\ell + \div(v_\ell \otimes v_\ell) + \nabla p_\ell - \nu \Delta v_{\ell}
&= \div \Big(  \RR_\ell +  R_{\rm commutator} \Big)\, ,  \\
\div v_\ell &= 0 \, ,
\end{align}
\end{subequations}
where  traceless symmetric commutator stress $\tilde R_{\rm commutator}$ is given by
\begin{align}
R_{\rm commutator} &= (v_\ell \mathring \otimes v_\ell) - ((v_q \mathring \otimes v_q)\ast_x \phi_{\ell}) \ast_t \varphi_{\ell} \, .\label{eq:NSE:R:commutator} 
\end{align}
Note that this definition is the same as the one for Euler in \eqref{e:euler_reynolds_ell:a}, since the Laplacian commutes with mollification.

The parameter $\ell$ has to be chosen similarly to the Euler case, e.g.~so that it obeys
\begin{align}
\ell \lambda_q^4  \leq \lambda_{q+1}^{-\alpha} \qquad \mbox{and} \qquad \ell^{-1} \leq \lambda_{q+1}^{2\alpha}.
\label{eq:ell:cond:NSE}
\end{align}
for $0 < \alpha \ll 1$ as in Remark~\ref{rem:b:beta:NSE}.  This choice is permitted because $\alpha b > 4$. 
In particular, we may define $\ell$ as the geometric mean of the two bounds imposed by \eqref{eq:ell:cond:NSE}
\begin{align}
\ell = \lambda_{q+1}^{- \sfrac{3 \alpha}{2}} \lambda_q^{-2} \,.
\label{eq:ell:def:NSE}
\end{align}
With this choice for $\ell$, the same arguments which gave us \eqref{eq:Rc:bound}--\eqref{eq:v:ell:C0} in the case of Euler, may be used in conjunction with the inductive assumptions \eqref{eq:inductive:NSE} to yield  bounds for the commutator stress and the mollified velocity. We have
\begin{align}
\norm{\RR_{\rm commutator}}_{L^1} \les \norm{\RR_{\rm commutator}}_{C^0} \les \ell \norm{v_q}_{C^1_{x,t}} \norm{v_q}_{C^0_{x,t}}  \les \ell   \lambda_q^4 \les  \lambda_{q+1}^{-\alpha} \, .\label{eq:Rc:bound:L1}
\end{align}
For the mollified velocity we have the bounds 
\begin{align}
 \norm{v_q - v_\ell}_{L^2} &\les  \norm{v_q - v_\ell}_{C^0} \les  \ell \norm{v_q}_{C^1}\les \ell \lambda_q^4 \les \lambda_{q+1}^{-  \alpha}   \ll \delta_{q+1}^{\sfrac12}
 \label{eq:V_q_ell_est:L2} \, ,\\
  \norm{v_\ell}_{L^2} &\leq \norm{v_q}_{L^2} \leq 1 - \delta_q^{\sfrac 12} 
  \label{eq:V_ell_est:L2} \, , \\
    \norm{v_\ell}_{C^N_{x,t}}&  \les \ell^{-N+1} \norm{v_q}_{C^1_{x,t}} \les \lambda_q^4  \ell^{-N+1} \les  \lambda_{q+1}^{- \alpha} \ell^{-N} \label{eq:V_ell_est:C1} \,,
\end{align}
where $N\geq 1$ and we have taken $\beta$ sufficiently small, in terms of $b$ and $\alpha$.  

\subsection{Intermittent jets}\label{s:intermittent}

Recall the geometric Lemma~\ref{l:linear_algebra} discussed earlier. For the Navier-Stokes construction we do not require two sets of wave-vectors $\Lambda_0$ and $\Lambda_1$, and instead choose just one of them, which we label as $\Lambda$. Moreover, we consider the constant $C_\Lambda$ in the definition of $M$ to equal $C_\Lambda = 8 |\Lambda| (1+ 8 \pi^3)^{\sfrac 12}$, where $|\Lambda|$ is the cardinality of the set $\Lambda$. We also recall the vectors $A_\xi \in\mathbb S^2\cap \mathbb Q^3$ and the constant $n_*$ introduced in Lemma~\ref{l:linear_algebra}.

Let $\Phi:\mathbb R^2\rightarrow \mathbb R^2$ be a smooth function with support contained in a ball of radius $1$. We normalize $\Phi$  such that $\phi=-\Delta  \Phi$ obeys
\begin{equation}\label{e:phi_normalize}
 \frac{1}{4\pi^2}\int_{\R^2} \phi^2(x_1,x_2)\,dx_1dx_2=1\,.
\end{equation}
We remark that by definition $\phi$ has mean zero. Define $\psi:\mathbb R\rightarrow \mathbb R$ to be a smooth, mean zero function with support in the ball of radius $1$ satisfying
\begin{equation}\label{e:psi_normalize}
\frac{1}{2\pi}\int_{\R} \psi^2(x_3)\,dx_3=1\,.
\end{equation}
For parameters $ r_{\perp}, r_{\|}>0$  such that
\[ 
r_{\perp}\ll  r_{\|}\ll 1\,,
\]
we define $\phi_{ r_{\perp}}$, $\Phi_{ r_{\perp}}$ and $\psi_{ r_{\|}}$ to be the rescaled cutoff functions
\begin{align}
\phi_{ r_{\perp}}(x_1,x_2)=\frac{1}{ r_{\perp}} \phi\left(\frac{x_1}{ r_{\perp}},\frac{x_2}{ r_{\perp}}\right),\quad \Phi_{ r_{\perp}}(x_1,x_2)=\frac{1}{ r_{\perp}} \Phi\left(\frac{x_1}{ r_{\perp}},\frac{x_2}{ r_{\perp}}\right),\quad 
\psi_{ r_{\|}}(x_3)=\frac{1}{ r_{\|}^{\sfrac12}} \psi\left(\frac{x_3}{ r_{\|}}\right) \, .
\label{eq:phi:psi:normalize}
\end{align}
With this rescaling we have $\phi_{ r_{\perp}}=- r_{\perp}^2\Delta \Phi_{ r_{\perp}}$, the functions $\phi_{r_\perp}$ and $\Phi_{r_\perp}$ are supported in the ball of radius $r_\perp$ in $\R^2$, $\psi_{r_\|}$ is supported in the ball of radius $r_\|$ in $\R$, and we maintain the normalizations $\norm{\phi_{r_\perp}}_{L^2}^2 = 4\pi$ and $\norm{\psi_{r_\|}}_{L^2}^2 = 2\pi$.
Lastly, by an abuse of notation, we {\em periodize} $\phi_{r_\perp}$, $\Phi_{ r_{\perp}}$, and $\psi_{ r_{\|}}$ so that the functions are treated as periodic functions defined on $\T^2$, $\T^2$ and $\T$ respectively. These periodic functions (rescaled and tilted version of them) form the building blocks for our intermittent jets, which are defined next.

Consider a large \emph{real number} $\lambda$ such that $\lambda r_{\perp}\in \mathbb N$, and a large {\em time oscillation parameter} $\mu>0$.  For every $\xi \in \Lambda$ we introduce the shorthand notation
\begin{subequations}
\label{eq:phi:psi:rotate}
\begin{align}
\psi_{(\xi)}(x,t)&:=\psi_{\xi, r_{\perp}, r_{\|},\lambda,\mu} (x,t) :=\psi_{ r_{\|}}(n_* r_{\perp}\lambda (x\cdot \xi+\mu t)),
\label{eq:psi:rotate}\\
\Phi_{(\xi)}(x)&:=\Phi_{\xi, r_{\perp},\lambda}(x):=\Phi_{ r_{\perp}}(n_* r_{\perp}\lambda (x-\alpha_{\xi})\cdot A_\xi,n_* r_{\perp}\lambda (x-\alpha_\xi) \cdot (\xi\times A_{\xi}))\\
\phi_{(\xi)}(x)&:=\phi_{\xi, r_{\perp},\lambda}(x):=\phi_{ r_{\perp}}(n_* r_{\perp}\lambda (x-\alpha_{\xi})\cdot A_\xi,n_* r_{\perp}\lambda (x-\alpha_\xi)\cdot (\xi\times A_{\xi}))\,.
\end{align}
\end{subequations}
where $\alpha_{\xi}\in\mathbb R^3$ are shifts which ensure that the functions $\{ \Phi_{(\xi)} \}_{\xi \in \Lambda}$  have mutually disjoint support. In order for such shifts $\alpha_{\xi}$ to exist, it is sufficient to assume that $ r_{\perp}$ smaller than a universal constant, which depending only on the geometry of the finite set  $\Lambda $. It is important to note that by \eqref{eq:phi:psi:normalize}, the function $\psi_{(\xi)}$ oscillates at frequency proportional to $r_\perp r_\|^{-1} \lambda \ll \lambda$, whereas $\phi_{(\xi)}$ and $\Phi_{(\xi)}$ oscillate at frequency proportional to $\lambda$.

With this notation, the {\em intermittent jets} $W_{(\xi)} \colon \T^3 \times \R \to \R^3$ are then defined as 
\begin{align}
W_{(\xi)}(x,t)
:=W_{\xi, r_{\perp}, r_{\|},\lambda,\mu}(x,t)  
:= \xi \, \psi_{(\xi)}(x,t) \, \phi_{(\xi)}(x)\,.
\label{eq:jet}
\end{align}
In view of \eqref{e:Mignolet} and the condition $ r_\perp \lambda \in \N$, we have that $W_{(\xi)}$ has zero mean on $\T^3$ and is $\left(\sfrac{\T}{ r_{\perp} \lambda}\right)^3$ periodic. Moreover, by our choice of $\alpha_\xi$, we have that
\begin{align}
W_{(\xi)} \otimes W_{(\xi')} \equiv 0 \qquad \mbox{whenever} \qquad \xi \neq \xi' \in   \Lambda \,,
\label{eq:Messi}
\end{align}
i.e.~the $\{ W_{(\xi)} \}_{\xi \in \Lambda}$ have mutually disjoint support. As a consequence of the normalizations \eqref{e:phi_normalize} and \eqref{e:psi_normalize}, the rescaling \eqref{eq:phi:psi:normalize}, the fact that translations and orthogonal transformations are volume preserving, and the fact that rescaling a periodic function does not alter its $L^p$ norms, we have
\begin{align*}
\fint_{\T^3} W_{(\xi)}(x,t) \otimes W_{(\xi)}(x,t)\,dx=\xi\otimes \xi\,. 
\end{align*}
As a consequence, using Lemma \ref{l:linear_algebra} and the cancelation \eqref{eq:Messi}, we have that
\begin{equation}\label{e:Reynolds_cancellation}
 \sum_{\xi   \in \Lambda }\gamma_{\xi}^2(R)  \fint_{\T^3}  W_{(\xi)}(x,t)\otimes   W_{(\xi)}(x,t)\,dx=R\,,
\end{equation}
for every symmetric matrix $R$ satisfying $\abs{R-\Id}\leq \sfrac12$.

The essential identity obeyed by the intermittent jets is
\begin{align}\label{eq:useful:2}
\div\left(W_{(\xi)}\otimes W_{(\xi)}\right) =2(W_{(\xi)}\cdot \nabla \psi_{(\xi)})\phi_{(\xi)}\xi=\frac{1}{\mu}\phi^2_{(\xi)}\partial_t \psi^2_{(\xi)}\xi = \frac{1}{\mu} \partial_t \left(\phi^2_{(\xi)}  \psi^2_{(\xi)}\xi\right)\,.
\end{align}
Identity \eqref{eq:useful:2} follows from the fact that by construction we have that $W_{(\xi)}$ is a scalar multiple of $\xi$, by \eqref{eq:psi:rotate} we have
\[
(\xi \cdot \nabla) \psi_{(\xi)} = \frac{1}{\mu} \partial_t \psi_{(\xi)}\,,
\]
and lastly, because $\phi_{(\xi)}$ is time-independent. 

\begin{figure}[h!]
\begin{center}
\includegraphics[width=0.4\textwidth]{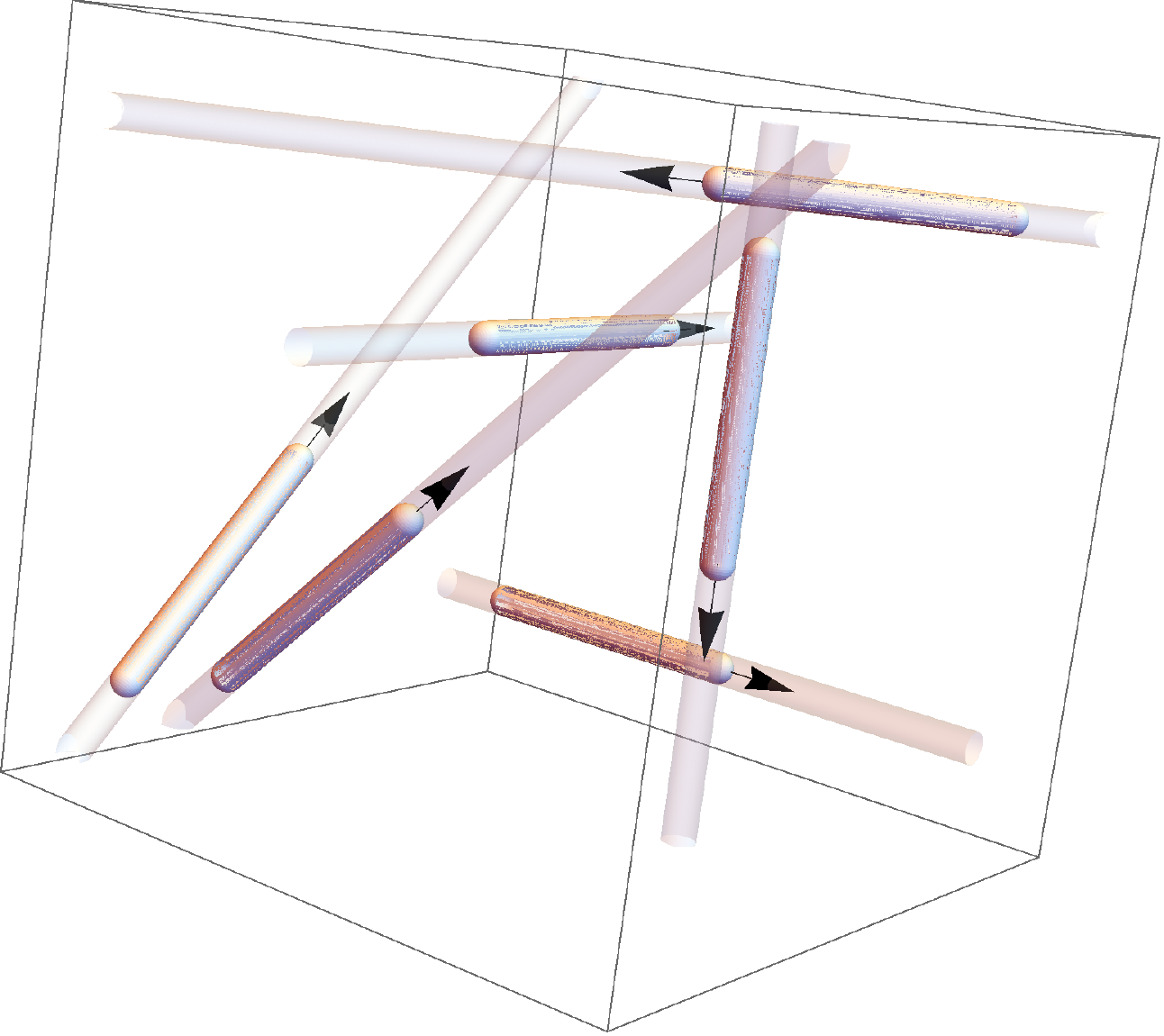}
\end{center}
\caption{{\small The intermittent jets in one of the $(\sfrac{\T}{r_\perp \lambda})^3$ periodic boxes.}}
\end{figure}
Lastly, we note that the intermittent jets $W_{(\xi)}$ {\em are not divergence free}, however assuming $ r_{\perp}\ll  r_{\|}$ they may be corrected by a small term, such that the sum with the corrector is divergence free. To see this, let 
us  define $V_{(\xi)}: \T^3\times \R \rightarrow \R^3$ by
\begin{align*}
V_{(\xi)}(x,t) :=V_{\xi, r_{\perp}, r_{\|},\lambda,\mu}(x,t) := \frac{1}{ n_*^2 \lambda^2} \xi \, \psi_{(\xi)}(x,t) \, \Phi_{(\xi)}(x) \, ,
\end{align*}
where we still use the notation from \eqref{eq:phi:psi:rotate}. A computation then shows that 
\begin{align}\label{e:div_correct_id}
\curl \curl V_{(\xi)} -  W_{(\xi)} = \underbrace{\frac{1}{ n_*^2 \lambda^2}\curl\left(\Phi_{(\xi)}\curl\left(\psi_{(\xi)}\xi\right)\right)}_{\equiv 0}+ \frac{1}{ n_*^2 \lambda^2}
\nabla\psi_{(\xi)}\times \curl\left(\Phi_{(\xi)}\xi\right) =:   W_{(\xi)}^{(c)}\, .
\end{align}
Thus, by the definition of $W_{(\xi)}^{(c)}$, we have
\[\div\left(W_{(\xi)}+W^{(c)}_{(\xi)}\right)\equiv 0\,.\]
Moreover, since we have $r_{\perp}\ll  r_{\|}$ the correction $W^{(c)}_{(\xi)}$ is comparatively small (say in $L^2$) compared to $W_{(\xi)}$. This follows from the estimates obeyed by the first two terms on the left side of \eqref{e:W_Lp_bnd}.

Having established the important geometric properties of the intermittent jets, namely \eqref{e:Reynolds_cancellation} and \eqref{eq:useful:2}, we summarize the bounds obeyed by the intermittent jets, their incompressibility correctors, and their corresponding building blocks. By construction we have that the directions of oscillation for the functions defined in \eqref{eq:phi:psi:rotate} are orthogonal. Therefore, we may use Fubini to combine the estimates obeyed by $\psi_{(\xi)}$ (which is a 1D function) and $\phi_{(\xi)},\Phi_{(\xi)}$ (which are 2D functions), to obtain estimates for the 3D functions $W_{(\xi)}$ and $V_{(\xi)}$. We claim that for $N,M\geq 0$ and  $p\in [1,\infty]$ the following bounds hold: 
\begin{align}
\norm{\nabla^N\partial_t^M\psi_{(\xi)}}_{L^{p}}
&\les  r_{\|}^{\sfrac{1}{p}-\sfrac12}\left(\frac{ r_{\perp}\lambda}{ r_{\|}}\right)^{\!\! N} \left(\frac{ r_{\perp}\lambda\mu}{ r_{\|}}\right)^{\!\!M}
\label{e:phi_Lp_bnd} \\
\norm{\nabla^N\phi_{(\xi)}}_{L^{p}}+\norm{\nabla^N\Phi_{(\xi)}}_{L^{p}}
&\les  r_{\perp}^{\sfrac{2}{p}-1}\lambda^N 
\label{e:psi_Lp_bnd}\\
\norm{\nabla^N\partial_t^M W_{(\xi)}}_{L^{p}} 
+ \frac{r_\|}{r_\perp} \norm{\nabla^N\partial_t^M W_{(\xi)}^{(c)}}_{L^{p}}
+\lambda^2\norm{\nabla^N\partial_t^M V_{(\xi)}}_{L^{p}}
&\les  r_{\perp}^{\sfrac{2}{p}-1} r_{\|}^{\sfrac{1}{p}-\sfrac12}\lambda^N\left(\frac{ r_{\perp}\lambda\mu}{ r_{\|}}\right)^{\!\! M}
\label{e:W_Lp_bnd}
\end{align}
where implicit constants may depend on $p$, $N$ and $M$, but are independent of $\lambda$, $r_\perp$, $r_\|$, $\mu$. 
For $N,M=0$, the corresponding $L^p$ estimates follow directly from simple scaling arguments and the bounds implied by \eqref{eq:phi:psi:normalize}  on the unit periodic box. For the derivative estimates we have assumed that 
\[ r_{\|}^{-1}\ll  r_{\perp}^{-1}\ll \lambda\,\]
holds, in order to identify the largest frequency of oscillation.

\subsection{The perturbation}\label{sec:perturbation}

In this section we will construct the perturbation $w_{q+1}$ which defines $v_{q+1} = v_\ell + w_{q+1}$.

\subsubsection{Parameter choices}
\label{sec:NSE:parameters}
In Section~\ref{s:intermittent} above, we have worked with abstract parameters $\lambda, \mu, r_\perp$, and $r_\|$. They had to obey the heuristic bounds $\lambda^{-1} \ll r_\perp \ll r_\| \ll 1$, and $\lambda r_\perp$ had to be a natural number.  
It is convenient to now fix these parameters, all in terms of $\lambda_{q+1}$. At this stage in the proof not all these choices are motivated, and some of the choices only become justified based on arguments in Section~\ref{sec:NSE:stress}  below. However, in order to specify the bounds on the velocity increment $w_{q+1}$, and the new velocity $v_{q+1}$, it is convenient to specify these parameters at this stage already.

Recall that $\lambda_{q}$ is defined in \eqref{eq:parameters:def} as $2\pi a^{(b^q)}$.
We define
\begin{subequations}
\label{eq:NSE:parameters}
\begin{align}
r_\| &= \lambda_{q+1}^{-\sfrac 47} \label{eq:NSE:r:parallel:def}\\
r_\perp &=  r_\|^{-\sfrac 14} \lambda_{q+1}^{-1} (2\pi)^{-\sfrac 17} =  \lambda_{q+1}^{-\sfrac 67} (2\pi)^{-\sfrac 17} \label{eq:NSE:r:perp:def}\\
\mu &=\lambda_{q+1} r_\| r_\perp^{-1} = \lambda_{q+1}^{\sfrac 97 } (2\pi)^{\sfrac 17}\label{eq:NSE:mu:def}
\end{align}
\end{subequations}
In \eqref{eq:NSE:r:perp:def} we have introduced a strange power of $2\pi$ in order to ensure that 
\[
\lambda_{q+1}r_{\perp} = ((2\pi)^{-1} \lambda_{q+1})^{\sfrac 17} = a^{\sfrac{(b^{q+1})}{7}}\in\mathbb N \, ,
\] 
which ensures the correct periodicity of $W_{(\xi)}$, $V_{(\xi)}$, $\Phi_{(\xi)}$, $\phi_{(\xi)}$ and $\psi_{(\xi)}$. For this purpose, it  is sufficient to ensure that $b$ is a multiple of $7$. 
For the remainder of the Navier-Stokes section we use these  fixed choices of $\lambda, r_{\perp},r_{\|}$, and $\mu$,  for the short hand notation $W_{(\xi)}$, $V_{(\xi)}$, $\Phi_{(\xi)}$, $\phi_{(\xi)}$ and $\psi_{(\xi)}$ introduced in Section~\ref{s:intermittent}. 
We also recall that the parameter $\ell$ was chosen to obey \eqref{eq:ell:cond:NSE} above.

\subsubsection{Amplitudes}
Let $\chi$ be a smooth function such that
\[
\chi(z)=\begin{cases}1&\mbox{if} \quad 0\leq z \leq 1\\ z &\mbox{if}\quad z\geq 2\end{cases} \]
and with $z \leq 2 \chi(z) \leq 4 z$ for $z\in (1,2)$. For such a function $\chi$ we then define 
\[
\rho (x,t)= 4 c_R \delta_{q+1}\chi\left((c_R \delta_{q+1})^{-1}\abs{\RR_{\ell}(x,t)}\right)\,.
\]
This simplified definition of $\rho$ was introduced in~\cite{LuoTiti18}. The main properties of $\rho$ are that pointwise in $(x,t)$ we have
\begin{align}
\left| \frac{\RR_\ell(x,t)}{\rho(x,t)} \right| = \frac{1}{4} \frac{(c_R \delta_{q+1})^{-1} |\RR_\ell(x,t)|}{\chi((c_R \delta_{q+1})^{-1} |\RR_\ell (x,t)|)}  \leq \frac 12
\label{eq:rho:is:useful}
\end{align}
and moreover for any $p \in [1,\infty]$ we have
\begin{align}
\label{eq:rho:Lp}
\norm{\rho}_{L^p} \leq 16  \left(c_R (8 \pi^3)^{\sfrac 1p} \delta_{q+1}+\norm{\RR_{\ell}}_{L^p} \right)
\end{align}
Moreover, using the Sobolev embedding $W^{4,1} \subset C^0$, from \eqref{eq:Rq:L1:NSE}, standard mollification estimates and repeated applications of the chain rule (see e.g.~\cite[Proposition C.1]{BDLISZ15}), we obtain that 
\begin{align}
\label{eq:rho:CN}
\norm{\rho}_{C^N_{x,t}} \les  \delta_{q+1}  \ell^{-4}  \ell^{-  5N}  \les  \ell^{- 4- 5N} 
\end{align}
for any $N\geq 0$, where the implicit constant depends on $N$.

Next we define the amplitude functions
\begin{equation}
a_{(\xi)} (x,t) := a_{\xi,q+1}(x,t):= \;\rho(x,t)^{\sfrac 12} \; \gamma_{\xi}\left(\Id - \frac{\RR_{\ell}(x,t)}{ \rho(x,t)}\right) \, .
\label{eq:a:oxi:def}
\end{equation}
In view of \eqref{eq:rho:is:useful}, the matrix $\Id - \rho^{-1} \RR_\ell$ lies in the domain of definition of the functions $\gamma_\xi$ from Lemma~\ref{l:linear_algebra}.
We note that as a consequence of \eqref{e:Reynolds_cancellation},  \eqref{eq:a:oxi:def}, and of Lemma ~\ref{l:linear_algebra}, we have 
\begin{align}
\sum_{\xi  \in \Lambda} a_{(\xi)}^2 \fint_{\mathbb T^3}   W_{(\xi)}\otimes   W_{(\xi)}  dx=
 \rho\, \Id - \RR_{\ell}\,,
\label{e:WW_id}
\end{align}
which justifies the definition of the amplitude functions $a_{(\xi)}$.

The amplitudes $a_{(\xi)}$ have a good $L^2$ norm. Indeed, using \eqref{eq:rho:Lp}, the bound \eqref{eq:Rq:L1:NSE},  the fact that the mollifier has mass $1$, the definition of the constant $M$ in Lemma~\ref{l:linear_algebra}, and choosing the constant $C_\gamma$ as discussed in the first paragraph of Section~\ref{s:intermittent}, we obtain
\begin{align}
\norm{a_{(\xi)}}_{L^2}  \leq \norm{\rho}_{L^1}^{\sfrac 12} \norm{\gamma_{\xi}}_{C^0(B_{\sfrac 12}(\Id))} \leq   \frac{4 c_R^{\sfrac 12}  \left( (8\pi^3)  + 1\right)^{\sfrac 12} M}{8|\Lambda|   ((8\pi^3) + 1)^{\sfrac 12}}   \delta_{q+1}^{\sfrac 12} \leq \frac{\delta_{q+1}^{\sfrac 12}}{2|\Lambda| } 
\label{eq:a:xi:L2}
\end{align}
by choosing $c_R$ sufficiently small, in terms of the universal constants $M$ and $|\Lambda|$.
Lastly,  we note that since $\rho$ is bounded from below by $4 c_R \delta_{q+1}$, using the chain rule and standard mollification estimates, similarly to \eqref{eq:rho:CN}, we have 
\begin{align}
\norm{a_{(\xi)}}_{C^N_{x,t}}  \les \delta_{q+1}^{\sfrac 12} \ell^{-2} \ell^{- 5N} \les \ell^{-2-5N}
\label{eq:a:CN:NSE}
\end{align}
for all $N\geq 0$. 

\subsubsection{Principal part of the perturbation, incompressibility and temporal correctors}

The {\em principal part of $w_{q+1}$} is defined as
\begin{equation}
w^{(p)}_{q+1}:=\sum_{\xi\in \Lambda}a_{(\xi)} \;  W_{(\xi)}\,.
\label{eq:w:q+1:p:def}
\end{equation}
Note that in view of \eqref{eq:Messi} we have that the summands in \eqref{eq:w:q+1:p:def} have mutually disjoint supports. Also, by \eqref{e:WW_id} we have that the low-frequency part of $w_{q+1}^{(p)} \otimes w_{q+1}^{(p)}$ cancels the mollification stress $\RR_\ell$.  More precisely, we have
\begin{align}
w_{q+1}^{(p)} \otimes w_{q+1}^{(p)} + \RR_\ell =  \sum_{\xi \in \Lambda} a_{(\xi)}^2 \Proj_{\neq 0} \left( W_{(\xi)} \otimes W_{(\xi)} \right) + \rho \Id \, .
\label{eq:oscillation:identity:1}
\end{align}

In order to fix the fact that $w_{q+1}^{(p)}$ is not divergence free, we define an {\em incompressibility corrector} by
\begin{align}
w^{(c)}_{q+1}
&:=\sum_{\xi\in \Lambda} \curl\left(\nabla a_{(\xi)}\times V_{(\xi)}
\right)+\nabla a_{(\xi)} \times \curl V_{(\xi)} + a_{(\xi)} W_{(\xi)}^{(c)} 
\, .
\label{eq:w:q+1:c:def}
\end{align}
The above definition is motivated by a computation similar to \eqref{e:div_correct_id}, which guarantees that
\begin{align}\label{e:curl_form}
w_{q+1}^{(p)}+w_{q+1}^{(c)}= \sum_{\xi\in \Lambda} \curl\curl(a_{(\xi)}V_{(\xi)})\,,
\end{align}
and thus ensures
\begin{align*}
\div\left(w_{q+1}^{(p)}+w_{q+1}^{(c)}\right)\equiv 0\,.
\end{align*}

In addition to the incompressibility corrector $w_{q+1}^{(c)}$, we introduce a  {\em temporal corrector} $w_{q+1}^{(t)}$, which is defined by
\begin{equation}\label{e:temporal_corrector}
w^{(t)}_{q+1}:=-\frac{1}{ \mu} \sum_{\xi\in \Lambda }\mathbb P_{H}\mathbb P_{\neq 0}\left(a_{(\xi)}^2\phi^2_{(\xi)} \psi^2_{(\xi)}\xi \right)\, ,
\end{equation}
where as before $\Proj_H$ is the Helmholtz-Leray projector. By construction, $\div w_{q+1}^{(t)} = 0$, and the choice of $\mu$ ensures that this corrector is much smaller than the principal corrector $w_{q+1}^{(p)}$.  The purpose of introducing the temporal corrector only becomes apparent when considering the oscillation error in Section~\ref{sec:NSE:stress}. More precisely, we use the key identity \eqref{eq:useful:2} obeyed by the intermittent jets, and the fact that $\Id - \Proj_H = \nabla \Delta^{-1} \div$, to rewrite
\begin{align}
&\partial_t w_{q+1}^{(t)} + \sum_{\xi \in \Lambda}  \Proj_{\neq 0} \left( a_{(\xi)}^2 \div \left( W_{(\xi)} \otimes W_{(\xi)} \right) \right)  
\notag\\
&\quad = -\frac{1}{ \mu} \sum_{\xi\in \Lambda }\mathbb P_{H}\mathbb P_{\neq 0} \partial_t \left(a_{(\xi)}^2\phi^2_{(\xi)} \psi^2_{(\xi)}\xi \right)  +  \frac{1}{\mu} \sum_{\xi \in \Lambda}  \Proj_{\neq 0} \left( a_{(\xi)}^2 \partial_t \left(\phi^2_{(\xi)}  \psi^2_{(\xi)}\xi\right) \right) 
\notag\\
&\quad = \underbrace{(\Id - \Proj_H) \frac{1}{\mu} \sum_{\xi\in \Lambda} \Proj_{\neq 0}\partial_t \left(a_{(\xi)}^2\phi^2_{(\xi)} \psi^2_{(\xi)}\xi \right)}_{=:\nabla P}  -  \frac{1}{\mu} \sum_{\xi \in \Lambda}\Proj_{\neq 0}  \left( \partial_t a_{(\xi)}^2 \left(\phi^2_{(\xi)}  \psi^2_{(\xi)}\xi\right) \right) \, .
\label{eq:oscillation:identity:2}
\end{align}
Identity \eqref{eq:oscillation:identity:2} shows the essential role played by temporal oscillations in our construction.

\subsubsection{The velocity increment and verification of the inductive estimates}
\label{sec:velocity:size}

The total velocity increment $w_{q+1}$ is defined by
\begin{equation}\label{eq:w:q+1:def}
w_{q+1}:=w_{q+1}^{(p)}+w_{q+1}^{(c)}+w_{q+1}^{(t)}\,,
\end{equation}
and is by construction mean zero and divergence-free. 
The new velocity field $v_{q+1}$ is defined as
\begin{equation}\label{eq:v:q+1:def}
v_{q+1} = v_q + w_{q+1} \, .
\end{equation}
In this section we verify that the inductive estimates \eqref{eq:vq:L2:NSE}, \eqref{eq:vq:C1:NSE} hold with $q$ replaced by $q+1$, and that \eqref{eq:increment:L2} is satisfied.

In order to efficiently estimate the $L^2$ norm of the principal part of the perturbation $w_{q+1}^{(p)}$ we need to use the fact that the amplitudes $a_{(\xi)}$ oscillate at a much lower frequency than the intermittent jets $W_{(\xi)}$. For this purpose, we recall the $L^p$ de-correlation estimate~\cite[Lemma 3.6]{BV}:
\begin{lemma}
\label{lem:Lp:decorrelate}
Fix integers $N,\kappa \geq 1$ and let $\zeta > 1$ be such that
\begin{align}
\label{eq:Lp:independence:assume}
\frac{2\pi\sqrt{3} \zeta}{\kappa} \leq \frac{1}{3} \qquad \mbox{and} \qquad \zeta^{4} \frac{(2\pi\sqrt{3} \zeta)^N}{\kappa^N}  \leq 1\,.
\end{align}
Let $p \in \{1,2\}$, and let $f$ be a $\T^3$-periodic function such that there exists a constant $C_f>0$ such that
\begin{align*}
\|D^j f\|_{L^p} \leq C_f\zeta^j \,,
\end{align*}
holds for  all $0 \leq j \leq N+4$. 
In addition, let $g$ be a $(\sfrac{\T}{\kappa})^{3}$-periodic function. Then we have that 
\begin{align*}
 \|f g \|_{L^p} \lesssim C_f \|g\|_{L^p} \,,
\end{align*}
where the implicit constant is universal.
\end{lemma}
We wish to apply the above de-correlation Lemma~\ref{lem:Lp:decorrelate} in $L^2$ with $f = a_{(\xi)}$ and $g = W_{(\xi)}$, which is by construction $(\T/\kappa)^3$ periodic with $\kappa = \lambda_{q+1} r_\perp$. For this purpose, note that by \eqref{eq:a:xi:L2} and \eqref{eq:a:CN:NSE} we obtain that 
\[
\norm{D^j a_{(\xi)}}_{L^2} \leq \frac{1}{2|\Lambda|} \delta_{q+1}^{\sfrac 12} \ell^{-8j} \, ,
\]
and thus we can take $C_f = \frac{1}{2|\Lambda|} \delta_{q+1}^{\sfrac 12}$ and $\zeta = \ell^{-8}$. 
Since by \eqref{eq:ell:cond:NSE} we have $\zeta \leq \lambda_{q+1}^{16 \alpha}$, whereas by \eqref{eq:NSE:r:perp:def} we have that $\lambda_{q+1} r_\perp \approx \lambda_{q+1}^{\sfrac 17}$. Thus, since $\alpha$ is sufficiently small, condition \eqref{eq:Lp:independence:assume} is satisfied for any $N\geq 1$, and thus Lemma~\ref{lem:Lp:decorrelate} is applicable. Combining the resulting estimate with the normalization of $W_{(\xi)}$ as $\norm{W_{(\xi)}}_{L^2} = 1$,  we obtain
\begin{align}
\norm{w_{q+1}^{(p)}}_{L^2} \leq \sum_{\xi \in \Lambda}  \frac{2}{2|\Lambda|} \delta_{q+1}^{\sfrac 12} \norm{W_{(\xi)}}_{L^2}
\leq \frac{1}{2} \delta_{q+1}^{\sfrac 12}.
\label{eq:w:q+1:p:L2}
\end{align}
For the correctors $w_{q+1}^{(c)}$ and $w_{q+1}^{(t)}$, and for bounds on the other $L^p$ norms of $w_{q+1}^{(p)}$ we may afford slightly less precise bounds (which do not appeal to Lemma~\ref{lem:Lp:decorrelate}), which follow directly from \eqref{e:W_Lp_bnd}, \eqref{eq:a:CN:NSE}, and the parameter choices in \eqref{eq:NSE:parameters}:
\begin{subequations}
\label{eq:correctors:Lp}
\begin{align}
\norm{w_{q+1}^{(p)}}_{L^p} 
&\les \sum_{\xi \in \Lambda} \norm{a_{(\xi)}}_{C^0}  \norm{W_{(\xi)}}_{L^p} \notag\\
&\qquad \les \delta_{q+1}^{\sfrac 12}  \ell^{-2}   r_{\perp}^{\sfrac{2}{p}-1} r_{\|}^{\sfrac{1}{p}-\sfrac12}  \\
\norm{w_{q+1}^{(c)}}_{L^p} 
&\les \sum_{\xi \in \Lambda} \norm{a_{(\xi)}}_{C^0}  \norm{W_{(\xi)}^{(c)}}_{L^p}  + \norm{a_{(\xi)}}_{C^2}  \norm{V_{(\xi)}}_{W^{1,p}} \notag\\
&\qquad \les \delta_{q+1}^{\sfrac 12}  \ell^{-12} r_{\perp}^{\sfrac{2}{p}-1} r_{\|}^{\sfrac{1}{p}-\sfrac12}  \left( r_{\perp} r_{\|}^{-1} +   \lambda_{q+1}^{-1} \right) \les \delta_{q+1}^{\sfrac 12}  \ell^{-12} r_{\perp}^{\sfrac{2}{p}} r_{\|}^{\sfrac{1}{p}-\sfrac32}  \\
\norm{w_{q+1}^{(t)}}_{L^p} 
&\les \mu^{-1} \sum_{\xi \in \Lambda} \norm{a_{(\xi)}}_{C^0}^2  \norm{\phi_{(\xi)}}_{L^{2p}}^2 \norm{\psi_{(\xi)}}_{L^{2p}}^2 \notag\\
&\qquad \les \delta_{q+1}   \ell^{-4}    r_{\perp}^{\sfrac{2}{p}-1} r_{\|}^{\sfrac{1}{p}-\sfrac 12}  (\mu^{-1} r_\perp^{-1} r_\|^{-\sfrac 12} )\les  \delta_{q+1}   \ell^{-4} r_{\perp}^{\sfrac{2}{p}-1} r_{\|}^{\sfrac{1}{p}-2} \lambda_{q+1}^{-1}.
\end{align}
\end{subequations}
Combining \eqref{eq:w:q+1:p:L2} with the last two estimates of \eqref{eq:correctors:Lp} for $p=2$, and using the parameter choices \eqref{eq:NSE:parameters}, we obtain that for a sufficiently large constant $C$ (which is independent of $q$), we have
\begin{align}
\norm{w_{q+1}}_{L^2} 
&\leq \delta_{q+1}^{\sfrac 12} \left( \frac{1}{2}  + C \ell^{-12} r_\perp r_\|^{-1} + C \ell^{-4} r_\|^{-\sfrac 32} \lambda_{q+1}^{-1} \right) \notag\\
&\leq \delta_{q+1}^{\sfrac 12} \left( \frac{1}{2}  + C \lambda_{q+1}^{24 \alpha -\sfrac 27 }   + C \lambda_{q+1}^{8 \alpha - \sfrac 17}   \right)  \leq  \frac{3}{4} \delta_{q+1}^{\sfrac 12}\,.
\label{eq:L2:temp:1}
\end{align}
In the last inequality we have used that $\alpha$ is sufficiently small, that $\beta$ is sufficiently small in terms of $\alpha$, and that by letting $a$ be sufficiently large we have $C \lambda_{q+1}^{-\sfrac{1}{10}} \leq \sfrac{1}{8}$. Moreover, from \eqref{eq:V_q_ell_est:L2}  we have that 
\begin{align}
\norm{v_q - v_\ell}_{L^2} \les \lambda_{q+1}^{-\alpha} \leq \lambda_{q+1}^{-\sfrac{\alpha}{2}} \leq  \frac{1}{4}  \delta_{q+1}^{\sfrac 12}
\label{eq:L2:temp:2}
\end{align}
once we ensure that $\beta$ is sufficiently small, in terms of $\alpha$, and we take $a$ to be sufficiently large. Combining \eqref{eq:L2:temp:1} and \eqref{eq:L2:temp:2} we obtain that 
\begin{align*}
\norm{v_{q+1}- v_q}_{L^2} \leq \norm{w_{q+1}}_{L^2} + \norm{v_\ell - v_q}_{L^2} \leq   \delta_{q+1}^{\sfrac 12},
\end{align*}
and thus \eqref{eq:increment:L2} is satisfied. 

The bound \eqref{eq:vq:L2:NSE} at level $q+1$ follows from \eqref{eq:V_ell_est:L2}, \eqref{eq:L2:temp:1}, and the inequality $-\delta_q^{\sfrac 12} +   \delta_{q+1}^{\sfrac 12} \leq - \delta_{q+1}^{\sfrac 12}$. The latter simply requires that $a$ is taken sufficiently large.

In addition, taking either a spatial or a temporal derivative,   using the bounds \eqref{e:phi_Lp_bnd}--\eqref{e:W_Lp_bnd}, the estimate for $\ell^{-1}$ in \eqref{eq:ell:cond:NSE}, and the parameter choices \eqref{eq:NSE:parameters} we may prove   
\begin{subequations}
\label{eq:correctors:C1}
\begin{align}
\norm{w_{q+1}^{(p)}}_{C^1_{x,t}} &\les  \ell^{-7}   r_{\perp}^{ -1} r_{\|}^{ -\sfrac12} \lambda_{q+1} \left( 1+ \frac{r_\perp \mu}{r_\|} \right) \les  \ell^{-7}   r_{\perp}^{ -1} r_{\|}^{ -\sfrac12} \lambda_{q+1}^2    \\
\norm{w_{q+1}^{(c)}}_{C^1_{x,t}} &\les  \ell^{-17}   r_{\|}^{ -\sfrac32}  \lambda_{q+1} \left( 1+ \frac{r_\perp \mu}{r_\|} \right)  \les \ell^{-17}   r_{\|}^{ -\sfrac32}  \lambda_{q+1}^2  \\
\norm{w_{q+1}^{(t)}}_{C^1_{x,t}} &\les \ell^{-9} r_{\perp}^{ -1} r_{\|}^{ -2} \lambda_{q+1}^{-1} \lambda_{q+1}^{1+\alpha} \left( 1+ \frac{r_\perp \mu}{r_\|} \right) \les \ell^{-9} r_{\perp}^{ -1} r_{\|}^{ -2} \lambda_{q+1}^{1+\alpha} \, .
\end{align}
\end{subequations}
In the last inequality we have paid an extra power of $\lambda_{q+1}^{\alpha}$ because $\Proj_{H} \Proj_{\neq 0}$ is not a bounded operator on $C^0$.
Combining the $C^1_{x,t}$ estimate for $v_\ell$ with the above obtained bound \eqref{eq:V_ell_est:C1}, we obtain from \eqref{eq:v:q+1:def} that 
\begin{align*}
\norm{v_{q+1}}_{C^{1}_{x,t}} \les \norm{v_\ell}_{C^{1}_{x,t}} + \norm{w_{q+1}}_{C^{1}_{x,t}} 
&\les \lambda_{q+1}^{\alpha} + \lambda_{q+1}^{14 \alpha + \sfrac{22}{7}} +  \lambda_{q+1}^{34 \alpha + \sfrac{20}{7}}  +  \lambda_{q+1}^{19 \alpha + 3} \notag\\
&\leq \lambda_{q+1}^{2 \alpha} + \lambda_{q+1}^{35 \alpha + \sfrac{22}{7}} \leq \lambda_{q+1}^4
\end{align*}
once we take $\alpha$ to be sufficiently small. This proves \eqref{eq:vq:C1:NSE} at level $q+1$.

\subsection{Reynolds stress}
\label{sec:NSE:stress}

\subsubsection{Decomposition of the new Reynolds stress}

Recall that $v_{q+1} = w_{q+1} + v_{\ell}$, where $v_{\ell}$ is defined in Section~\ref{sec:mollify:NSE}. Subtracting from \eqref{e:Navier-Reynolds} at level $q+1$ the system \eqref{e:NSE_reynolds_ell}, we thus obtain (compare to \eqref{eq:R:q+1:A} in the case of Euler)
\begin{align}
\div \RR_{q+1} - \nabla  p_{q+1}
&= \underbrace{- \nu \Delta w_{q+1} +  \partial_t (w^{(p)}_{q+1} +w^{(c)}_{q+1} )  + \div( v_{\ell} \otimes w_{q+1} + w_{q+1} \otimes v_{\ell})}_{\div(R_{\rm linear}) + \nabla p_{\rm linear}} \notag\\
&\qquad + \underbrace{\div\left((w_{q+1}^{(c)}+w_{q+1}^{(t)}) \otimes w_{q+1}+ w_{q+1}^{(p)} \otimes (w_{q+1}^{(c)}+w_{q+1}^{(t)}) \right)}_{\div(R_{\rm corrector}) + \nabla p_{\rm corrector}} \notag\\
&\qquad + \underbrace{\div(w_{q+1}^{(p)} \otimes w_{q+1}^{(p)} + \RR_{\ell})+\partial_t w^{(t)}_{q+1}}_{\div(R_{\rm oscillation}) + \nabla p_{\rm oscillation}}   + \, \div (R_{\rm commutator}) - \nabla p_\ell   \,.
\label{eq:tilde:R:1}
\end{align}
Here $R_{\rm commutator}$ is as defined in \eqref{eq:NSE:R:commutator} and we have used the inverse divergence operator from \eqref{eq:RSZ} to define
\begin{align}
R_{\rm linear} 
&= - \nu \RSZ \Delta w_{q+1} + \RSZ \partial_t (w_{q+1}^{(p)} + w_{q+1}^{(c)}) + v_\ell \mathring \otimes w_{q+1} + w_{q+1} \mathring \otimes v_\ell 
\label{eq:NSE:R:linear} \\
R_{\rm corrector} 
&= (w_{q+1}^{(c)}+w_{q+1}^{(t)}) \mathring\otimes w_{q+1}+ w_{q+1}^{(p)} \mathring\otimes (w_{q+1}^{(c)}+w_{q+1}^{(t)}) 
\label{eq:NSE:R:corrector}
\end{align}
where $p_{\rm linear} = 2 v_\ell \cdot w_{q+1}$ and $p_{\rm corrector} =  |w_{q+1}|^2 - |w_{q+1}^{(p)}|^2$. The remaining stress $R_{\rm oscillation}$ and corresponding pressure $p_{\rm oscillation}$ are defined as follows.
By \eqref{eq:oscillation:identity:1} and \eqref{eq:oscillation:identity:2} we have that 
\begin{align*}
 & \div(w_{q+1}^{(p)} \otimes w_{q+1}^{(p)} + \RR_{\ell})+\partial_t w^{(t)}_{q+1}
 \notag\\
 &\qquad =  \sum_{\xi \in \Lambda} \div\left( a_{(\xi)}^2 \Proj_{\neq 0} \left( W_{(\xi)} \otimes W_{(\xi)} \right)  \right) + \nabla \rho+\partial_t w^{(t)}_{q+1} \notag\\
 &\qquad =  \sum_{\xi \in \Lambda} \Proj_{\neq 0} \left( \nabla a_{(\xi)}^2 \Proj_{\neq 0} \left( W_{(\xi)} \otimes W_{(\xi)} \right)  \right) + \nabla \rho +   \sum_{\xi \in \Lambda} \Proj_{\neq 0} \left( a_{(\xi)}^2 \div \left( W_{(\xi)} \otimes W_{(\xi)} \right)  \right) +\partial_t w^{(t)}_{q+1} \notag\\
 &\qquad =  \sum_{\xi \in \Lambda} \Proj_{\neq 0} \left( \nabla a_{(\xi)}^2 \Proj_{\neq 0} \left( W_{(\xi)} \otimes W_{(\xi)} \right)  \right) + \nabla \rho + \nabla P -  \frac{1}{\mu} \sum_{\xi \in \Lambda}\Proj_{\neq 0}  \left( \partial_t a_{(\xi)}^2 \left(\phi^2_{(\xi)}  \psi^2_{(\xi)}\xi\right) \right) \, ,
\end{align*}
where the pressure term $P$ is as defined by the first term on the right side of \eqref{eq:oscillation:identity:2}. Therefore, we define $p_{\rm oscillation} = \rho + P$, and let 
\begin{align}
R_{\rm oscillation} 
&=  \sum_{\xi \in \Lambda} \RSZ  \left( \nabla a_{(\xi)}^2 \Proj_{\neq 0} \left( W_{(\xi)} \otimes W_{(\xi)} \right)  \right)   -  \frac{1}{\mu} \sum_{\xi \in \Lambda} \RSZ  \left( \partial_t a_{(\xi)}^2 \left(\phi^2_{(\xi)}  \psi^2_{(\xi)}\xi\right) \right) \, .
\label{eq:NSE:R:oscillation}
\end{align}
(We recall there that $\RSZ = \RSZ \Proj_{\neq 0}$.)
Now all the terms in \eqref{eq:tilde:R:1} are well defined. We have $p_{q+1} = p_\ell - p_{\rm oscillation}  - p_{\rm corrector} - p_{\rm linear} $ and  
\begin{align}
\RR_{q+1} = R_{\rm linear} + R_{\rm corrector} + R_{\rm oscillation} + R_{\rm commutator} \,
\label{eq:R:q+1:BB}
\end{align}
which are stresses defined in \eqref{eq:NSE:R:linear}, \eqref{eq:NSE:R:corrector}, \eqref{eq:NSE:R:oscillation}, and respectively \eqref{eq:NSE:R:commutator}.

\subsubsection{Estimates for the new Reynolds stress}

We need to estimate the new stress $\RR_{q+1}$ in $L^1$. However, Calder\'on-Zygmund operators such as $\nabla \RSZ$ just fail to be bounded on $L^1$. For this purpose, we introduce an integrability parameter 
\begin{align}
p \in (1,2] \, \qquad \mbox{such that} \qquad p-1 \ll 1.
\end{align}
We may now use that  Calder\'on-Zygmund operators are bounded on $L^p$. Recalling \eqref{eq:NSE:parameters}, we fix $p$ to obey
\begin{align}
r_{\perp}^{\sfrac{2(1-p)}{p}} r_{\|}^{\sfrac{(1-p)}{p}} \leq (2\pi)^{\sfrac{1}{7}} \lambda_{q+1}^{\sfrac{16 (p-1)}{(7p)}} \leq \lambda_{q+1}^{\alpha} 
\label{eq:NSE:p:choice}
\end{align}
where $0 < \alpha \ll 1$ is as in Remark~\ref{rem:b:beta:NSE}. For instance, we may take $p = \sfrac{32}{(32-7 \alpha)} > 1$.

{\bf Linear error. \,}
Recalling identity~\eqref{e:curl_form}, using that $\nu \in (0,1]$,   that $\RSZ \curl$ is a bounded operator on $L^p$, and appealing to the Sobolev embedding $H^3 \subset C^0$, the linear stress defined in \eqref{eq:NSE:R:linear}  obeys
\begin{align*}
\norm{R_{\rm linear}}_{L^p}  
&\les \nu \norm{\RSZ \Delta w_{q+1} }_{L^p} + \norm{ v_{\ell} \mathring \otimes w_{q+1} + w_{q+1} \mathring \otimes v_{\ell} }_{L^p} +\norm{ \RSZ \partial_t (w^{(p)}_{q+1} +w^{(c)}_{q+1} )}_{L^p} 
\notag\\
& \les
\norm{\nabla w_{q+1}}_{L^p} + \norm{v_{\ell}}_{L^{\infty}}\norm{w_{q+1}}_{L^{p}}  +\sum_{\xi\in \Lambda } \norm{\partial_t  \curl(a_{(\xi)}V_{(\xi)})}_{L^p}
\notag\\
&\les
(1+\norm{v_{\ell}}_{C^1}) \sum_{\xi\in \Lambda }\norm{ a_{(\xi)}}_{C^1} \norm{W_{(\xi)}}_{W^{1,p}} \notag\\
&\qquad + \sum_{\xi\in \Lambda }  \left(\norm{a_{(\xi)}}_{C^1} \norm{\partial_t V_{(\xi)}}_{W^{1,p}}+ \norm{\partial_t a_{(\xi)}}_{C^1}  \norm{V_{(\xi)}}_{W^{1,p}}\right) \, .
\end{align*}
By appealing to \eqref{eq:V_ell_est:C1}, \eqref{e:W_Lp_bnd}, \eqref{eq:a:CN:NSE}, and \eqref{eq:NSE:p:choice} we obtain from the above estimate that
\begin{align}
\norm{R_{\rm linear}}_{L^p} 
&\les  \ell^{-8}   r_{\perp}^{\sfrac{2}{p}-1} r_{\|}^{\sfrac{1}{p}-\sfrac12}\lambda_{q+1} +  \ell^{-7}  r_{\perp}^{\sfrac{2}{p} } r_{\|}^{\sfrac{1}{p}-\sfrac32} \mu +  \ell^{-12}  r_{\perp}^{\sfrac{2}{p}-1} r_{\|}^{\sfrac{1}{p}-\sfrac12}\lambda_{q+1}^{-1}\notag\\
&\les  \ell^{-12}   r_{\perp}^{\sfrac{2(1-p)}{p}} r_{\|}^{\sfrac{(1-p)}{p}}  \left(r_\perp r_\|^{\sfrac 12} \lambda_{q+1} + r_\perp^2 r_\|^{-\sfrac 12} \mu  \right) \notag\\
&\les  \ell^{-12} \lambda_{q+1}^{\alpha} \left(r_\perp r_\|^{\sfrac 12} \lambda_{q+1} + r_\perp^2 r_\|^{-\sfrac 12} \mu  \right)   \,.
\label{eq:NSE:R:linear:1}
\end{align}

{\bf Corrector error. \,}
Recall the definitions \eqref{eq:w:q+1:p:def}, \eqref{eq:w:q+1:c:def}, \eqref{e:temporal_corrector}, and the estimates \eqref{eq:correctors:Lp}. From the parameter choices in \eqref{eq:NSE:parameters}, we see that the worse $L^{2p}$ estimate is the one for for the principal corrector. 

Therefore, using that $\delta_{q+1} \leq 1$, from the parameter inequalities  \eqref{eq:correctors:Lp} is used, and the choice of $p$ in \eqref{eq:NSE:p:choice}, it follows that the corrector stress defined in \eqref{eq:NSE:R:corrector} obeys 
\begin{align}
\norm{ R_{\rm corrector}}_{L^p} 
&\leq   \norm{ w_{q+1}^{(c)}+w_{q+1}^{(t)} }_{L^{2p}}   \norm{ w_{q+1}}_{L^{2p}} + \norm{w_{q+1}^{(p)}  }_{L^{2p}}   \norm{ w_{q+1}^{(c)}+w_{q+1}^{(t)} }_{L^{2p}}
\notag \\
&\les \left( \ell^{-12} r_\perp^{\sfrac 1p} r_\|^{\sfrac{1}{2p} - \sfrac 32} +  \ell^{-4} \lambda_{q+1}^{-1} r_\perp^{\sfrac{1}{p}-1} r_\|^{\sfrac{1}{2p} - 2} \right)\ell^{-2}   r_{\perp}^{\sfrac{1}{p}-1} r_{\|}^{\sfrac{1}{2p}-\sfrac12}
\notag\\
&\les \ell^{-14} \lambda_{q+1}^{\alpha} \left(r_\perp r_\|^{-1}    + \lambda_{q+1}^{-1}  r_\|^{-\sfrac{3}{2}}   \right)     
\,.
\label{eq:NSE:R:linear:2}
\end{align}

{\bf Oscillation error. \,}
The oscillation error defined in \eqref{eq:NSE:R:oscillation} has two terms, let us call $R_{\rm oscillation}^{(x)}$ the one when the spatial gradient lands on $a_{(\xi)}$ and let us denote by $R_{\rm oscillation}^{(t)}$ the one when temporal derivatives land on $a_{(\xi)}$. 

Let us treat $R_{\rm oscillation}^{(x)}$ first. Note that by definition $W_{(\xi)}$ is $(\sfrac{\T}{r_\perp \lambda_{q+1}})^3$ periodic. Therefore, so is $W_{(\xi)} \otimes W_{(\xi)}$ and we obtain that the  minimal active frequency in $\Proj_{\neq 0} (W_{(\xi)} \otimes W_{(\xi)})$ is given by $r_\perp \lambda_{q+1}$. Equivalently, we have 
\[
\Proj_{\neq 0} (W_{(\xi)} \otimes W_{(\xi)}) = \Proj_{\geq \sfrac{r_\perp \lambda_{q+1}}{2}} (W_{(\xi)} \otimes W_{(\xi)}).
\]
On the other hand, the amplitude term $\nabla a_{(\xi)}^2$ oscillates at a much lower frequency (proportional to  $\ell^{-5}$), so that we expect the inverse divergence operator $\RSZ$ to gain a factor of $r_\perp \lambda_{q+1}$. To make this intuition precise, let us recall \cite[Lemma B.1]{BV}.
\begin{lemma}
\label{lem:comm:1}
Fix parameters $1 \leq \zeta < \kappa$, $p \in (1,2]$, and assume there exists a $N \in \N$ such that 
\[
 \zeta^N \leq \kappa^{N-2} \,.
\]
Let $a \in C^N(\T^3)$ be such that there exists $C_a >0$ with
\begin{align}
\norm{D^j a}_{C^0} \leq C_a \zeta^j
\label{eq:Merlot:1}
\end{align}
for all $0 \leq j \leq N$. Assume furthermore that $f\in L^p(\T^3)$ is such that $\int_{\T^3} a(x) \Proj_{\geq \kappa} f(x) dx = 0$.
Then we have
\begin{align}
\norm{ |\nabla|^{-1} (a \; \Proj_{\geq \kappa} f)}_{L^p} 
\les C_a   \frac{\norm{f}_{L^p}}{\kappa}
\label{eq:Merlot:2}
\end{align}
where the implicit constant depends only on $p$ and $N$.
\end{lemma}
From \eqref{eq:a:CN:NSE} we see that $\nabla a_{(\xi)}^2$ obeys \eqref{eq:Merlot:1} with $C_a = \ell^{-9}$, and $\zeta = \ell^{-5}$.  Since $\ell^{-5} \leq \lambda_{q+1}^{10\alpha} \ll \lambda_{q+1}^{\sfrac 17} \approx r_\perp \lambda_{q+1}$, we are justified to use Lemma~\ref{lem:comm:1} (with any $N\geq 3$), combined with the bound \eqref{e:W_Lp_bnd}, to conclude that 
\begin{align}
\norm{R_{\rm oscillation}^{(x)}}_{L^p} 
&\leq \sum_{\xi \in \Lambda} \norm{\RSZ  \left( \nabla a_{(\xi)}^2 \Proj_{\geq \sfrac{r_\perp \lambda_{q+1}}{2}} \left( W_{(\xi)} \otimes W_{(\xi)} \right)  \right)}_{L^p}
\notag\\
&\les \ell^{-9}  \frac{\norm{W_{(\xi)} \otimes W_{(\xi)}}_{L^p}}{r_\perp \lambda_{q+1}}
 \les \ell^{-9}   \frac{\norm{W_{(\xi)}}_{L^{2p}}^2}{r_\perp \lambda_{q+1}}
\notag\\
&\les \ell^{-9}   r_\perp^{\sfrac{2(1-p)}{p}} r_{\|}^{\sfrac{(1-p)}{p}} (r_\perp^{-1} \lambda_{q+1}^{-1}) \notag \\
&\les \ell^{-9}   \lambda_{q+1}^\alpha (r_\perp^{-1} \lambda_{q+1}^{-1})   
\, .
\label{eq:NSE:R:linear:3}
\end{align}
On the other hand, for the second term on the right side of \eqref{eq:NSE:R:oscillation} we just use that $\RSZ$ is a bounded operator on $L^p$ (we don't use that is a smoothing operator), Fubini's theorem to integrate along the orthogonal directions of $\phi_{(\xi)}$ and $\psi_{(\xi)}$ via \eqref{e:phi_Lp_bnd} and \eqref{e:psi_Lp_bnd}, and the bound \eqref{eq:a:CN:NSE} for the amplitude functions, to conclude   
\begin{align}
\norm{ R_{\rm oscillation}^{(t)}}_{L^p} 
&\les \mu^{-1} \sum_{\xi \in \Lambda} \norm{\partial_t a_{(\xi)}^2}_{C^0} \norm{\phi_{(\xi)}}_{L^{2p}}^2 \norm{\psi_{(\xi)}}_{L^{2p}}^2 
\notag\\
&\les \mu^{-1}  \ell^{-9}   r_{\perp}^{\sfrac{2(1-p)}{p}} r_{\|}^{\sfrac{(1-p)}{p}} \notag \\
&\les \mu^{-1}  \ell^{-9}  \lambda_{q+1}^\alpha \, .
\label{eq:NSE:R:linear:4}
\end{align}

\subsubsection{Proof of \eqref{eq:Rq:L1:NSE} at level $q+1$}
\label{eq:NSE:end:of:proof}

The stress $\RR_{q+1}$ as defined in \eqref{eq:R:q+1:BB} above, may now be estimated by combining \eqref{eq:Rc:bound:L1}, \eqref{eq:NSE:R:linear:1}, \eqref{eq:NSE:R:linear:2}, \eqref{eq:NSE:R:linear:3}, \eqref{eq:NSE:R:linear:4}, the bound for $\ell$ given by \eqref{eq:ell:cond:NSE}, the parameter choices \eqref{eq:NSE:parameters}, and choosing $\alpha$ and $\beta b$ to be sufficiently small. Since $L^p(\T^3) \subset L^1(\T^3)$, we obtain
\begin{align*} 
\norm{\RR_{q+1}}_{L^1} 
&\les \norm{R_{\rm linear}}_{L^p} + \norm{R_{\rm corrector}}_{L^p} + \norm{R_{\rm oscillation}}_{L^p} + \norm{R_{\rm commutator}}_{L^1} \notag\\
&\les \lambda_{q+1}^{29 \alpha - \frac 17} +  \lambda_{q+1}^{-\alpha}
\leq \lambda_{q+1}^{30 \alpha - \frac 17} +  \lambda_{q+1}^{-\sfrac{\alpha}{2}} \leq \delta_{q+2} \, .
\end{align*}
In the second to last above inequality we have used a power of $\lambda_{q+1}^\alpha$ to absorb the implicit $q$ independent constant, while in the last inequality we have used that by Remark~\ref{rem:b:beta:NSE}, we have $\alpha \leq \frac{1}{240}$, and $ \beta b < \alpha$. This proves that \eqref{eq:Rq:L1:NSE} holds at level $q+1$, thereby concluding the proof of Proposition~\ref{prop:iteration:NSE}.

\section{Open problems}\label{sec:open}

Prior to the work \cite{BV}, a proof of non-uniqueness of weak solutions to the Navier-Stokes equations via convex integration was widely believed to be infeasible via the techniques of convex integration. The use of intermittency in the context of convex integration, significantly widens the scope of applicability of convex integration to nonlinear PDE~\cite{ModenaSZ17,ModenaSZ18,Luo18,LuoTiti18,BCV18,Dai18,CheskidovLuo19}. Previously, in the language of Gromov, convex integration achieved flexibility of PDE via low regularity. The paper \cite{BV} demonstrates that flexibility may also be attained via low integrability. This leads us to the following rather open ended problem.

\begin{open_problem}
Given the expanded applicability of intermittent convex integration techniques, in what new  contexts can one apply convex integration?\label{p:Buffon}  
\end{open_problem}

It is interesting to observe that while Onsager's conjecture was originally stated in the context of the 3D Euler equations, the arguments used to prove Part~\eqref{Onsager:a} are not dimension dependent. Indeed, the arguments of \cite{ConstantinETiti94,CCFS08} apply equally for any dimension $d\geq 2$. It is then natural to ask whether Onsager's conjecture holds in 2D. Currently, the best known result in 2D is that there exist non-conservative weak solutions to the Euler equations lying in any H\"older space with H\"older exponent less that $\sfrac15$ \cite{ChDLSz12,Choffrut13,BSV16,Novack18}.

\begin{open_problem}
Can one prove Onsager's conjecture in 2D, i.e.~for every $\beta<\sfrac13$ can one construct non-conservative weak solutions to the 2D Euler equations lying in the H\"older space $C^\beta_{x,t}$.
\label{p:RobertoCarlos}
\end{open_problem}
 
Motivated in part by the successful application of entropy conditions to conservations laws, the \emph{local energy inequality} (see~\eqref{eq:KHM}) 
\begin{align}
\partial_t\left(\frac{\abs{v}^2}{2}\right)+\div\left(\left(\frac{\abs{v}^2}{2}+p\right)v\right) = - D(v) \leq 0
\label{eq:suitable:weak:Euler}
\end{align}
which is meant in a distributional sense,
has been proposed as a possible admissibility criterium to recover uniqueness for the Euler equations. In \cite{DeLellisSzekelyhidi10}, De Lellis and Sz\'ekelyhidi Jr., demonstrated that unlike the entropy condition for the Burgers equation, the local energy inequality does not uniquely select a solution to the initial value problem. Nevertheless, owing to the fact that the local energy inequality is strictly stronger than the regular energy inequality, it remains of interest to determine whether for every H\"older exponent $<\sfrac13$, there exists non-conservative weak solutions to the Euler equations satisfying the local energy inequality. This would in effect verify a stronger version of Onsager's conjecture. We refer the reader to \cite{Isett17b} for recent progress in this direction.

\begin{open_problem}
For any $\beta<\sfrac13$, can one construct non-conservative weak solutions  $v\in C^{\beta}_{x,t}$ satisfying the local energy inequality \eqref{eq:suitable:weak:Euler} ?
\label{p:Maldini}
\end{open_problem}

It also a natural question to ask, what happens at the critical exponent $\sfrac13$? As was mentioned in Section~\ref{sec:Euler:C1/3}, in  this direction,  Isett proved in \cite{Isett17} the existence of non-conservative weak solutions lying in the intersection of all H\"older spaces $C^\beta_{x,t}$ for $\beta<\sfrac13$. More specifically, for any $B>\sfrac43$, a non-conservative weak solution $v$ can be constructed satisfying
\[\sup_{|\ell|>0} \abs{ \delta v(x,t;\ell ) } \abs{\ell}^{-\sfrac13+B\sqrt{\frac{\log\log \abs{\ell}^{-1}}{\log\abs{\ell}^{-1}}}} < \infty\,.\]
\begin{open_problem}
Does there exist non-conservative weak solutions to the Euler equations lying in the critical space $L^3_t B^{1/3}_{3,\infty,x}$ identified in~\cite{CCFS08} ?
\label{p:Cafu}
\end{open_problem}

One may also consider Onsager's conjecture in the context of $L^2$ based Sobolev spaces. In \cite{SulemFrisch75} (see also \cite{CCFS08}), the authors prove that kinetic energy is conserved for any weak solution $v\in H^{\sfrac56}$. It is an open problem to determine whether or not this result is sharp. Such a result has important implications for the physical theory of intermittency. As mentioned in Section~\ref{sec:intermittent}, as a consequence of intermittency, one expects the second order structure function to satisfy $\zeta_2>2/3$. It would be interesting to determine the largest possible deviation from the Kolmogorov prediction which is sustained by the Euler equations. In terms of weak solutions realizations, this corresponds to the existence of non-conservative weak solutions $u\in L^{\infty}_t H_x^{\beta}$ for $\beta>\sfrac13$.

\begin{open_problem}
For every $\beta<\sfrac56$, does there exist non-conservative weak solution to the Euler equations that lie in $H^\beta$?
\label{p:Beckenbauer}
\end{open_problem}

As was mentioned Section \ref{ss:euler_results}, in the context of the Euler equations, the Onsager exponent $\alpha_O=\sfrac13$ is not the sole interesting threshold exponent. In the context of the 2D Euler equations, a classical problem tracing back to the work of Yudovich is to determine where non-uniqueness occurs in the class of weak solutions with vorticity bounded in $L^p$ for $p<\infty$ \cite{Yudovich63}. If such non-uniqueness holds for all $p<\infty$ this would demonstrate that the uniqueness threshold exponent is $\alpha_U=1$ and would constitute a proof that the threshold of regularity for which kinetic energy is conserved does not coincide with the threshold of regularity for which uniqueness holds. We refer the reader to a recent paper of Vishik in this direction \cite{Vishik18}.

\begin{open_problem}
Given $p<\infty$, can one demonstrate non-uniqueness in the class of weak solutions to the 2D Euler equations with vorticity bounded in $L^p$? More generally, for either the 2D or 3D Euler equations, can one demonstrate that $\alpha_O\neq \alpha_U$?
\label{p:CRonaldo}
\end{open_problem}

In view of the partial regularity results of Leray~\cite{Leray34}, Scheffer~\cite{Scheffer76}, the Navier-Stokes inequality results~\cite{Scheffer85,Scheffer87,Ozanski18}, and the recent work \cite{BCV18} it is natural to investigate the limits of partial regularity in the context of the Navier-Stokes equations. This leads us to following open problem.

\begin{open_problem}
For any $0<d_{\Sigma_T}<1$, can one demonstrate the existence of weak solutions to the Navier-Stokes equations with non-empty singular set $\Sigma_T$ with Hausdorff dimension less than $d_{\Sigma_T}$?
\label{p:Zidanne}
\end{open_problem}

We note that for weak solutions (either of the type described in Definition \ref{d:weak_sol} or of Leray-Hopf type in Definitions \ref{d:leray-hopf} and \ref{d:loc_leray-hopf}), regularity and uniqueness is implied if one of the Lady{\v{z}}enskaja-Prodi-Serrin conditions is satisfied, i.e.\ the solution is bounded in a space $L^p_t L^q_x$ for  $\sfrac 2p + \sfrac 3q \leq 1$ \cite{KiselevLadyzhenskaya57,Prodi59,Serrin62,EscauriazaSerginSverak03,FabesJonesRiviere72,FurioliLemarieRieussetTerraneo00,LionsMasmoudi01,LemarieRieusset02,Kukavica06b}. The converse statement is however open:

\begin{open_problem}
For every $p,q$ satisfying $\sfrac 2p + \sfrac 3q > 1$, can one demonstrate the non-uniqueness of weak solutions (in the sense of Definition \ref{d:weak_sol}) bounded in $L_t^pL_x^q$.
\label{p:Ronaldinho}
\end{open_problem}

Perhaps the most fundamental question regarding weak solutions in fluid dynamics is to verify the  famous conjecture of Lady{\v{z}}enskaja \cite{Lady67} regarding the non-uniqueness of Leray-Hopf solutions to the Navier-Stokes equations. As already mentioned, \v{S}ver\'ak and Jia proved conditional non-uniqueness, assuming a spectral assumption. Although there is compelling numerical evidence \cite{GuillodSverak17} to support this assumption, the assumption appears to be remarkably difficult to verify analytically. While a non-uniqueness result involving Leray-Hopf solutions satisfying the local energy inequality (Definition \ref{d:loc_leray-hopf}) appears to be out of reach of methods involving convex integration, a proof of non-uniqueness satisfying the regular energy inequality (Definition \ref{d:leray-hopf}) via convex integration remains plausible.

\begin{open_problem}
Are the solutions to the initial value problem for the Navier-Stokes equation unique in the class of Leray-Hopf (in the sense of either Definition \ref{d:leray-hopf} or \ref{d:loc_leray-hopf}) weak solutions?
\label{p:Messi}
\end{open_problem}

The last problem is motivated by a desire to better understand the mechanisms of anomalous dissipation. In spite of the extensive experimental evidence supporting \eqref{eq:anomaly}, and the fact that $\eps > 0 $ is the fundamental ansatz of the Kolmogorov and Onsager theories of fully developed hydrodynamic turbulence,  to date this has not been proven in a mathematically rigorous context.\footnote{Compare this question to Theorem~\ref{thm:NSE:Euler} and Remark~\ref{rem:vanishing:viscosity}.} 
\begin{open_problem}
Can one produce a sequence of Leray-Hopf solutions $\{v^\nu\}_{\nu>0}$ of the forced Navier-Stokes system \eqref{eq:NSE}, with smooth forcing acting only at large ($\nu$-independent) scales, for which the corresponding sequence $\{ \eps^\nu \}_{\nu >0}$ (defined in \eqref{eq:epsilon:nu} by a long-time and space average), has a finite and non-zero limit? Moreover, can one show that $v^\nu$ converges to a dissipative weak solution of the Euler equations?
\label{p:Ronaldo}
\end{open_problem}

{\small

}

\end{document}